\documentclass[reqno]{amsart}
\usepackage{amsmath,amsrefs,amssymb,color}

\numberwithin{equation}{section}

\DeclareMathOperator{\divergence}{div}
\DeclareMathOperator{\bigO}{O}
\DeclareMathOperator{\smallo}{o}
\DeclareMathOperator{\Scal}{Scal}
\DeclareMathOperator{\Weyl}{Weyl}
\DeclareMathOperator{\Vol}{Vol}
\DeclareMathOperator{\Isom}{Isom}
\newcommand{\id}{id}
\newcommand{\loc}{loc}
\newcommand{\sym}{sym}
\newcommand{\eps}{\varepsilon}
\newcommand{\R}{\mathbb{R}}
\renewcommand{\S}{\mathbb{S}}
\newcommand{\N}{\mathbb{N}}

\newcommand{\<}{\left<}
\renewcommand{\>}{\right>}
\renewcommand{\[}{\left[}
\renewcommand{\]}{\right]}
\renewcommand{\(}{\left(}
\renewcommand{\)}{\right)}

\newtheorem{theorem}{Theorem}[section]
\newtheorem{corollary}[theorem]{Corollary}
\newtheorem{proposition}[theorem]{Proposition}
\newtheorem{lemma}[theorem]{Lemma}
\newtheorem{step}{Step}

\begin{document}

\title[Stability and instability results for sign-changing solutions]{Stability and instability results for sign-changing solutions to second-order\\critical elliptic equations}

\author{Bruno Premoselli}

\address{Bruno Premoselli, Universit\'e Libre de Bruxelles, Service d'analyse, CP 218, Boulevard du Triomphe, B-1050 Bruxelles, Belgique.}
\email{bruno.premoselli@ulb.be}

\author{J\'er\^ome V\'etois}

\address{J\'er\^ome V\'etois, Department of Mathematics and Statistics, McGill University, 805 Sherbrooke Street West, Montreal, Quebec H3A 0B9, Canada}
\email{jerome.vetois@mcgill.ca}

\date{February 28, 2022}

\thanks{Published in {\it Journal de Math\'ematiques Pures et Appliqu\'ees} {\bf167} (2022), 257–-293.}

\thanks{The first author was supported by the FNRS CdR grant J.0135.19, the Fonds Th\'elam and an ARC Avanc\'e 2020 grant. The second author was supported by the Discovery Grant RGPIN-2016-04195 from the Natural
Sciences and Engineering Research Council of Canada.}

\begin{abstract}
On a smooth, closed Riemannian manifold $\(M,g\)$ of dimension $n\ge3$, we consider the stationary Schr\"odinger equation $\Delta_gu+h_0u=\left|u\right|^{2^*-2}u$, where $\Delta_g:=-\text{div}_g\nabla$, $h_0\in C^1\(M\)$ and $2^*:=\frac{2n}{n-2}$. We prove that, up to perturbations of the potential function $h_0$ in $C^1\(M\)$, the sets of sign-changing solutions that are bounded in $H^1\(M\)$ are precompact in the $C^2$ topology. We obtain this result under the assumptions that $\(M,g\)$ is locally conformally flat, $n\ge7$ and $h_0\ne\frac{n-2}{4\(n-1\)}\Scal_g$ at all points in $M$, where $\Scal_g$ is the scalar curvature of the manifold. We then provide counterexamples in every dimension $n\ge3$ showing the optimality of these assumptions.
\end{abstract}

\maketitle

\section{Introduction}\label{Intro}

Let $\(M,g\)$ be a smooth, closed (i.e. compact and without boundary) Riemannian manifold of dimension $n\ge3$. We consider the stationary Schr\"odinger equation
\begin{equation}\label{IntroEq1}
\Delta_gu+h_0u=\left|u\right|^{2^*-2}u\quad\text{in }M,
\end{equation}
where $\Delta_g:=-\text{div}_g\nabla$ is the Laplace--Beltrami operator on functions, $h_0\in C^1\(M\)$ and $2^* =\frac{2n}{n-2}$ is the critical exponent for the embeddings of the Sobolev space $H^1\(M\)$ into Lebesgue's spaces. In this article, we are interested in a question of stability of sign-changing solutions to \eqref{IntroEq1} under perturbations of the potential function $h_0$. More precisely, we are interested in the question of existence of families of functions $\(h_\eps\)_{0<\eps\ll1}$ in $C^1\(M\)$ and $\(u_\eps\)_{0<\eps\ll1}$ in $C^2\(M\)$ such that $h_\eps\to h_0$ in $C^1\(M\)$ as $\eps\to0$, $\(u_\eps\)_{\eps}$ is bounded in $H^1\(M\)$ and blows up at least at some point $\xi_0\in M$ as $\eps\to0$ (i.e. there exists a family of points $\(\xi_\eps\)_{0<\eps\ll1}$ such that $\xi_\eps\to\xi_0$ and $\left|u_\eps\(\xi_\eps\)\right|\to\infty$ as $\eps\to0$) and for each $\eps>0$, $u_\eps$ is a sign-changing solution to the equation
\begin{equation}\label{IntroEq2}
\Delta_gu_\eps+h_\eps u_\eps=\left|u_\eps\right|^{2^*-2}u_\eps\quad\text{in }M.
\end{equation}
By ``sign-changing solutions'' we mean solutions that are allowed to change sign. Stability questions of this type have been extensively studied in recent years, and most of the focus has been on positive solutions. A possible reference in book form on stability questions for this type of equations is by Hebey \cite{He2}. In the case of positive solutions to the Yamabe equation (i.e. when $u_\eps>0$ and $h_\eps=h_0=c_n\Scal_g$ in $M$ for all $\eps$, where $c_n:=\frac{n-2}{4\(n-1\)}$ and $\Scal_g$ is the scalar curvature of the manifold), we refer to the work by Brendle \cite{Br}, Brendle and Marques \cite{BM1}, Druet \cite{Dr2}, Khuri, Marques and Schoen \cite{KMS}, Li and Zhang \cites{LZha1,LZha2}, Li and Zhu \cite{LZhu}, Marques \cite{M} and Schoen \cites{Sc2,Sc3} (see also the survey article by Brendle and Marques \cite{BM2}). (In)Stability results for positive solutions in the case of more general potential functions can be found, among others, in the articles by Druet \cite{Dr1}, Esposito, Pistoia and V\'etois \cite{EPV}, Premoselli and Thizy \cites{PT} and Robert and V\'etois \cites{RV4}. 

\smallskip
As regards the existence of sign-changing blowing-up solutions, a well-studied case is the Yamabe equation on the standard sphere; see in this case the historic reference by Ding \cite{Di} and the more recent work by del Pino, Musso, Pacard and Pistoia \cites{dPMPP1,dPMPP2}, Musso and Medina \cite{MM}, Musso and Wei \cite{MW} and Medina, Musso and Wei \cite{MMW}. Few compactness and stability results for sign-changing solutions of equations of type \eqref{IntroEq1} have been obtained so far, which was mainly due to the lack of a satisfactory asymptotic description for energy-bounded sequences of solutions to \eqref{IntroEq2}. Possible references on the topic include the work by Deng, Musso and Wei \cite{DMW}, Pistoia and V\'etois \cite{PiV}, Premoselli \cite{P}, Premoselli and V\'etois \cites{PrV1}, Robert and V\'etois \cites{RV2,RV3} and V\'etois \cite{V}. Finally, other possible references on the existence and multiplicity of sign-changing solutions to the Yamabe equation include the articles by Ammann and Humbert \cite{AH}, Clapp \cite{C}, Clapp and Fern\'andez \cite{CF}, Clapp, Pistoia and Tavares \cite{CPT}, Fernandez and Petean \cite{FP} and Gursky and P\'erez-Ayala \cite{GP}. Needless to say, we do not pretend to any exhaustivity in these lists.

\smallskip
In this article we obtain new stability results extending the one in V\'etois \cite{V} and we construct explicit counter-examples in all dimensions to show their sharpness. Our results highlight new and unexpected blow-up phenomena that are specific to the sign-changing regime.

\subsection{A stability result for locally conformally flat manifolds}

We obtain the following stability result:

\begin{theorem}\label{Th1}
Let $(M,g)$ be a smooth, closed, locally conformally flat Riemannian manifold of dimension $n\ge7$ and $\(h_\eps\)_{0\le\eps\ll1}$ be a family of functions in $C^1\(M\)$ such that $h_\eps\to h_0$ in $C^1\(M\)$ and
\begin{equation}\label{Th1Eq1}
h_0\(x\)\ne c_n\Scal_g\(x\)\quad\forall x\in M.
\end{equation}
Then every family of solutions $\(u_\eps\)_{0<\eps\ll1}$ to \eqref{IntroEq2} that is bounded in $H^1\(M\)$ converges up to a subsequence in $C^2\(M\)$.
\end{theorem}

This result extends previous results obtained by Druet \cite{Dr1} in the case where $u_\eps>0$ in $M$ for all $\eps$ and V\'etois \cite{V} in the case where $h_0<c_n\Scal_g$ in $M$. A simple corollary is that when the assumptions of Theorem~\ref{Th1} are satisfied, for a fixed $C^1$ function $h_0$ satisfying \eqref{Th1Eq1}, every bounded set in $H^1\(M\)$ of (possibly sign-changing) solutions to \eqref{IntroEq1} is precompact in $C^2\(M\)$. 

\smallskip
Remark that if we remove the $H^1$ bound on $\(u_\eps\)_\eps$ while assuming \eqref{Th1Eq1}, then this result does not remain true. Indeed, existence results of families of solutions with unbounded $H^1$ norm have been obtained by V\'etois \cite{V} in the case where $h_0<c_n\Scal_g$ in $M$ and $\(M,g\)$ is locally conformally flat and by Chen, Wei and Yan \cite{CWY} when $n\ge5$, $h_0$ is constant and $h_0>n\(n-2\)/4$ (the solutions obtained in \cite{CWY} are positive; see also the article by V\'etois and Wang \cite{VW} for the $n=4$ counterpart and the article by Premoselli and V\'etois \cite{PrV1} for a construction of sign-changing blowing-up solutions in this case). In the case where $h_0 \equiv c_n \Scal_g$, at least on the standard sphere of dimension $n \ge 3$, it follows from Ding's result \cite{Di} that an analogue of Theorem~\ref{Th1} cannot hold true.

\smallskip
The statement of Theorem \ref{Th1} is very similar to the analogous result obtained by Druet \cite{Dr1} for positive solutions. However, in the sign-changing case, very general sign-changing bubbles can a priori arise in the asymptotic behavior of energy-bounded families of solutions, hence it is actually quite surprising that the condition \eqref{Th1Eq1} remains sufficient to ensure the pre-compactness of such families. This is crucially related to our conformal flatness assumption as is made clear in the proof. In the case of a non-locally conformally flat manifold, it is likely that local necessary conditions for the blow-up of families of solutions to \eqref{IntroEq2} would involve the potential function $h_0$ together with an interplay between the local geometry of the manifold and the arising sign-changing bubbles. It is therefore still unclear what would be a natural generalization of Theorem \ref{Th1} to the non-locally conformally flat setting. 

\smallskip
The assumptions $n \ge 7$ and \eqref{Th1Eq1} in Theorem~\ref{Th1} are sharp. We prove this in Theorems~\ref{Th2} and~\ref{Th3} below. To better understand these counterexamples it is instructive to restrict the focus of Theorem~\ref{Th1} to specific blow-up configurations. In this case, more precise information on the blow-up can be obtained. In the next result, we consider solutions to \eqref{IntroEq2} of the form 
\begin{equation}\label{IntroEq3}
u_\eps=u_0-\(\frac{\sqrt{n\(n-2\)}\mu_\eps}{\mu_\eps^2+d_g\(\cdot,\xi_\eps\)^2}\)^{\frac{n-2}{2}}+\smallo\(1\)\text{ in }H^1\(M\)\text{ as }\eps\to0,
\end{equation}
where $u_0$ is a solution to \eqref{IntroEq1}, $d_g$ is the geodesic distance with respect to the metric $g$ and $\(\mu_\eps\)_{0<\eps\ll1}$ and $\(\xi_\eps\)_{0<\eps\ll1}$ are two families of positive numbers and points in $M$, respectively, such that $\mu_\eps\to0$ and $\xi_\eps\to\xi_0$ as $\eps\to0$ for some point $\xi_0\in M$. In this case, the location of the blow-up point $\xi_0$ is constrained. More precisely, we obtain the following:

\begin{theorem}\label{Th5}
Let $\(M,g\)$ be a smooth, closed Riemannian manifold of dimension $n\ge3$, $h_0\in C^{0,\vartheta}\(M\)$, $\vartheta\in\(0,1\)$, and $u_0$ be a solution to \eqref{IntroEq1}. Assume that $\varphi_0\ne0$ for all points in M, where
\begin{equation}\label{Th5Eq1}
\varphi_0:=\left\{\begin{aligned}&u_0&&\text{if }n\in\left\{3,4,5\right\}\\&h_0+2u_0-c_n\Scal_g&&\text{if }n=6\\&h_0-c_n\Scal_g&&\text{if }n\ge7.\end{aligned}\right.
\end{equation} 
Then there does not exist any families of functions $\(h_\eps\)_{0<\eps\ll1}$ in $C^{0,\vartheta}\(M\)$ and $\(u_\eps\)_{0<\eps\ll1}$ in $C^2\(M\)$, positive numbers $\(\mu_\eps\)_{0<\eps\ll1}$ and points $\(\xi_\eps\)_{0<\eps\ll1}$ in $M$ satisfying \eqref{IntroEq2} and \eqref{IntroEq3} and such that $\mu_\eps\to0$, $\xi_\eps\to\xi_0$ and $h_\eps\to h_0$ in $C^{0,\vartheta}\(M\)$ as $\eps\to0$.
\end{theorem}

In the case of positive solutions (and possibly with multiple bubbles), this result is due to Druet \cite{Dr1}. In particular, it follows from Theorem~\ref{Th5} that if a family $\(u_\eps\)_{0<\eps \ll1}$ of solutions to \eqref{IntroEq2} blows up as in \eqref{IntroEq3}, then $\varphi_0\(\xi_0\)=0$, where $\varphi_0$ is as in \eqref{Th5Eq1}. In dimensions $3,4$ and $5$, for configurations of type \eqref{IntroEq3}, the weak limit $u_0$ thus vanishes at the blow-up point $\xi_0$. Since $u_0$ is now a sign-changing solution to \eqref{IntroEq1}, this does not imply that $u_0$ vanishes on the whole manifold, unlike in the case of positive solutions. Note also that the assumption \eqref{IntroEq3} allows us, in Theorem~\ref{Th5}, to lift the local conformal flatness assumption from Theorem~\ref{Th1} since there cannot be any cancellation phenomena between several sign-changing bubbles in this case. 

\smallskip
Our proofs of Theorems~\ref{Th1}  and ~\ref{Th5} rely on the pointwise estimates for sign-changing solutions recently obtained by Premoselli \cite{P}. The proof of Theorem~\ref{Th1} that we give here is conceptually different from the proof in case $h_0 < c_n \Scal_g$ in the article by V\'etois \cite{V}, which was based on subtle comparison arguments inspired from those in the Euclidean case by Devillanova and Solimini \cite{DS}.

\smallskip
In what follows, we provide counterexamples showing the optimality of the conditions \eqref{Th1Eq1} and $n\ge7$ in Theorem~\ref{Th1}. These counterexamples are obtained by using the Lyapunov--Schmidt reduction method. We obtain two different types of counterexamples, the first in dimensions $n\ge6$ and the second in dimensions $n\in\left\{3,4,5\right\}$, which are inspired from Theorem~\ref{Th5}.

\subsection{Counterexamples: the case of dimensions $n\ge6$}

In the case where $n \ge 6$, we construct solutions of \eqref{IntroEq2} of the form 
$$u_\eps=u_0-\(\frac{\sqrt{n\(n-2\)}\mu_\eps}{\mu_\eps^2+d_g\(\cdot,\xi_\eps\)^2}\)^{\frac{n-2}{2}}+\smallo\(1\)\text{ in }H^1\(M\)\text{ as }\eps\to0,$$
where $u_0$ is a solution to \eqref{IntroEq1}, $d_g$ is the geodesic distance with respect to the metric $g$ and $\(\mu_\eps\)_{0<\eps\ll1}$ and $\(\xi_\eps\)_{0<\eps\ll1}$ are two families of positive numbers and points in $M$, respectively, such that $\mu_\eps\to0$ and $\xi_\eps\to\xi_0$ as $\eps\to0$ for some point $\xi_0\in M$. As Theorem~\ref{Th5} shows, a necessary condition for such solutions to exist is that $\varphi_0\(\xi_0\)=0$, where $\varphi_0$ is as in \eqref{Th5Eq1}.

\smallskip
We say that $u_0$ is nondegenerate if there does not exist any functions $\psi\in C^2\(M\)\backslash\left\{0\right\}$ such that 
$$\Delta_g\psi+h_0\psi=\(2^*-1\)\left|u_0\right|^{2^*-2}\psi\quad\text{in }M.$$
We obtain the following result, which, except in the case where $n=6$ and $u_0\(\xi_0\)>0$, is a straightforward adaptation of a result by Robert and V\'etois \cite{RV4}:

\begin{theorem}\label{Th2}
Let $(M,g)$ be a smooth, closed Riemannian manifold of dimension $n\ge6$, $h_0,h\in C^2\(M\)$, $u_0$ be a nondegenerate solution to \eqref{IntroEq1} and $\varphi_0$ be the function defined in \eqref{Th5Eq1}. Assume that there exists a point $\xi_0\in M$ such that the following two conditions hold true:
\begin{enumerate}
\item[(i)]$\varphi_0\(\xi_0\)=\left|\nabla\varphi_0\(\xi_0\)\right|=0$ and $\left\{\begin{aligned}&D^2\varphi_0\(\xi_0\)<0\text{ if }n=6\text{ and }u_0\(\xi_0\)>0\\&D^2\varphi_0\(\xi_0\)\text{ is nondegenerate otherwise,}\end{aligned}\right.$
\item[(ii)]$K_0\(\xi_0\)\times\left\{\begin{aligned}&\[h-2\(\Delta_g+h_0-2\left|u_0\right|\)^{-1}\(hu_0\)\]\(\xi_0\) &&\text{if }n=6\\&h\(\xi_0\)&&\text{if }n\ge7\end{aligned}\right\}>0$, where
$$K_0:=\left\{\begin{aligned}&\Delta_g\varphi_0+\frac{c_6}{6}\left|\Weyl_g\right|_g^2&&\hbox{if }n=6 \text{ and }u_0\(\xi_0\)<0\\
&u_0&&\hbox{if }n=6\text{ and }u_0\(\xi_0\)\ge0\\
&-u_0&&\hbox{if }7\leq n\leq 9\\
&\Delta_g\varphi_0+\frac{c_{10}}{6}\left|\Weyl_g\right|_g^2-672u_0&&\hbox{if }n=10\\
&\Delta_g\varphi_0+\frac{c_n}{6}\left|\Weyl_g\right|_g^2&&\hbox{if }n\geq 11.\end{aligned}\right.$$
\end{enumerate}\noindent Then there exists a family of solutions $\(u_\eps\)_{0<\eps\ll1}$ to \eqref{IntroEq2} with $h_\eps:=h_0+\eps h$, which satisfies \eqref{IntroEq3} with $\xi_\eps\to\xi_0$ as $\eps\to0$.
\end{theorem}

Theorem~\ref{Th2} shows in particular that the assumption \eqref{Th1Eq1} in Theorem~\ref{Th1} is sharp when $n \ge 7$. Note that we do not assume here that $\(M,g\)$ is locally conformally flat. We refer to Section~\ref{addresults} for extensions of~Theorem~\ref{Th2} and further results in the case where $n \ge 6$. Our results in this case are very close to those obtained in the case of positive solutions in the article by Robert and V\'etois \cite{RV4}. This is not the case for our results in dimensions $n\in\left\{3,4,5\right\}$, which we discuss in what follows. 

\subsection{Counterexamples: the case of dimensions $n\in\left\{3,4,5\right\}$}

We now consider the case where $n\in\left\{3,4,5\right\}$. Our counterexamples in this case are again of the form \eqref{IntroEq3}, but this time the solution $u_0$ changes sign and the point $\xi_0$ is such that $u_0\(\xi_0\)=0$. Here again, the condition $u_0\(\xi_0\)=0$ is necessary for the existence of solutions of type \eqref{IntroEq3} to \eqref{IntroEq1} when $n\in\left\{3,4,5\right\}$ in view of Theorem~\ref{Th5}. Additional difficulties arise in this low-dimensional case. They are related on the one hand to the degeneracy issues in the Lyapunov--Schmidt procedure and the construction of suitable initial functions $h_0$ and $u_0$ and on the other hand to the higher precision required in the expansions due to the condition $u_0\(\xi_0\)=0$.

\smallskip
We obtain our result when $n\in\left\{3,4,5\right\}$ under some additional symmetry assumptions. We recall the following definition, which was used for a different construction in the case of positive solutions by Morabito, Pistoia and Vaira \cite{MPV}: we say that $\(M,g\)$ is symmetric with respect to a point $\xi_0\in M$ if there exists an isometry $\sigma:M\to M$ such that $\sigma\(\xi_0\)=\xi_0$ and $d\sigma_{\xi_0}\(v\)=-v$ for all vectors $v\in T_{\xi_0}M$. Moreover, if $\(M,g\)$ is symmetric with respect to $\xi_0$, then we say that a function $u:M\to M$ is symmetric with respect to $\xi_0$ if $u\circ\sigma=u$ in $M$.

\smallskip
Let us also recall (see for example the article by Li and Zhu \cite{LZhu}) that in dimension $n=3$, for  every smooth metric $g$ on $M$ and every function $h\in C^0\(M\)$, if the operator $\Delta_g+h$ has empty kernel, then there exists a function $m_{g,h}\in C^0\(M\)$, usually called the mass of  the operator $\Delta_g+h$, such that the Green's function $G_{g,h}$ of $\Delta_g+h$ satisfies
\begin{equation}\label{IntroEq5}
G_{g,h}\(x,\xi\)=\(4\pi d_g\(x,\xi\)\)^{-1}+m_{g,h}\(\xi\)+\smallo\(1\)
\end{equation}
as $x\to\xi$. A similar result holds true in dimensions $n\in\left\{4,5\right\}$ when the potential function touches the function $c_n\Scal_g$ at order 1 at some point on the manifold. More precisely, a straightforward extension of a result by Robert and V\'etois \cite{RV4}*{Proposition~8.1} gives that for every smooth metric $g$ on $M$ and every function $h\in C^{1,\theta}\(M\)$, $\vartheta\in\(0,1\)$, if the operator $\Delta_g+h$ has empty kernel and 
$$\(h-c_n\Scal_g\)\(\xi\)=\left|\nabla\(h-c_n\Scal_g\)\(\xi\)\right|=0$$ 
for some point $\xi\in M$, then there exists a function $m_{g,h}\in C^0\(M\)$ such that
\begin{equation}\label{IntroEq6}
G_{g,h}\(x,\xi\)=\(\(n-2\)\omega_{n-1}\)^{-1}\Lambda_\xi\(x\)d_{g_\xi}\(x,\xi\)^{2-n}+m_{g,h}\(\xi\)+\smallo\(1\)
\end{equation}
as $x\to\xi$, where $\omega_{n-1}$ is the volume of the standard unit sphere of dimension $n-1$ and $g_\xi:=\Lambda_\xi^{2^*-2}g$ is a conformal metric to $g$ yielding conformal normal coordinates at the point $\xi$ as in the article by Lee and Parker \cite{LP} (see \eqref{Sec3Eq1}).

\smallskip
With these definitions, we obtain the following:

\begin{theorem}\label{Th3}
Let $\(M,g\)$ be a smooth, closed Riemannian manifold of dimension $n\in\left\{3,4,5\right\}$, $h_0\in C^{1,\theta}\(M\)$, $\vartheta\in\(0,1\)$, and $u_0\not\equiv0$ be a nondegenerate solution to \eqref{IntroEq1}. Assume that there exists a point $\xi_0\in M$ such that $\(M,g\)$, $h_0$ and $u_0$ are symmetric with respect to $\xi_0$ and $u_0\(\xi_0\)=0$. Assume moreover that
\begin{itemize}
\item if $n=3$, then $m_{g, \tilde{h}_0}\(\xi_0\)\ne0$, 
\item if $n\in\left\{4,5\right\}$, then either $h_0\(\xi_0\) \neq c_n\Scal_g\(\xi_0\)$ or $[h_0\(\xi_0\) = c_n\Scal_g\(\xi_0\)$ and $m_{g, \tilde{h}_0}\(\xi_0\)\ne0]$,
\end{itemize} 
where $\tilde{h}_0:=h_0-\(2^*-1\)\left|u_0\right|^{2^*-2}$. Then for every $h\in C^{0,\vartheta'}\(M\)$, $\vartheta'\in\(0,1\)$, that is symmetric with respect to $\xi_0$ and such that
\begin{multline}\label{Th3Eq1}
\big[\big(\Delta_g+\tilde{h}_0\big)^{-1}\(hu_0\)\big]\(\xi_0\)\times\\
\left\{\begin{aligned}&-m_{g, \tilde{h}_0}\(\xi_0\)&&\text{if }n=3\text{ or } [ n=4,5 \text{ and } h_0\(\xi_0\)=c_n\Scal_g\(\xi_0\)]\\&h_0\(\xi_0\)-c_n\Scal_g\(\xi_0\)&&\text{if }n\in\left\{4,5\right\}\text{ and }h_0\(\xi_0\)\ne c_n\Scal_g\(\xi_0\)\end{aligned}\right\}>0,
\end{multline}
there exists a family of solutions $\(u_\eps\)_{0<\eps\ll1}$ to \eqref{IntroEq2} with $h_\eps:=h_0+\eps h$, which satisfies \eqref{IntroEq3} with $\xi_\eps\to\xi_0$ as $\eps\to0$.
\end{theorem}

Note that $\tilde{h}_0\in C^{1,\theta}\(M\)$ for some $\vartheta\in\(0,1\)$ when $n\in\left\{3,4,5\right\}$ and the symmetry assumptions on $\(M,g\)$, $h_0$ and $u_0$ imply that $\nabla h_0\(\xi_0\)=\nabla u_0\(\xi_0\)=\nabla \Scal_g\(\xi_0\)=0$, so the mass $m_{g, \tilde{h}_0}$ is well defined when $h_0\(\xi_0\)=c_n\Scal_g\(\xi_0\)$. Theorem~\ref{Th3} is in stark contrast with the results by Druet \cite{Dr1} for positive solutions. It shows that, in the sign-changing case, blowing-up solutions can exist in dimensions $3$, $4$ and $5$ even when \eqref{Th1Eq1} holds true or the weak limit $u_0$ is not identically zero. Note also that in the case where $n=3$ or [$n\in\left\{4,5\right\}$ and $h_0 \(\xi_0\)=c_n\Scal_g\(\xi_0\)$], the weak limit surprisingly appears in a global constraint via the mass $m_{g, \tilde{h}_0}$. Like in Theorem \ref{Th2}, we do not assume in Theorem \ref{Th3} that $\(M,g\)$ is locally conformally flat.

\smallskip
We refer to Section~\ref{addresults} for existence results of triples $\(M,g\)$, $u_0$ and $h_0$ satisfying the assumptions of Theorem~\ref{Th3}. The non-degeneracy of the solution $u_0$, in particular, is not easily obtained. We construct our examples in the context of $2$-symmetric manifolds, which is defined in Section~\ref{addresults}. By applying Theorem~\ref{Th3} to suitable $2$-symmetric configurations, we obtain in Theorem~\ref{Th4} a general existence result of blowing-up families of solutions of the form \eqref{IntroEq3}. Further results are then obtained in Corollaries~\ref{Cor3} and~\ref{Cor4}.

\subsection{Organization of the article} Our article is organized as follows. Section~\ref{addresults} contains extensions of Theorems~\ref{Th2} and ~\ref{Th3} and further results. In Section~\ref{Sec3}, we prove Theorems~\ref{Th1} and~\ref{Th5}. In Section~\ref{Sec4}, we prove Theorem~\ref{Th2} and Corollaries~\ref{Cor1} and~\ref{Cor2}. In Section~\ref{Sec5}, we  prove Theorems~\ref{Th3} and~\ref{Th4} and Corollaries~\ref{Cor3} and~\ref{Cor4}. 

\section{Additional results}\label{addresults}

\subsection{The case of dimensions $n \ge 6$.} 

By using Theorem~\ref{Th2} together with a perturbation argument, we can prove the following:

\begin{corollary}\label{Cor1}
Let $\(M,g\)$ be a smooth, closed Riemannian manifold of dimension $n\ge6$, $h_0\in C^p\(M\)$, $1\le p\le\infty$, $u_0$ be a positive solution to \eqref{IntroEq1} and $\varphi_0$ be the function defined in \eqref{Th5Eq1}. Assume that there exists a point $\xi_0\in M$ such that $\varphi_0\(\xi_0\)=\left|\nabla\varphi_0\(\xi_0\)\right|=0$.  In the case where $n=6$, assume moreover that $\xi_0$ is a local maximum point of the function $\varphi_0$. Then there exist families of functions $\(h_\eps\)_{0<\eps\ll1}$ in $C^p\(M\)$ and $\(u_\eps\)_{0<\eps\ll1}$ in $C^2\(M\)$ satisfying \eqref{IntroEq2} and \eqref{IntroEq3} with $\xi_\eps\to\xi_0$ as $\eps\to0$ and such that $h_\eps\to h_0$ in $C^p\(M\)$ as $\eps\to0$.
\end{corollary}

Remark that in Corollary~\ref{Cor1} we do not assume that $u_0$ is non-degenerate. We establish this property in the proof by a suitable choice of the perturbation $\(h_\eps\)_{0<\eps\ll 1}$. By combining Corollary~\ref{Cor1} with well-known results by Aubin \cite{A}, Schoen \cite{Sc1} and Trudinger \cite{T}, we then obtain the following:

\begin{corollary}\label{Cor2}
Let $\(M,g\)$ be a smooth, closed Riemannian manifold of dimension $n\ge6$, $\xi_0\in M$, $1\le p\le\infty$ and $\varphi\in C^p\(M\)$ be such that $\varphi\(\xi_0\)=\left|\nabla\varphi\(\xi_0\)\right|=0$ and the operator $\Delta_g+c_n\Scal_g+\varphi$ is coercive in $H^1\(M\)$. In the case where $n=6$, assume moreover that $\xi_0$ is a local maximum point of the function $\varphi_0$. Then the following two assertions hold true:
\begin{enumerate}
\item[(i)] If $\varphi<0$ at some point in $M$, then there exist families of functions $\(h_\eps\)_{0\le\eps\ll1}$ in $C^p\(M\)$ and $\(u_\eps\)_{0\le\eps\ll1}$ in $C^2\(M\)$ satisfying \eqref{IntroEq2} and \eqref{IntroEq3} with $\xi_\eps\to\xi_0$ as $\eps\to0$ and such that $h_\eps\to h_0$ in $C^p\(M\)$ as $\eps\to0$, $u_0>0$ in $M$ and $\varphi_0=\varphi$ in $M$, where $\varphi_0$ is as in \eqref{Th5Eq1}.
\item[(ii)] If $\(M,g\)$ is not conformally diffeomorphic to the standard sphere, then there exists a positive constant $\eps_0$ depending only on $\(M,g\)$ such that if $\left\|\max\(\varphi,0\)\right\|_{L^1\(M\)}<\eps_0$, then there exist families of functions $\(h_\eps\)_{0\le\eps\ll1}$ in $C^p\(M\)$ and $\(u_\eps\)_{0\le\eps\ll1}$ in $C^2\(M\)$ satisfying \eqref{IntroEq2} and \eqref{IntroEq3} with $\xi_\eps\to\xi_0$ as $\eps\to0$ and such that $h_\eps\to h_0$ in $C^p\(M\)$ as $\eps\to0$, $u_0>0$ in $M$ and $\varphi_0=\varphi$ in $M$, where $\varphi_0$ is as in \eqref{Th5Eq1}.
\end{enumerate}
\end{corollary}

Corollaries~\ref{Cor1} and~\ref{Cor2} are proven in Section~\ref{Sec4}.

\subsection{The case of dimensions $n \in \{3,4,5\}$.}

To obtain examples of situations where we can apply Theorem~\ref{Th3}, we now introduce some additional symmetry assumptions. We say that $\(M,g\)$ is $2$-symmetric with respect to a point $\xi_0\in M$ if $M$ is connected and there exist $3$ isometries $\sigma_0,\sigma_1,\sigma_2:M\to M$ and a domain $\Omega\subset M$ such that 
$$\left\{\begin{aligned}&\text{$\sigma_0\circ\sigma_0=\sigma_1\circ\sigma_1=\sigma_2\circ\sigma_2=\id_M$ in $M$,}\\&\text{$\sigma_0\circ\sigma_1=\sigma_1\circ\sigma_0$, $\sigma_0\circ\sigma_2=\sigma_2\circ\sigma_0$ and $\sigma_1\circ\sigma_2=\sigma_2\circ\sigma_1$ in $M$,}\\&\text{$\xi_0\in\Gamma_{\sigma_0}\cap\Gamma_{\sigma_1}\cap\Gamma_{\sigma_2}$ and $d\(\sigma_0\circ\sigma_1\circ\sigma_2\)_{\xi_0}\(v\)=-v$ for all vectors $v\in T_{\xi_0}M$,}\\&\text{$\sigma_0\(\Omega\)  = \Omega$ and $M=\Omega\sqcup\sigma_1\(\Omega\)\sqcup\sigma_2\(\Omega\)\sqcup\sigma_1\circ\sigma_2\(\Omega\)\sqcup\(\Gamma_{\sigma_1}\cup\Gamma_{\sigma_2}\)$}\\\end{aligned}\right.$$
where $\sqcup$ is the disjoint union and $\Gamma_{\sigma_i}:=\left\{x\in M:\,\sigma_i\(x\)=x\right\}$ for $i\in\left\{0,1,2\right\}$. Remark that it easily follows from these properties that $\xi_0\in\partial\Omega$. Moreover, if $\(M,g\)$ is $2$-symmetric with respect to $\xi_0$, then we say that a function $h:M\to M$ is $2$-symmetric with respect to $\xi_0$ if $h\circ\sigma_0=h\circ\sigma_1=h\circ\sigma_2=h$ in $M$. Clearly, if $\(M,g\)$ and $h$ are $2$-symmetric with respect to $\xi_0$, then by taking $\sigma:=\sigma_0\circ\sigma_1\circ\sigma_2$, we obtain that $\(M,g\)$ and $h$ are symmetric with respect to $\xi_0$.

\smallskip
A first example of a $2$-symmetric manifold is given by the standard unit sphere $\(\S^n,g_0\)$ with $n \ge 2$ and 
$$\left\{\begin{aligned}&\xi_0:=\(0,\dotsc,0,1\),\,\Omega:=\left\{x\in\S^n:\,x_1>0\text{ and }x_2>0\right\},\\&\sigma_0\(x\):=\(x_1,x_2,-x_3,\dotsc,-x_n,x_{n+1}\),\,\sigma_1\(x\):=\(-x_1,x_2,x_3,\dotsc,x_n,x_{n+1}\)\\&\text{and }\sigma_2\(x\):=\(x_1,-x_2,x_3,\dotsc,x_n,x_{n+1}\)\end{aligned}\right.$$
for all points $x:=\(x_1,\dotsc,x_{n+1}\)\in\S^n$. Another example is given by the product of spheres $\S^{n_1}\times\S^{n_2}\subset\R^{n_1+1}\times\R^{n_2+1}$ equipped with the standard product metric and where
$$\left\{\begin{aligned}&\xi_0:=\(\(0,\dotsc,0,1\),\(0,\dotsc,0,1\)\),\\
&\Omega:=\left\{\(x,y\)\in\S^{n_1}\times\S^{n_2}:\,x_1>0\text{ and }y_1>0\right\},\\&
\sigma_0\(x,y\):=\(\(x_1,-x_2,\dotsc,-x_{n_1},x_{n_1+1}\),\(y_1,-y_2,\dotsc,-y_{n_1},y_{n_1+1}\)\),\\&\sigma_1\(x,y\):=\(\(-x_1,x_2,\dotsc,x_{n_1},x_{n_1+1}\),y\)\\
&\text{and }\sigma_2\(x,y\):=\(x,\(-y_1,y_2,\dotsc,y_{n_1},y_{n_1+1}\)\)\end{aligned}\right.$$
for all points $\(x,y\):=\(\(x_1,\dotsc,x_{n_1+1}\),\(y_1,\dotsc,y_{n_2+1}\)\)\in\S^{n_1}\times\S^{n_2}$.

\smallskip
For $2$-symmetric manifolds, we are able to reduce the construction of functions $h_0$ and $u_0$ satisfying the assumptions of Theorem~\ref{Th3} to constructing a positive solution $u_0$ to the problem
\begin{equation}\label{IntroEq7}
\left\{\begin{aligned}&\Delta_gu_0+h_0u_0=u_0^{2^*-1}\quad\text{in }\Omega\\&u_0\in H^1_{0,\sigma_0}\(\Omega\),\end{aligned}\right.
\end{equation}
where 
$$H^1_{0,\sigma_0}\(\Omega\):=\left\{u\in H^1_0\(\Omega\),:\,u\circ\sigma_0=u\text{ in }\Omega\right\}.$$
More precisely, we obtain the following:

\begin{theorem}\label{Th4}
Let $\(M,g\)$ be a smooth, closed Riemannian manifold of dimension $n\in\left\{3,4,5\right\}$ and $h_0\in C^p\(M\)$, $1\le p\le\infty$. Assume that there exists a point $\xi_0\in M$ such that $\(M,g\)$ and $h_0$ are $2$-symmetric with respect to $\xi_0$ and there exists a positive solution $u_0\in H^1_{0,\sigma_0}\(\Omega\)$ to \eqref{IntroEq7}. Then we can extend $u_0$ as a sign-changing solution to \eqref{IntroEq1} and there exist families of functions $\(h_\eps\)_{0<\eps\ll1}$ in $C^p\(M\)$ and $\(u_\eps\)_{0<\eps\ll1}$ in $C^2\(M\)$ satisfying \eqref{IntroEq2} and \eqref{IntroEq3} with $\xi_\eps\to\xi_0$ as $\eps\to0$ and such that $h_\eps\to h_0$ in $C^p\(M\)$ as $\eps\to0$.
\end{theorem}

By combining Theorem~\ref{Th4} with a similar result as in the article by Aubin \cite{A}, we obtain the following:

\begin{corollary}\label{Cor3}
Let $\(M,g\)$ be a smooth, closed Riemannian manifold of dimension $n\in\left\{3,4,5\right\}$ and $h_0\in C^p\(M\)$, $1\le p\le\infty$. Assume that there exists a point $\xi_0\in M$ such that $\(M,g\)$ and $h_0$ are $2$-symmetric with respect to $\xi_0$. Assume moreover that the operator $\Delta_g+h_0$ is coercive in $H^1_{0,\sigma_0}\(\Omega\)$ and $h_0\(\xi_1\)<c_n\Scal_g\(\xi_1\)$ at some point $\xi_1\in\Omega\cap\Gamma_{\sigma_0}$ (remark that $\Omega\cap\Gamma_{\sigma_0}\ne\emptyset$; see Proposition~\ref{Pr1}). Then there exists a solution $u_0\not\equiv0$ to \eqref{IntroEq1} such that $u_0\(\xi_0\)=0$ and there exist families of functions $\(h_\eps\)_{0<\eps\ll1}$ in $C^p\(M\)$ and $\(u_\eps\)_{0<\eps\ll1}$ in $C^2\(M\)$ satisfying \eqref{IntroEq2} and \eqref{IntroEq3} with $\xi_\eps\to\xi_0$ as $\eps\to0$ and such that $h_\eps\to h_0$ in $C^p\(M\)$ as $\eps\to0$.
\end{corollary}

In view of Theorem~\ref{Th1} and Druet's results \cite{Dr1} for positive solutions, we are  particularly interested, when $n\in\{3,4,5\}$, in obtaining blowing-up counterexamples when $h_0 \ge c_n\Scal_g$ in $M$. For this, we assume moreover that there exists a group $G$ of isometries of $M$ such that
\begin{equation}\label{IntroEq8}
\left\{\begin{aligned}&\text{the set }\left\{\sigma\(\xi\):\,\sigma\in G\right\}\text{ is infinite for all points }\xi\in\overline\Omega,\\&h_0\circ\sigma=h_0\text{ in }\Omega\text{ and }\sigma\(\Omega\)=\Omega\text{ for all isometries }\sigma\in G.\end{aligned}\right.
\end{equation}
For example, given any $2$-symmetric manifold $\(M^{n_1},g_1\)$ of dimension $n_1\ge2$, letting $\(\S^{n_0},g_0\)$ be the standard unit sphere of dimension $n_0\ge1$, it is not difficult to see that the product manifold $\(\S^{n_0}\times M^{n_1},g_0\times g_1\)$ is $2$-symmetric when equipped with 
$$\left\{\begin{aligned}&\xi_0:=\(\(0,\dotsc,0,1\),\xi_0^{n_1}\),\,\Omega:=\S^{n_0}\times\Omega^{n_1},\\&
\sigma_0\(x,y\):=\(\(-x_1,\dotsc,-x_{n_0},x_{n_0+1}\),\sigma_0^{n_1}\(y\)\),\\&\sigma_1\(x,y\):=\(x,\sigma_1^{n_1}\(y\)\)\text{ and }\sigma_2\(x,y\):=\(x,\sigma_2^{n_1}\(y\)\)\end{aligned}\right.$$
for all points $x:=\(x_1,\dotsc,x_{n_0+1}\)\in\S^{n_0}$ and $y\in M^{n_1}$, where $\xi_0^{n_1}$, $\Omega^{n_1}$, $\sigma_0^{n_1}$, $\sigma_1^{n_1}$ and $\sigma_2^{n_1}$ correspond to $M^{n_1}$ and $\xi_0$, $\Omega$, $\sigma_0$, $\sigma_1$ and $\sigma_2$ correspond to $\S^{n_0}\times M^{n_1}$. In this case, letting $\Isom\(\S^{n_0},g_0\)$ be the group of isometries of $\(\S^{n_0},g_0\)$, we can define a group 
$$G=\left\{\left\{\begin{array}{ccc}\S^{n_0}\times M^{n_1}&\to&\S^{n_0}\times M^{n_1}\\\(x,y\)&\mapsto&\(\sigma\(x\),y\)\end{array}\right\}:\,\sigma\in \Isom\(\S^{n_0},g_0\)\right\}$$
of isometries of $\S^{n_0}\times M^{n_1}$ which clearly satisfies \eqref{IntroEq8}. 

\smallskip
By using Theorem~\ref{Th4} together with this additional assumption, we obtain the following:

\begin{corollary}\label{Cor4}
Let $\(M,g\)$ be a smooth, closed Riemannian manifold of dimension $n\in\left\{3,4,5\right\}$ and $h_0\in C^p\(M\)$, $1\le p\le\infty$. Assume that there exists a point $\xi_0\in M$ such that $\(M,g\)$ and $h_0$ are $2$-symmetric with respect to $\xi_0$. Assume moreover that the operator $\Delta_g+h_0$ is coercive in $H^1_{0,\sigma_0}\(\Omega\)$ and there exists a group $G$ of isometries of $M$ satisfying \eqref{IntroEq8}. Then there exists a solution $u_0\not\equiv0$ to \eqref{IntroEq1} such that $u_0\(\xi_0\)=0$ and there exist families of functions $\(h_\eps\)_{0<\eps\ll1}$ in $C^p\(M\)$ and $\(u_\eps\)_{0<\eps\ll1}$ in $C^2\(M\)$ satisfying \eqref{IntroEq2} and \eqref{IntroEq3} with $\xi_\eps\to\xi_0$ as $\eps\to0$ and such that $h_\eps\to h_0$ in $C^p\(M\)$ as $\eps\to0$.
\end{corollary}

Theorems~\ref{Th3} and~\ref{Th4} and Corollaries~\ref{Cor3} and~\ref{Cor4} are proven in Section~\ref{Sec5}.

\section{Proof of Theorem~\ref{Th1}}\label{Sec3}

This section is devoted to the proof of Theorems~\ref{Th1} and~\ref{Th5}. We first prove Theorem ~\ref{Th1}.  

\proof[Proof of Theorem~\ref{Th1}]
Assume by contradiction that there exists a family of solutions $\(u_\eps\)_{0\le\eps\ll1}$ to \eqref{IntroEq2} that is bounded in $H^1\(M\)$ and such that no subsequences of $\(u_\eps\)_{\eps}$ converge in $C^2\(M\)$. By standard elliptic theory, we obtain that up to a subsequence, $\(u_\eps\)_{\eps}$ blows up as $\eps\to0$. We can then apply the main results in the articles by Premoselli \cite{P} and Struwe \cite{St}, which give that there exist a subsequence $\(u_\alpha\)_{\alpha\in\N}$ of $\(u_\eps\)_\eps$, a solution $u_0\in C^2\(M\)$ to \eqref{IntroEq1}, a number $k\in\N$ and for every $i\in\left\{1,\dotsc,k\right\}$, a sequence of positive numbers $\(\mu_{i,\alpha}\)_{\alpha\in\N}$, a sequences of points $\(\xi_{i,\alpha}\)_{\alpha\in\N}$ in $M$ and a nontrivial solution $V_i\in D^{1,2}\(\R^n\)$ to the equation
$$\Delta_{\delta_0}V_i=\left|V_i\right|^{2^*-2}V_i\quad\text{in }\R^n,$$
where $\delta_0$ is the Euclidean metric, such that 
$$u_\alpha=u_0+\sum_{i=1}^kV_{i,\alpha}+\smallo\(1\)\quad\text{in }H^1\(M\)$$
and
\begin{equation}\label{Th1Eq2}
\left|u_\alpha-u_0-\sum_{i=1}^kV_{i,\alpha}\right|=\smallo\(1+\sum_{i=1}^kB_{i,\alpha}\)\quad\text{ uniformly in }M
\end{equation}
as $\alpha\to\infty$, where
$$B_{i,\alpha}\(x\):=\(\frac{\sqrt{n\(n-2\)}\mu_{i,\alpha}}{\mu_{i,\alpha}^2+d_g\(x,\xi_{i,\alpha}\)^2}\)^{\frac{n-2}{2}}\quad\forall x\in M$$
and the bubbling profiles $V_{i,\alpha}$ are defined as follows:
\begin{itemize}
\item In the case where  $u_0 \not \equiv 0$ or $\ker\(\Delta_g + h_0\) \neq  \left\{0\right\}$, we let
\begin{equation*}
V_{i, \alpha}\(x \):=   \chi \( d_{g}\(\xi_{i,\alpha}, x\) \) \mu_{i,\alpha}^{-\frac{n-2}{2}} V_i \( \frac{1}{\mu_{i,\alpha}} \exp_{\xi_{i,\alpha}}^{-1}\(x\) \)\quad\forall x\in M,
\end{equation*} 
where $\exp_{\xi_{i,\alpha}}$ is the exponential map at the point $\xi_{i,\alpha}$ with respect to the metric $g$, we identify $T_{\xi_{i,\alpha}}M$ with $\R^n$ and $\chi$ is a smooth, nonnegative cutoff function in $\[0,\infty\)$ such that $\chi=1$ in $\[0,r_0/2\)$ and $\chi=0$ in $\[r_0,\infty\)$ for some number $r_0>0$ smaller than the injectivity radius of the metric $g$. 
\item In the case where $u_0 \equiv 0$ and $\ker(\Delta_g + h_0) =  \left\{0\right\}$, we define
$$ F\(\xi_{i,\alpha}, x\):=  \(n-2\)\omega_{n-1} d_g\(\xi_{i,\alpha}, x\)^{n-2} G_{g,h_0}\(\xi_{i,\alpha}, x\)\quad\forall x\in M,$$
where $\omega_{n-1}$ is the volume of the standard unit sphere of dimension $n-1$ and $G_{g,h_0}$ is the Green's function of the operator $\Delta_g + h_0$. We then let
\begin{multline*} \label{defbulle2}
V_{i, \alpha}\(x \):=  \chi \( d_{g}\(\xi_{i,\alpha}, x\) \)F\(\xi_{i,\alpha}, x\) \mu_{i,\alpha}^{-\frac{n-2}{2}}V_i \( \frac{1}{\mu_{i,\alpha}} \exp_{\xi_{i,\alpha}}^{-1}\(x\) \)  \\
 +\(1 -  \chi \( d_{g}\(\xi_{i,\alpha}, x\) \) \)\(n-2\)\omega_{n-1} \lambda_i \mu_{i, \alpha}^{\frac{n-2}{2}} G_{g,h_0}\(\xi_{i,\alpha}, x\)\quad\forall x\in M,
\end{multline*}
where $\exp_{\xi_{i,\alpha}}$ and $\chi$ are as in the first case and $\lambda_i \in \mathbb{R}$ is such that
$$ V_i\(y\) =  \lambda_i\left|y\right|^{2-n} + \bigO \( \left|y\right|^{1-n}\) \textrm{ as } |y| \to + \infty,$$
where $y \in \mathbb{R}^n$ and $\left|y\right|:=\left|y\right|_{\delta_0}$.
\end{itemize}
Moreover, the sequences $\(\mu_{i,\alpha}\)_{\alpha\in\N}$ and $\(\xi_{i,\alpha}\)_{\alpha\in\N}$ satisfy
\begin{equation}\label{Th1Eq3}
\frac{d_g\(\xi_{i,\alpha},\xi_{j,\alpha}\)^2}{\mu_{i,\alpha}\mu_{j,\alpha}}+\frac{\mu_{i,\alpha}}{\mu_{j,\alpha}}+\frac{\mu_{j,\alpha}}{\mu_{i,\alpha}}\longrightarrow\infty
\end{equation}
as $\alpha\to\infty$, for all $i,j\in\left\{1,\dotsc,k\right\}$ such that $i\ne j$. By renumbering and passing to a subsequence if necessary, we may assume that $\mu_{1,\alpha}\ge\dotsb\ge\mu_{k,\alpha}$ for all $\alpha\in\N$ and $\xi_{1,\alpha}\to\xi_{1,0}$ as $\alpha\to\infty$ for some point $\xi_{1,0}\in M$. Since $\(M,g\)$ is locally conformally flat, by choosing $r_0$ smaller if necessary, we may assume moreover that there exists a positive function $\Lambda\in C^\infty\(M\)$ such that the metric $\hat{g}:=\Lambda^{2^*-2}g$ is flat in $B_g\(\xi_{1,\alpha},r_0\)$ for all $i\in\left\{1,\dotsc,k\right\}$ and $\alpha\in\N$, where $B_g\(\xi_{1,\alpha},r_0\)$ is the ball of center $\xi_{1,\alpha}$ and radius $r_0$ with respect to the metric $g$. By using the conformal invariance property of the conformal Laplacian together with the change of variable $x=\widehat{\exp}_{\xi_{1,\alpha}}\(\mu_{1,\alpha}x\)$ for $y\in B\(0,r_0/\mu_{1,\alpha}\)$, where $\widehat{\exp}_{\xi_{1,\alpha}}$ is the exponential map at the point $\xi_{1,\alpha}$ with respect to the metric $\hat{g}$ and $B\(0,r_0/\mu_{1,\alpha}\)$ is the ball of center 0 and radius $r_0/\mu_{1,\alpha}$ with respect to the Euclidean metric, the equation \eqref{IntroEq2} can then be rewritten as
$$\Delta_{\delta_0}\hat{u}_\alpha+\mu_{1,\alpha}^2\hat{h}_\alpha\hat{u}_\alpha=\left|\hat{u}_\alpha\right|^{2^*-2}\hat{u}_\alpha\quad\text{in }B\(0,r_0/\mu_{1,\alpha}\),$$
where 
$$\hat{u}_\alpha\(y\):=\mu_{1,\alpha}^{\frac{n-2}{2}}\big(\Lambda^{-1}u_\alpha\big)\big(\widehat{\exp}_{\xi_{1,\alpha}}\(\mu_{1,\alpha} y\)\big)$$
and
$$\hat{h}_\alpha\(y\):=\big(\Lambda^{2-2^*}\(h_\alpha-c_n\Scal_g\)\big)\big(\widehat{\exp}_{\xi_{1,\alpha}}\(\mu_{1,\alpha} y\)\big)$$
for all points $y\in B\(0,r_0/\mu_{1,\alpha}\)$. For every $\varrho>0$, by using a Pohozaev-type identity in $B\(0,\varrho/ \sqrt{\mu_{1,\alpha}}\)$ (see for example the article by Marques \cite{M}*{Formula~(2.8)}), we then obtain
\begin{multline}\label{Th1Eq4}
\int_{B\(0,\varrho/\sqrt{\mu_{1,\alpha}}\)}\(\hat{h}_\alpha+\frac{1}{2}\big<\nabla\hat{h}_\alpha,\cdot\big>_{\delta_0}\)\hat{u}_\alpha^2\,dv_{\delta_0}=\mu_{1,\alpha}^{-5/2}\int_{\partial B\(0,\varrho/\sqrt{\mu_{1,\alpha}}\)}\Big(\frac{\varrho}{2}\left|\nabla\hat{u}_\alpha\right|_{\delta_0}^2\\
-\varrho\(\partial_\nu\hat{u}_\alpha\)^2-\frac{n-2}{2}\sqrt{\mu_{1,\alpha}}\hat{u}_\alpha\partial_\nu\hat{u}_\alpha-\frac{\varrho}{2^*}\left|\hat{u}_\alpha\right|^{2^*}+\frac{\varrho}{2}\mu_{1,\alpha}^2\hat{h}_\alpha\hat{u}_\alpha^2\Big)d\sigma_{\delta_0},
\end{multline}
where $\nu$ and $d\sigma_{\delta_0}$ are the outward unit normal vector and volume element, respectively, of the Euclidean metric on $\partial B\(0,\varrho/\sqrt{\mu_{1,\alpha}}\)$. To estimate the terms in the right-hand side of \eqref{Th1Eq4}, we use another change of variable. We define
$$\check{u}_\alpha\(y\):=\big(\Lambda^{-1}u_\alpha\big)\big(\widehat{\exp}_{\xi_{1,\alpha}}\(\sqrt{\mu_{1,\alpha}}\,y\)\big)$$ 
and
$$\check{h}_\alpha\(y\):=\big(\Lambda^{2-2^*}\(h_\alpha-c_n\Scal_g\)\big)\big(\widehat{\exp}_{\xi_{1,\alpha}}\(\sqrt{\mu_{1,\alpha}}y\)\big)$$
for all points $y\in B\(0,r_0/\sqrt{\mu_{1,\alpha}}\)$. Then we have
\begin{equation}\label{Th1Eq5}
\Delta_{\delta_0}\check{u}_\alpha+\mu_{1,\alpha}\check{h}_\alpha\check{u}_\alpha=\mu_{1,\alpha}\left|\check{u}_\alpha\right|^{2^*-2}\check{u}_\alpha\quad\text{in }B\(0,r_0/\sqrt{\mu_{1,\alpha}}\)
\end{equation}
and
\begin{multline}\label{Th1Eq6}
\int_{\partial B\(0,\varrho/\sqrt{\mu_{1,\alpha}}\)}\bigg(\frac{\varrho}{2}\left|\nabla\hat{u}_\alpha\right|_{\delta_0}^2-\varrho\(\partial_\nu\hat{u}_\alpha\)^2-\frac{n-2}{2}\sqrt{\mu_{1,\alpha}}\hat{u}_\alpha\partial_\nu\hat{u}_\alpha-\frac{\varrho}{2^*}\left|\hat{u}_\alpha\right|^{2^*}\\
+\frac{\varrho}{2}\mu_{1,\alpha}^2\hat{h}_\alpha\hat{u}_\alpha^2\bigg)d\sigma_{\delta_0}=\mu_{1,\alpha}^{\frac{n-1}{2}}\int_{\partial B\(0,\varrho\)}\\
\(\frac{\varrho}{2}\left|\nabla\check{u}_\alpha\right|_{\delta_0}^2-\varrho\(\partial_\nu\check{u}_\alpha\)^2-\frac{n-2}{2}\check{u}_\alpha\partial_\nu\check{u}_\alpha-\frac{\varrho}{2^*}\mu_{1,\alpha}\left|\check{u}_\alpha\right|^{2^*}+\frac{\varrho}{2}\mu_{1,\alpha}\check{h}_\alpha\hat{u}_\alpha^2\)d\sigma_{\delta_0},
\end{multline}
 Passing to a subsequence if necessary, we may assume that there exists $\varrho_0>0$ such that for every $i\in\left\{2,\dotsc,k\right\}$, either  $d_{\hat{g}}\(\xi_{1,\alpha},\xi_{i,\alpha}\)=\smallo\(\sqrt{\mu_{1,\alpha}}\)$ as $\alpha\to\infty$ or $d_{\hat{g}}\(\xi_{1,\alpha},\xi_{i,\alpha}\)\ge2\varrho_0\sqrt{\mu_{1,\alpha}}$ for all $\alpha\in\N$. We then let $I_1$ be the set of indices $i\in\left\{1,\dotsc,k\right\}$ such that $d_{\hat{g}}\(\xi_{1,\alpha},\xi_{i,\alpha}\)=\smallo\(\sqrt{\mu_{1,\alpha}}\)$ as $\alpha\to\infty$, $I_2:=\left\{1,\dotsc,k\right\}\backslash I_1$ and 
$$\varrho_\alpha:=2\mu_{1,\alpha}^{-1/2}\max_{i\in I_1}d_{\hat{g}}\(\xi_{1,\alpha},\xi_{i,\alpha}\)$$
for all $\alpha\in\N$. By using \eqref{Th1Eq2} and remarking that 
\begin{equation}\label{Th1Eq7}
\left|V_{i,\alpha}\(x\)\right|=\bigO\(B_{i,\alpha}\(x\)\)
\end{equation}
uniformly with respect to $x\in M$, $i\in\left\{1,\dotsc,k\right\}$ and $\alpha\in\N$, we then obtain
\begin{align}\label{Th1Eq8}
&\left|\check{u}_\alpha\(y\)\right|=\bigO\(1+\sum_{i=1}^k\mu_{i,\alpha}^{\frac{n-2}{2}}d_{\hat{g}}\big(\exp_{\xi_{1,\alpha}}\(\sqrt{\mu_{1,\alpha}}\,y\),\xi_{i,\alpha}\big)^{2-n}\)\nonumber\allowdisplaybreaks\\
&\qquad=\bigO\(1+\sum_{i\in I_1}\(\frac{\mu_{i,\alpha}}{\mu_{1,\alpha}}\)^{\frac{n-2}{2}}\(\left|y\right|-\frac{\varrho_\alpha}{2}\)^{2-n}+\sum_{i\in I_2}\(\frac{\mu_{i,\alpha}}{\mu_{1,\alpha}}\)^{\frac{n-2}{2}}\(2\varrho_0-\left|y\right|\)^{2-n}\)\nonumber\\
&\qquad=\bigO\big(\left|y\right|^{2-n}\big)
\end{align}
uniformly with respect to $\alpha\in\N$ and $y\in B\(0,\varrho_0\)\backslash B\(0,\varrho_\alpha\)$. By standard elliptic theory and since $\varrho_\alpha\to0$ and $h_\alpha\to h_0$ in $C^1\(M\)$ as $\alpha\to\infty$, it follows from \eqref{Th1Eq5} and \eqref{Th1Eq8} that $\(\check{u}_\alpha\)_\alpha$ is bounded in $C^1_{\loc}\(B\(0,\varrho_0\)\backslash\left\{0\right\}\)$. In particular, we obtain
\begin{multline}\label{Th1Eq9}
\int_{\partial B\(0,\varrho_0\)}\bigg(\frac{\varrho_0}{2}\left|\nabla\check{u}_\alpha\right|_{\delta_0}^2-\varrho_0\(\partial_\nu\check{u}_\alpha\)^2-\frac{n-2}{2}\check{u}_\alpha\partial_\nu\check{u}_\alpha-\frac{\varrho_0}{2^*}\mu_{1,\alpha}\left|\check{u}_\alpha\right|^{2^*}\\
+\frac{\varrho_0}{2}\mu_{1,\alpha}\check{h}_\alpha\hat{u}_\alpha^2\bigg)d\sigma_{\delta_0}=\bigO\(1\)
\end{multline}
uniformly with respect to $\alpha\in\N$. We now consider the terms in the left-hand side of \eqref{Th1Eq4}. Since $h_\alpha\to h_0$ in $C^1\(M\)$ as $\alpha\to\infty$, we obtain
\begin{equation}\label{Th1Eq10}
\hat{h}_\alpha+\frac{1}{2}\big<\nabla\hat{h}_\alpha,\cdot\big>_{\delta_0}=\big(\Lambda^{2-2^*}\(h_0- c_n\Scal_g\)\big)\(\xi_{1,\alpha}\)+\smallo\(1\)
\end{equation}
as $\alpha\to\infty$, uniformly in $B\(0,\varrho_0/\sqrt{\mu_{1,\alpha}}\)$. By putting together \eqref{Th1Eq1}, \eqref{Th1Eq4}, \eqref{Th1Eq6}, \eqref{Th1Eq9} and \eqref{Th1Eq10}, we then obtain
\begin{equation}\label{Th1Eq11}
\int_{B\(0,\varrho_0/\sqrt{\mu_{1,\alpha}}\)}\hat{u}_\alpha^2\,dv_{\delta_0}=\bigO\(\mu_{1,\alpha}^{\frac{n-6}{2}}\)
\end{equation}
uniformly with respect to $\alpha\in\N$. On the other hand, for every $R>0$, by remarking that $B_g\(\xi_{1,\alpha},R\mu_{1,\alpha}\)\subset B_{\hat{g}}\(\xi_{1,\alpha},\varrho_0\sqrt{\mu_{1,\alpha}}\)$, we obtain
\begin{equation}\label{Th1Eq12}
\int_{B\(0,\varrho_0/\sqrt{\mu_{1,\alpha}}\)}\hat{u}_\alpha^2\,dv_{\delta_0}\ge C\int_{B\(0,R\)}\tilde{u}_\alpha^2\,dv_{\tilde{g}_\alpha}
\end{equation}
for some constant $C>0$ independent of $\alpha$, where
$$\tilde{u}_\alpha\(y\):=\mu_{1,\alpha}^{\frac{n-2}{2}}u_\alpha\big(\exp_{\xi_{1,\alpha}}\(\mu_{1,\alpha} y\)\big)\quad\text{and}\quad \tilde{g}_\alpha\(y\):=\exp_{\xi_{1,\alpha}}^*g\(\mu_{1,\alpha}y\)$$
for all $y\in B\(0,R\)$. Since $\mu_{i,\alpha}\le\mu_{1,\alpha}$ for all $\alpha\in\N$ and $i\in\left\{2,\dotsc,k\right\}$, by using \eqref{Th1Eq3} and passing to a subsequence if necessary, we may assume that for every $i\in\left\{2,\dotsc,k\right\}$, we have $\mu_{1,\alpha}=\smallo\(d_g\(\xi_{1,\alpha},\xi_{i,\alpha}\)\)$ or $\mu_{i,\alpha}=\smallo\(\mu_{1,\alpha}\)$ as $\alpha\to\infty$. By passing again to a subsequence, we may assume moreover that there exist two numbers $R_1,R_2>0$ such that $R_1<R_2$ and for every $i\in\left\{2,\dotsc,k\right\}$, either $d_g\(\xi_{1,\alpha},\xi_{i,\alpha}\)\le R_1\mu_{1,\alpha}$ for all $\alpha\in\N$ or $d_g\(\xi_{1,\alpha},\xi_{i,\alpha}\)\ge R_2\mu_{1,\alpha}$ for all $\alpha\in\N$. We then let $J_1$ be the set of indices $i\in\left\{2,\dotsc,k\right\}$ such that $\mu_{1,\alpha}=\smallo\(d_g\(\xi_{1,\alpha},\xi_{i,\alpha}\)\)$ as $\alpha\to\infty$, $J_2$ be the set of indices $i\in\left\{2,\dotsc,k\right\}$ such that $\mu_{i,\alpha}=\smallo\(\mu_{1,\alpha}\)$ as $\alpha\to\infty$ and $d_g\(\xi_{1,\alpha},\xi_{i,\alpha}\)\le R_1\mu_{1,\alpha}$ for all $\alpha\in\N$, $J_3:=\left\{2,\dotsc,k\right\}\backslash \(J_1\cup J_2\)$ and $R'_1,R'_2>0$ be two numbers such that $R_1<R'_1<R'_2<R_2$. By using \eqref{Th1Eq2} and \eqref{Th1Eq7}, we then obtain 
\begin{align}\label{Th1Eq13}
\left|\tilde{u}_\alpha\(y\)-V_1\(y\)\right|&=\bigO\(\sum_{i=2}^k\(\frac{\mu_{1,\alpha}\mu_{i,\alpha}}{d_g\big(\exp_{\xi_{1,\alpha}}\(\mu_{1,\alpha}y\),\xi_{i,\alpha}\big)^2}\)^{\frac{n-2}{2}}\)+\smallo\(1\)\allowdisplaybreaks\nonumber\\
&=\bigO\Bigg(\sum_{i\in J_1}\(\frac{\mu_{1,\alpha}\mu_{i,\alpha}}{d_g\(\xi_{1,\alpha},\xi_{i,\alpha}\)^2}\)^{\frac{n-2}{2}}+\sum_{i\in J_2}\(\frac{\mu_{i,\alpha}}{\mu_{1,\alpha}}\)^{\frac{n-2}{2}}\(\left|y\right|-R_1\)^{2-n}\nonumber\\
&\qquad+\sum_{i\in J_3}\(\frac{\mu_{i,\alpha}}{\mu_{1,\alpha}}\)^{\frac{n-2}{2}}\(R_2-\left|y\right|\)^{2-n}\Bigg)+\smallo\(1\)=\smallo\(1\)
\end{align}
as $\alpha\to\infty$, uniformly with respect to $y\in B\(0,R'_2\)\backslash B\(0,R'_1\)$. By remarking that $\tilde{g}_\alpha\to\delta_0$ uniformly in $B\(0,R'_2\)$ as $\alpha\to\infty$, it follows from \eqref{Th1Eq13} that 
\begin{equation}\label{Th1Eq14}
\int_{B\(0,R'_2\)}\tilde{u}_\alpha^2\,dv_{\tilde{g}_\alpha}\ge \int_{B\(0,R'_2\)\backslash B\(0,R'_1\)}\tilde{u}_\alpha^2\,dv_{\tilde{g}_\alpha}\longrightarrow\int_{B\(0,R'_2\)\backslash B\(0,R'_1\)}V_1^2\,dv_{\delta_0}>0
\end{equation}
as $\alpha\to\infty$. Since $\mu_{1,\alpha}\to0$ as $\alpha\to\infty$, by putting together \eqref{Th1Eq11}, \eqref{Th1Eq12} and \eqref{Th1Eq14}, we then obtain a contradiction when $n\ge7$. This ends the proof of Theorem~\ref{Th1}.
\endproof

We now prove Theorem~\ref{Th5}, which shows that the conditions $\varphi_0\(\xi_0\)=0$ in Theorem~\ref{Th2} and $u_0\(\xi_0\)=0$ in Theorem~\ref{Th3} are necessary in the case of solutions of type \eqref{IntroEq3}. Note that in the statement of Theorem~\ref{Th5}, precisely in \eqref{IntroEq3}, $u_0$ can be any sign-changing solution to \eqref{IntroEq1}.

\smallskip
Let us first introduce a crucial ingredient, namely the conformal normal coordinates, which we use both in the proof of Theorem~\ref{Th5} and in our constructions of families of blowing-up solutions in the next sections. By using Lee and Parker's construction \cite{LP}, we obtain the existence of a smooth family of metrics $\(g_\xi\)_{\xi\in M}$, $g_\xi=\Lambda_\xi^{2^*-2}g$, such that for every point $\xi\in M$, the function $\Lambda_\xi$ is smooth, positive and satisfies
\begin{equation}\label{Sec3Eq1}
\Lambda_\xi\(\xi\)=1,\quad\nabla\Lambda_\xi\(\xi\)=0\quad\text{and}\quad dv_{g_\xi}\(x\)=\big(1+\bigO\big(\left|x\right|^N\big)\big)dv_{\delta_0}\(x\)
\end{equation}
in geodesic normal coordinates, where $dv_{g_\xi}$ and $dv_{\delta_0}$ are the volume elements of $g_\xi$ and the Euclidean metric $\delta_0$, respectively, and $N\in\N$ can be chosen arbitrarily large. In particular (see \cite{LP}), it follows from \eqref{Sec3Eq1} that
\begin{equation}\label{Sec3Eq2}
\Scal_{g_\xi}(\xi)=0,\quad\nabla \Scal_{g_\xi}(\xi)=0\quad\text{and}\quad\Delta_g \Scal_{g_\xi}(\xi)=\frac{1}{6}|\Weyl_g(\xi)|_g^2.
\end{equation}
In the case where $\(M,g\)$ is locally conformally flat, we may assume moreover that there exists a number $r_0>0$ such that for every point $\xi\in M$, the metric $g_\xi$ is flat in the ball of center $\xi$ and radius $r_0$.

\smallskip
As a first step in the proof of Theorem~\ref{Th5}, we obtain the following result, which is essentially contained in the article by Premoselli \cite{P}:

\begin{lemma}\label{Lem1}
Let $\(M,g\)$ be a smooth, closed Riemannian manifold of dimension $n\ge3$, $h_0\in C^{0,\vartheta}\(M\)$, $\vartheta\in\(0,1\)$, and $u_0$ be a solution to \eqref{IntroEq1}. Assume that there exist families of functions $\(h_\eps\)_{0<\eps\ll1}$\,in $C^{0,\vartheta}\(M\)$ and $\(u_\eps\)_{0<\eps\ll1}$ in $C^2\(M\)$, positive numbers $\(\mu_\eps\)_{0<\eps\ll1}$ and points $\(\xi_\eps\)_{0<\eps\ll1}$ in $M$ satisfying \eqref{IntroEq2} and \eqref{IntroEq3} and such that $\mu_\eps\to0$, $\xi_\eps\to\xi_0$ and $h_\eps\to h_0$ in $C^{0,\vartheta}\(M\)$ as $\eps\to0$. Let $\(g_\xi\)_\xi$, $g_\xi=\Lambda_\xi^{2^*-2}g$, be a family of conformal metrics to $g$ satisfying \eqref{Sec3Eq1}. Then 
\begin{equation}\label{Lem1Eq1}
u_\eps-u_0+B_\eps=\smallo\big(B_\eps+1\big)
\end{equation}
as $\eps\to0$, uniformly in $M$, where 
$$ B_\eps\(x\):= \Lambda_{\xi_\eps} \(x \)   \chi \( d_{g_{\xi_\eps}}\(\xi_{\eps}, x\) \)\(\frac{\sqrt{n\(n-2\)}\mu_\eps}{\mu_\eps^2+d_{g_{\xi_\eps}}\(\xi_\eps,x\)^2}\)^{\frac{n-2}{2}}\quad\forall x\in M $$
and  $\chi$ is a smooth, nonnegative cutoff function in $\[0,\infty\)$ such that $\chi=1$ in $\[0,r_0/2\)$ and $\chi=0$ in $\[r_0,\infty\)$ for some number $r_0>0$ smaller than the injectivity radii of the metrics $g_\xi$.
\end{lemma}

\proof[Proof of Lemma~\ref{Lem1}]
Since $u_\eps$ satisfies \eqref{IntroEq1} and \eqref{IntroEq3}, the main result in the article by Premoselli \cite{P} shows that 
\begin{equation}\label{Lem1Eq2}
u_\eps=u_0- \widetilde{B}_{\eps}+\smallo\(1+\(\frac{\mu_\eps}{\mu_\eps^2+d_{g}\(\xi_\eps,\cdot\)^2}\)^{\frac{n-2}{2}}\)\text{ as }\eps\to0,
\end{equation}
uniformly in $M$, where $\widetilde{B}_{\eps}$ is defined by 
$$ \widetilde{B}_\eps\(x\):=\chi \( d_{g}\(\xi_{\eps}, x\) \)\(\frac{\sqrt{n\(n-2\)}\mu_\eps}{\mu_\eps^2+d_{g}\(\xi_\eps,x\)^2}\)^{\frac{n-2}{2}}\quad\forall x\in M$$
if $u_0 \not \equiv 0$ or $\ker \(\Delta_g + h_0\) \neq \left\{0\right\}$, and 
\begin{multline*} 
 \widetilde{B}_\eps\(x\) :=  \chi \( d_{g}\(\xi_{\eps}, x\) \) F\(\xi_\eps, x \) \(\frac{\sqrt{n\(n-2\)}\mu_\eps}{\mu_\eps^2+d_{g}\(\xi_\eps,x\)^2}\)^{\frac{n-2}{2}} \\
 + \(1 -   \chi \( d_{g}(\xi_{\eps}, x) \) \)k_n\mu_{\eps}^{\frac{n-2}{2}} G_{g,h_0}\(\xi_\eps, x\)\quad\forall x\in M
\end{multline*} 
otherwise, where 
\begin{align}
&k_n:=n^{\frac{n-2}{4}}\(n-2\)^{\frac{n+2}{4}}\omega_{n-1},\label{kn}\allowdisplaybreaks\\
&F\(\xi_\eps, x\) :=  \(n-2\)\omega_{n-1} d_{g}\(\xi_\eps, x\)^{n-2} G_{g,h_0}\(\xi_{\eps}, x\) \quad\forall x\in M\nonumber
\end{align}
and $G_{g,h_0}$ is the Green's function of the operator $\Delta_g + h_0$. Since $\Lambda_{\xi_\eps}(\xi_\eps) = 1$, we have
$$ \Lambda_{\xi_\eps}(x) = 1 + \bigO \( d_g(\xi_\eps, x)\) $$
and
\begin{equation} \label{Lem1Eq3}
d_{g_{\xi_\eps}} \( x, \xi_\eps \)^2 = d_{g} \( x, \xi_\eps \)^2 + \bigO \( d_{g} \( x, \xi_\eps \)^3\) 
\end{equation}
uniformly with respect to $\eps >0 $ and $x \in M$. Independently, by definition of $G_{g,h_0}$, we have
$$ F\(\xi_\eps, x\) = 1 + \bigO\( d_g\(x, \xi_\eps\) \). $$
This shows that 
\begin{equation} \label{Lem1Eq4}
\begin{aligned}
\widetilde{B}_\eps(x) & = \big( 1+ \bigO\(d_g\(x, \xi_{\eps}\)\) \big) B_\eps(x) =  B_\eps(x)  + \smallo\( B_\eps(x) + 1 \)
\end{aligned}
\end{equation}
as $\eps\to0$, uniformly with respect to $x \in M$. The estimate \eqref{Lem1Eq1} then follows from \eqref{Lem1Eq2}, \eqref{Lem1Eq3} and \eqref{Lem1Eq4}.
\endproof

By using Lemma~\ref{Lem1}, we then obtain the following:

\begin{lemma}\label{Lem2}
Let $\(M,g\)$ be a smooth, closed Riemannian manifold of dimension $n\ge3$, $h_0\in C^{0,\vartheta}\(M\)$, $\vartheta\in\(0,1\)$, and $u_0$ be a solution to \eqref{IntroEq1}. Let $\(h_\eps\)_\eps$, $\(u_\eps\)_\eps$, $\(\mu_\eps\)_\eps$, $\(\xi_\eps\)_\eps$, $\(B_\eps\)_\eps$ and $\(g_\xi\)_\xi$ be as in Lemma~\ref{Lem1}. Then
\begin{multline}\label{Lem2Eq1}
\int_{B\(0,1/\sqrt{\mu_\eps}\)}\(\<\nabla\hat{u}_\eps,\cdot\>_{\delta_0}+\frac{n-2}{2}\hat{u}_\eps\)\(\(\Delta_{\hat{g}_\eps}-\Delta_{\delta_0}\)\hat{u}_\eps+\mu_\eps^2\hat{h}_\eps\hat{u}_\eps\)dv_{\delta_0}\\
=\frac{1}{2}n^{\frac{n-2}{4}}\(n-2\)^{\frac{n+6}{4}}\omega_{n-1}u_0\(\xi_0\)\mu_\eps^{\frac{n-2}{2}}+\smallo\(\mu_\eps^{\frac{n-2}{2}}\)
\end{multline}
as $\eps\to0$, where 
\begin{equation}\label{Lem2Eq2}
\hat{u}_\eps\(y\):=\mu_\eps^{\frac{n-2}{2}}\big(\Lambda_{\xi_\eps}^{-1}u_\eps\big)\(\exp_{\xi_\eps}\(\mu_\eps y\)\),\quad\hat{g}_\eps\(y\):=\exp_{\xi_\eps}^*g_{\xi_\eps}\(\mu_\eps y\)
\end{equation}
and
\begin{equation}\label{Lem2Eq3}
\hat{h}_\eps\(y\):=\big(c_n\Scal_{\hat{g}_\eps}+\Lambda_{\xi_\eps}^{2-2^*}\(h_\eps-c_n\Scal_g\)\big)\(\exp_{\xi_\eps}\(\mu_\eps y\)\)
\end{equation}
for all points $y\in B\(0,1/\sqrt{\mu_\eps}\)$, $\exp_{\xi_\eps}$ is the exponential map at the point $\xi_\eps$ with respect to the metric $g_{\xi_\eps}$ and we identify $T_{\xi_\eps}M$ with $\R^n$.
\end{lemma}

\proof[Proof of Lemma~\ref{Lem2}]
By using the conformal invariance property of the conformal Laplacian together with the change of variable $x=\exp_{\xi_\eps}\(\mu_\eps y\)$ for $y\in B\(0,r_0/\mu_\eps\)$, where $r_0$ is as in Lemma~\ref{Lem1}, the equation \eqref{IntroEq2} can be rewritten as
$$\Delta_{\hat{g}_\eps}\hat{u}_\eps+\mu_\eps^2\hat{h}_\eps\hat{u}_\eps=\left|\hat{u}_\eps\right|^{2^*-2}\hat{u}_\eps\quad\text{in }B\(0,r_0/\mu_\eps\),$$
where $\hat{u}_\eps$ and $\hat{h}_\eps$ are as in \eqref{Lem2Eq2} and \eqref{Lem2Eq3}. By using a Pohozaev-type identity (see for example the article by Marques \cite{M}*{formula~(2.7)}), we then obtain
\begin{multline}\label{Lem2Eq4}
\int_{B\(0,1/\sqrt{\mu_\eps}\)}\(\<\nabla\hat{u}_\eps,\cdot\>_{\delta_0}+\frac{n-2}{2}\hat{u}_\eps\)\(\(\Delta_{\hat{g}_\eps}-\Delta_{\delta_0}\)\hat{u}_\eps+\mu_\eps^2\hat{h}_\eps\hat{u}_\eps\)dv_{\delta_0}=\frac{1}{\sqrt{\mu_\eps}}\\
\times\int_{\partial B\(0,1/\sqrt{\mu_\eps}\)}\(\frac{n-2}{2}\sqrt{\mu_\eps}\hat{u}_\eps\partial_\nu\hat{u}_\eps+\(\partial_\nu\hat{u}_\eps\)^2-\frac{1}{2}\left|\nabla\hat{u}_\eps\right|_{\delta_0}^2+\frac{1}{2^*}\left|\hat{u}_\eps\right|^{2^*}\)d\sigma_{\delta_0},
\end{multline}
where $\nu$ and $d\sigma_{\delta_0}$ are the outward unit normal vector and volume element, respectively, of the metric induced by $\delta_0$ on $\partial B\(0,1/\sqrt{\mu_\eps}\)$. Similarly as in \eqref{Th1Eq5}-\eqref{Th1Eq9}, to estimate the terms in the right-hand side of \eqref{Lem2Eq4}, we use another change of variable, which gives
\begin{equation}\label{Lem2Eq5}
\Delta_{\check{g}_\eps}\check{u}_\eps+\mu_\eps\check{h}_\eps\check{u}_\eps=\mu_\eps\left|\check{u}_\eps\right|^{2^*-2}\check{u}_\eps\quad\text{in }B\(0,r_0/\sqrt{\mu_\eps}\)
\end{equation}
and
\begin{multline}\label{Lem2Eq6}
\int_{\partial B\(0,1/\sqrt{\mu_\eps}\)}\(\frac{n-2}{2}\sqrt{\mu_\eps}\hat{u}_\eps\partial_\nu\hat{u}_\eps+\(\partial_\nu\hat{u}_\eps\)^2-\frac{1}{2}\left|\nabla\hat{u}_\eps\right|_{\delta_0}^2+\frac{1}{2^*}\left|\hat{u}_\eps\right|^{2^*}\)d\sigma_{\delta_0}\\
=\mu_\eps^{\frac{n-1}{2}}\int_{\partial B\(0,1\)}\(\frac{n-2}{2}\check{u}_\eps\partial_\nu\check{u}_\eps+\(\partial_\nu\check{u}_\eps\)^2-\frac{1}{2}\left|\nabla\check{u}_\eps\right|_{\check{g}_\eps}^2+\frac{\mu_\eps}{2^*}\left|\check{u}_\eps\right|^{2^*}\)d\sigma_{\delta_0},
\end{multline}
where
$$\check{u}_\eps\(y\):=\big(\Lambda_{\xi_\eps}^{-1}u_\eps\big)\(\exp_{\xi_\eps}\(\sqrt{\mu_\eps}\,y\)\)\quad\text{and}\quad\check{g}_\eps(y) = \exp_{\xi_\eps}^*g_{\xi_\eps}\(\sqrt{\mu_\eps}\,y\).$$
By standard elliptic theory and since $\Lambda_{\xi_\eps}\(\xi_\eps\)=1$, it follows from \eqref{Lem1Eq1} and \eqref{Lem2Eq5} that 
\begin{equation}\label{Lem2Eq7}
\check{u}_\eps\to u_0\(\xi_0\)-\(n\(n-2\)\)^{\frac{n-2}{4}}{\left|\cdot\right|^{2-n}}\quad\text{in }C^1_{\loc}\(\R^n\backslash\left\{0\right\}\)
\end{equation}
as $\eps\to0$. It then follows from \eqref{Lem2Eq7} that
\begin{multline}\label{Lem2Eq8}
\int_{\partial B\(0,1\)}\(\frac{n-2}{2}\check{u}_\eps\partial_\nu\check{u}_\eps+\(\partial_\nu\check{u}_\eps\)^2-\frac{1}{2}\left|\nabla\check{u}_\eps\right|_{\check{g}_\eps}^2+\frac{\mu_\eps}{2^*}\left|\check{u}_\eps\right|^{2^*}\)d\sigma_{\delta_0}\\
=\frac{1}{2}n^{\frac{n-2}{4}}\(n-2\)^{\frac{n+6}{4}}\omega_{n-1}u_0\(\xi_0\)+\smallo\(1\)
\end{multline}
as $\eps\to0$. The estimate \eqref{Lem2Eq1} then follows from \eqref{Lem2Eq4}, \eqref{Lem2Eq6} and \eqref{Lem2Eq8}.
\endproof

We can now use Lemmas~\ref{Lem1} and~\ref{Lem2} to prove Theorem~\ref{Th5}.
\proof[Proof of Theorem~\ref{Th5}]
We proceed by contrapositive. Assume that there exist families of functions $\(h_\eps\)_{0<\eps\ll1}$ in $C^{0,\vartheta}\(M\)$ and $\(u_\eps\)_{0<\eps\ll1}$ in $C^2\(M\)$, positive numbers $\(\mu_\eps\)_{0<\eps\ll1}$ and points $\(\xi_\eps\)_{0<\eps\ll1}$ in $M$ satisfying \eqref{IntroEq2} and \eqref{IntroEq3} and such that $\mu_\eps\to0$, $\xi_\eps\to\xi_0$ and $h_\eps\to h_0$ in $C^{0,\vartheta}\(M\)$ as $\eps\to0$. Let $\(B_\eps\)_\eps$ and $\(g_\xi\)_\xi$ be as in Lemma~\ref{Lem1}. We now estimate the terms in the left-hand side of \eqref{Lem2Eq1}. It follows from \eqref{Lem1Eq1} together with standard elliptic estimates that 
\begin{equation}\label{Th5Eq2}
\left|\nabla^i\(\hat{u}_\eps+U_0\)\(y\)\right|=\smallo\big(\(1+\left|y\right|\)^{2-n-i}\big)\quad\text{for }i\in\left\{0,1,2\right\}
\end{equation}
as $\eps\to0$, uniformly with respect to $y\in B\(0,1/\sqrt{\mu_\eps}\)${, where
\begin{equation}\label{Th5Eq3}
U_0\(y\):=\(\frac{\sqrt{n\(n-2\)}}{1+\left|y\right|^2}\)^{\frac{n-2}{2}}.
\end{equation}}Since the function $U_0$ is radially symmetric about the origin, it follows from \eqref{Sec3Eq1} that 
\begin{equation}\label{Th5Eq4}
\(\Delta_{\hat{g}_\eps}-\Delta_{\delta_0}\)U_0\(y\)=\bigO\big(\mu_\eps^{N}\left|y\right|^{N-1}\left|\nabla U_0\(y\)\right|\big)=\bigO\(\frac{\mu_\eps^{N}\left|y\right|^{N}}{\(1+\left|y\right|\)^n}\)
\end{equation}
uniformly with respect to $y\in B\(0,1/\sqrt\mu_\eps\)$, where $\left|y\right|:=\left|y\right|_{\delta_0}$. By using \eqref{Th5Eq2} and \eqref{Th5Eq4}, we then obtain
\begin{multline}\label{Th5Eq5}
\(\Delta_{\hat{g}_\eps}-\Delta_{\delta_0}\)\hat{u}_\eps\(y\)=\bigO\Bigg(\frac{\mu_\eps^{N-1}\left|y\right|^{N}}{\(1+\left|y\right|\)^n}+\left|\(\hat{g}_\eps-\delta_0\)\(y\)\right||\nabla^2\(\hat{u}_\eps+U_0\)\(y\)|\\
+\left|\nabla\hat{g}_\eps\(y\)\right|\left|\nabla\(\hat{u}_\eps+U_0\)\(y\)\right|\Bigg)=\smallo\(\frac{\mu_\eps^2\left|y\right|}{\(1+\left|y\right|\)^{n-1}}\)
\end{multline}
as $\eps\to0$, uniformly with respect to $y\in B\(0,1/\sqrt{\mu_\eps}\)$. Since $h_\eps\to h_0$ in $C^{0,\vartheta}\(M\)$ as $\eps\to0$, by using \eqref{Sec3Eq2} and \eqref{Th5Eq2}, we obtain 
\begin{equation}\label{Th5Eq6}
\widehat{h}_\eps\(y\)\widehat{u}_\eps\(y\)=\(c_n\Scal_g-h_0\)\(\xi_0\)U_0\(y\)+\smallo\(U_0\(y\)\)
\end{equation}
uniformly with respect to $y\in B\(0,1/\sqrt\mu_\eps\)$. It follows from \eqref{Th5Eq2}, \eqref{Th5Eq5} and \eqref{Th5Eq6} that
\begin{multline}\label{Th5Eq7}
\int_{B\(0,1/\sqrt{\mu_\eps}\)}\(\<\nabla\hat{u}_\eps,\cdot\>_{\delta_0}+\frac{n-2}{2}\hat{u}_\eps\)\(\(\Delta_{\hat{g}_\eps}-\Delta_{\delta_0}\)\hat{u}_\eps+\mu_\eps^2\hat{h}_\eps\hat{u}_\eps\)dv_{\delta_0}\\
=\left\{\begin{aligned}&\bigO\(\mu_\eps^2\int_{B\(0,1/\sqrt{\mu_\eps}\)}\frac{dy}{\(1+\left|y\right|\)^{2n-4}}\)=\left\{\begin{aligned}&\bigO\(\mu_\eps^{3/2}\)&&\text{if }n=3\\&\bigO\(\mu_\eps^2\ln\(1/\mu_\eps\)\)&&\text{if }n=4\end{aligned}\right.\\&\frac{1}{2}n^{\frac{n-2}{2}}\(n-2\)^{\frac{n}{2}}a_n\(c_n\Scal_g-h_0\)\(\xi_0\)\mu_\eps^2+\smallo\(\mu_\eps^2\)\quad\text{if }n\ge5\end{aligned}\right.
\end{multline}
as $\eps\to0$, where
\begin{equation}\label{Th5Eq8}
a_n:=\int_{\R^n}\frac{\left|y\right|^2-1}{\big(1+\left|y\right|^2\big)^{n-1}}\,dy.
\end{equation} 
By putting together \eqref{Lem2Eq1} and \eqref{Th5Eq7}, we then obtain
\begin{multline*}
\omega_{n-1}u_0\(\xi_0\)\mu_\eps^{\frac{n-2}{2}}+\smallo\(\mu_\eps^{\frac{n-2}{2}}\)\\
=\left\{\begin{aligned}&\bigO\big(\mu_\eps^{3/2}\big)&&\text{if }n=3\\&\bigO\(\mu_\eps^2\ln\(1/\mu_\eps\)\)&&\text{if }n=4\\&n^{\frac{n-2}{4}}\(n-2\)^{\frac{n-6}{4}}a_n\(c_n\Scal_g-h_0\)\(\xi_0\)\mu_\eps^2+\smallo\(\mu_\eps^2\)&&\text{if }n\ge5\end{aligned}\right.
\end{multline*}
as $\eps\to0$. In the case where $n=6$, straightforward computations give
\begin{equation}\label{Th5Eq9}
n^{\frac{n-2}{4}}\(n-2\)^{\frac{n-6}{4}}a_n=\frac{\omega_5}{2}\,.
\end{equation} 
It then follows from \eqref{Th5Eq8} and \eqref{Th5Eq9} that $\varphi_0\(\xi_0\)=0$, where $\varphi_0$ is as in \eqref{Th5Eq1}. By contrapositive, this proves Theorem~\ref{Th5}.
\endproof

\section{Proofs of Theorem~\ref{Th2} and Corollaries~\ref{Cor1} and~\ref{Cor2}}\label{Sec4}

This section is devoted to our constructions of blowing-up solutions in the case of dimensions $n\ge6$. 

\proof[Proof of Theorem~\ref{Th2}]
The proof of Robert and V\'etois \cite{RV4}*{Theorem~1.5} continues to apply with straightforward adaptations in the case where $n\ge7$ or [$n=6$ and $u_0\(\xi_0\)<0$]. Therefore, in what follows, we only consider the case where $n=6$ and $u_0\(\xi_0\)>0$. We let $\(g_\xi\)_{\xi\in M}$, $g_\xi=\Lambda_\xi^{2^*-2}g$, be a family of conformal metrics to $g$ satisfying \eqref{Sec3Eq1}. We let $\chi$ be a smooth, nonnegative cutoff function in $\[0,\infty\)$ such that $\chi=1$ in $\[0,r_0/2\)$ and $\chi=0$ in $\[r_0,\infty\)$ for some number $r_0>0$ smaller than the injectivity radii of the metrics $g_\xi$. For every $\mu>0$ and $x,\xi\in M$, we then define 
$$B_{\mu,\xi}\(x\):=\Lambda_\xi\(x\)\chi\(d_{g_\xi}\(x,\xi\)\)\(\frac{2\sqrt{6}\,\mu}{\mu^2+d_{g_\xi}\(x,\xi\)^2}\)^2.$$
Since $u_0$ is a nondegenerate solution to \eqref{IntroEq1}, we obtain the existence of a family of solutions $\(w_\eps\)_{0<\eps\ll1}$ in $C^2\(M\)$ to the equations
\begin{equation}\label{Th2Eq1}
\Delta_gw_\eps+\(h_0+\eps h\)w_\eps=\left|w_\eps\right|w_\eps\quad\text{in }M
\end{equation}
such that $w_\eps\to u_0$ in $C^2\(M\)$ as $\eps\to0$ and $w_\eps$ is differentiable with respect to $\eps$. By differentiating \eqref{Th2Eq1}, we obtain 
\begin{equation}\label{Th2Eq2}
\left.\frac{dw_\eps}{d\eps}\right|_{\eps=0}=-\big(\Delta_g+h_0-2\left|u_0\right|\big)^{-1}\(hu_0\).
\end{equation}
We then define
$$J_\eps\(u\):=\frac{1}{2}\int_M\(\left|\nabla u\right|_g^2+\(h_0+\eps h\)u^2\)dv_g-\frac{1}{3}\int_M\left|u\right|^3dv_g$$
for all $\eps>0$ and $u\in H^1\(M\)$. By using \eqref{Th2Eq1} together with straightforward computations, we obtain
\begin{multline}\label{Th2Eq3}
J_\eps\(w_\eps-B_{\mu,\xi}\)=\frac{1}{6}\int_{\R^n}\left|w_\eps\right|^3\,dv_g+J_\eps\(B_{\mu,\xi}\)+\int_M w_\eps B_{\mu,\xi}^2\,dv_g\\
-\frac{2}{3}\int_{A_{\eps,\mu,\xi}}\(3w_\eps-B_{\mu,\xi}\)B_{\mu,\xi}^2\,dv_g-\frac{2}{3}\int_{A'_{\eps,\mu,\xi}}\(3B_{\mu,\xi}-w_\eps\)w_\eps^2\,dv_g,
\end{multline}
where 
$$A_{\eps,\mu,\xi}:=\left\{x\in M:\,w_\eps\(x\)\ge B_{\mu,\xi}\(x\)\right\}$$ 
and 
$$A'_{\eps,\mu,\xi}:=\left\{x\in M:\,0<w_\eps\(x\)<B_{\mu,\xi}\(x\)\right\}.$$
By applying \cite{RV4}*{Proposition~7.2}, we obtain
\begin{equation}\label{Th2Eq4}
J_\eps\(B_{\mu,\xi}\)+\int_M w_\eps B_{\mu,\xi}^2\,dv_g=\frac{1}{6}\int_{\R^n}U_0^3\,dv_{\delta_0}+\frac{1}{2}\varphi_\eps\(\xi\)\mu^2\int_{\R^n}U_0^2\,dv_{\delta_0}+\smallo\(\mu^3\)
\end{equation}
as $\eps,\mu\to0$, uniformly with respect to $\xi\in M$, where $\delta_0$ is the Euclidean metric, {$U_0$ is as in \eqref{Th5Eq3} and
$$\varphi_\eps:=h_0-c_n\Scal_g+\eps h + 2 w_\eps.$$}Further estimates yield
\begin{equation}\label{Th2Eq5}
\int_{A_{\eps,\mu,\xi}}\(3w_\eps-B_{\mu,\xi}\)B_{\mu,\xi}^2\,dv_g=64\sqrt{6}\,\omega_5u_0\(\xi_0\)^{3/2}\mu^3+\smallo\(\mu^3\)
\end{equation}
and
\begin{equation}\label{Th2Eq6}
\int_{A'_{\eps,\mu,\xi}}\(3B_{\mu,\xi}-w_\eps\)w_\eps^2\,dv_g=64\sqrt{6}\,\omega_5u_0\(\xi_0\)^{3/2}\mu^3+\smallo\(\mu^3\)
\end{equation}
as $\eps,\mu\to0$, uniformly with respect to $\xi\in M$, where $\omega_5$ is the volume of the standard unit sphere of dimension $5$. It follows from \eqref{Th2Eq3}--\eqref{Th2Eq6} that
\begin{multline}\label{Th2Eq7}
J_\eps\(w_\eps-B_{\mu,\xi}\)=\frac{1}{6}\int_{\R^n}\left|w_\eps\right|^3\,dv_g+\frac{1}{6}\int_{\R^n}U_0^3\,dv_{\delta_0}+\frac{1}{2}\varphi_\eps\(\xi\)\mu^2\int_{\R^n}U_0^2\,dv_{\delta_0}\\
-\frac{256}{3}\sqrt{6}\,\omega_5u_0\(\xi_0\)^{3/2}\mu^3+\smallo\(\mu^3\)
\end{multline}
as $\eps,\mu\to0$, uniformly with respect to $\xi\in M$. On the other hand, by using \eqref{Th2Eq1}, we obtain
\begin{multline}\label{Th2Eq8}
\left\|\(\Delta_g+h_0+\eps h-\left|w_\eps-B_{\mu,\xi}\right|\)\(w_\eps-B_{\mu,\xi}\)\right\|_{L^{3/2}\(M\)}\\
=\left\|\(\Delta_g+h_0+\eps h+2w_\eps-B_{\mu,\xi}\)B_{\mu,\xi}\right\|_{L^{3/2}\(M\)}\\
+2\left\|\(2w_\eps-B_{\mu,\xi}\)B_{\mu,\xi}\right\|_{L^{3/2}\(A_{\eps,\mu,\xi}\)}+2\left\|w_\eps^2\right\|_{L^{3/2}\(A'_{\eps,\mu,\xi}\)}.
\end{multline}
By applying \cite{RV4}*{Proposition~7.1}, we obtain
\begin{equation}\label{Th2Eq9}
\left\|\(\Delta_g+h_0+\eps h+2w_\eps-B_{\mu,\xi}\)B_{\mu,\xi}\right\|_{L^{3/2}\(M\)}=\bigO\(\mu^2\ln\(1/\mu\)^{2/3}\)
\end{equation}
uniformly with respect to $\eps,\mu>0$ and $\xi\in M$. Moreover, we obtain
\begin{equation}\label{Th2Eq10}
\left\|\(2w_\eps-B_{\mu,\xi}\)B_{\mu,\xi}\right\|_{L^{3/2}\(A_{\eps,\mu,\xi}\)}=\bigO\(\left\|B_{\mu,\xi}\right\|_{L^{3/2}\(A_{\eps,\mu,\xi}\)}\)=\bigO\(\mu^2\ln\(1/\mu\)^{2/3}\)
\end{equation}
and
\begin{equation}\label{Th2Eq11}
\left\|w_\eps^2\right\|_{L^{3/2}(A'_{\eps,\mu,\xi})}=\bigO\(\left\|1\right\|_{L^{3/2}(A'_{\eps,\mu,\xi})}\)=\bigO\(\mu^3\)
\end{equation}
uniformly with respect to $\eps,\mu>0$ and $\xi\in M$. It follows from \eqref{Th2Eq8}--\eqref{Th2Eq11} that
\begin{equation}\label{Th2Eq12}
\left\|\(\Delta_g+h_0+\eps h-\left|w_\eps-B_{\mu,\xi}\right|\)\(w_\eps-B_{\mu,\xi}\)\right\|_{L^{3/2}\(M\)}=\bigO\(\mu^2\ln\(1/\mu\)^{2/3}\)
\end{equation}
uniformly with respect to $\eps,\mu>0$ and $\xi\in M$. By applying \cite{RV2}*{Theorem~1.1}, it follows from \eqref{Th2Eq12} that for small $\eps,\mu>0$ and all $\xi\in M$, there exists a function $\psi_{\eps,\mu,\xi}\in H^1\(M\)$ such that $\psi_{\eps,\mu,\xi}\to0$ in $H^1\(M\)$ and 
\begin{equation}\label{Th2Eq13}
J_\eps\(w_\eps-B_{\mu,\xi}+\psi_{\eps,\mu,\xi}\)=J_\eps\(w_\eps-B_{\mu,\xi}\)+\smallo\(\mu^3\)
\end{equation}
as $\eps,\mu\to0$, uniformly with respect to $\xi\in M$, and 
\begin{equation}\label{Th2Eq14}
D_uJ_\eps\(w_\eps-B_{\mu,\xi}+\psi_{\eps,\mu,\xi}\)=0\ \ \ \ \Longleftrightarrow\ \ \ \ D_{\(\mu,\xi\)}\[J_\eps\(u_0-B_{\mu,\xi}+\psi_{\eps,\mu,\xi}\)\]=0.
\end{equation} 
For every $\eps,t>0$ and $\tau\in\R^n$, we define
$$\mu_\eps\(t\):=t\eps\quad\text{and}\quad\xi_\eps\(\tau\):=\exp_{\xi_\eps}\(\sqrt\eps\tau\).$$
Since $\varphi_0\(\xi_0\)=\left|\nabla\varphi_0\(\xi_0\)\right|=0$, it then follows from \eqref{Th2Eq7} and \eqref{Th2Eq13} that
\begin{multline}\label{Th2Eq15}
\eps^{-3}\(J_\eps\(w_\eps-B_{\mu_\eps\(t\),\xi_\eps\(\tau\)}+\psi_{\eps,\mu_\eps\(t\),\xi_\eps\(\tau\)}\)-\frac{1}{6}\int_{\R^n}\left|w_\eps\right|^3\,dv_g-\frac{1}{6}\int_{\R^n}U_0^{3}dv_{\delta_0}\)\\
\longrightarrow F\(t,\tau\):=\(\frac{1}{4}D^2\varphi_0\(\xi_0\).\(\tau,\tau\)+\frac{1}{2}\(h+2\left.\frac{dw_\eps}{d\eps}\right|_{\eps=0}\)\(\xi_0\)\)t^2\int_{\R^n}U_0^2\,dv_{\delta_0}\\
-\frac{256}{3}\sqrt{6}\,\omega_5u_0\(\xi_0\)^{3/2}t^3
\end{multline}
as $\eps\to0$, uniformly with respect to $\xi\in M$ and $t$ in compact subsets of $\(0,\infty\)$. By using the assumptions (i) and (ii) together with \eqref{Th2Eq2}, we obtain that the function $F$ attains a strict maximum at $\(t_0,0\)$ for some number $t_0>0$. By using \eqref{Th2Eq15}, we then obtain that there exists a family of maximum points $\(\(t_\eps,\tau_\eps\)\)_{0<\eps\ll1}$ of the function $\(t,\tau\)\mapsto J_\eps\(w_\eps-B_{\mu_\eps\(t\),\xi_\eps\(\tau\)}+\psi_{\eps,\mu_\eps\(t\),\xi_\eps\(\tau\)}\)$ in $\[t_0/2,2t_0\]\times B\(0,1\)$ that converges to $\(t_0,0\)$ as $\eps \to 0$. It then follows from \eqref{Th2Eq14} that for small $\eps>0$, the function $u_\eps:=w_\eps-B_{\mu_\eps\(t_\eps\),\xi_\eps\(\tau_\eps\)}+\psi_{\eps,\mu_\eps\(t_\eps\),\xi_\eps\(\tau_\eps\)}$ is a solution to \eqref{IntroEq2}. It is easy to verify that $u_\eps$ satisfies \eqref{IntroEq3}. This ends the proof of Theorem~\ref{Th2}.
\endproof

\proof[Proof of Corollary~\ref{Cor1}]
Here again, the proof of Robert and V\'etois \cite{RV4}*{Theorem~1.3} continues to apply with straightforward adaptations in the case where $n\ge7$. Therefore, in what follows, we only consider the case where $n=6$. Remark that in this case, we can rewrite the equation satisfied by $u_0$ as
$$\Delta_gu_0+\(h_0+2u_0\)u_0=3u_0^2\quad\text{in }M.$$ 
By applying \cite{RV4}*{Lemma~13.1}, we obtain that there exist a neighborhood $\Omega_0$ of $\xi_0$ and families of functions $(\tilde{h}_\eps)_{0\le\eps\ll1}$ in $C^p\(M\)$ and $\(\tilde{u}_\eps\)_{0\le\eps\ll1}$ in $C^2\(M\)$ such that $\tilde{h}_\eps\to h_0$ in $C^p\(M\)$ and $\tilde{u}_\eps\to u_0$ in $C^2\(M\)$ as $\eps\to0$, and for every $\eps>0$, $\tilde{h}_\eps=h_0$ and $\tilde{u}_\eps=u_0$ in $\Omega_0$, and $\tilde{u}_\eps$ is a nondegenerate, positive solution to the equation
\begin{equation}\label{Cor1Eq1}
\Delta_g\tilde{u}_\eps+\big(\tilde{h}_\eps+2u_0\big)\tilde{u}_\eps=3\tilde{u}_\eps^2\quad\text{in }M.
\end{equation} 
For every $\eps>0$, since $\tilde{h}_\eps=h_0$ and $\tilde{u}_\eps=u_0$ in $\Omega_0$, $\varphi_0\(\xi_0\)=0$ and $\xi_0$ is a local maximum point of the function $\varphi_0$, a straightforward perturbation argument (together with a regularization argument in case $p=1$) yields the existence of a family of functions $(\tilde{h}_{\eps,\alpha})_{0\le\alpha\ll1}$ in $C^{\max\(p,2\)}\(M\)$ such that $\tilde{h}_{\eps,\alpha}\to\tilde{h}_\eps$ in $C^p\(M\)$ as $\alpha\to0$ and for every $\alpha>0$, the function $\widetilde{\varphi}_{\eps,\alpha}:=\tilde{h}_{\eps,\alpha}+2u_0-c_n\Scal_g$ satisfies $\widetilde{\varphi}_{\eps,\alpha}\(\xi_0\)=\left|\nabla\widetilde{\varphi}_{\eps,\alpha}\(\xi_0\)\right|=0$ and $D^2\widetilde{\varphi}_{\eps,\alpha}\(\xi_0\)<0$. Since $\tilde{u}_\eps$ is a nondegenerate, positive solution to \eqref{Cor1Eq1}, we then obtain that there exists a family of positive solutions $(\tilde{u}_{\eps,\alpha})_{0\le\alpha\ll1}$ to the equations
\begin{equation}\label{Cor1Eq2}
\Delta_g\tilde{u}_{\eps,\alpha}+\big(\tilde{h}_{\eps,\alpha}+2u_0\big)\tilde{u}_{\eps,\alpha}=3\tilde{u}_{\eps,\alpha}^2\quad\text{in }M
\end{equation} 
such that $\tilde{u}_{\eps,\alpha}\to\tilde{u}_\eps$ in $C^2\(M\)$ as $\alpha\to0$. Since $\tilde{h}_\eps\to h_0$ and $\tilde{h}_{\eps,\alpha}\to\tilde{h}_\eps$ in $C^p\(M\)$ as $\eps\to 0$ and $\alpha\to 0$, respectively, by standard elliptic theory, we obtain that $\tilde{u}_\eps\to u_0$ and $\tilde{u}_{\eps,\alpha}\to\tilde{u}_\eps$ in $C^{p+1}\(M\)$ as $\eps\to 0$ and $\alpha\to 0$, respectively. As is easily seen, the equation \eqref{Cor1Eq2} can be rewritten as
$$\Delta_g\tilde{u}_{\eps,\alpha}+\hat{h}_{\eps,\alpha}\tilde{u}_{\eps,\alpha}=\tilde{u}_{\eps,\alpha}^2\quad\text{in }M,\quad\text{where }\hat{h}_{\eps,\alpha}:=\tilde{h}_{\eps,\alpha}+2u_0-2\tilde{u}_{\eps,\alpha}\quad\text{in }M.$$
Since $\tilde{h}_{\eps,\alpha}\to\tilde{h}_\eps$ and $\tilde{u}_{\eps,\alpha}\to\tilde{u}_\eps$ in $C^p\(M\)$ as $\alpha\to 0$, we obtain that $\hat{h}_{\eps,\alpha}\to \tilde{h}_\eps$ in $C^p\(M\)$. By applying again \cite{RV4}*{Lemma~13.1}, we then obtain that for every $\eps,\alpha>0$, there exist a neighborhood $\Omega_{\eps,\alpha}$ of $\xi_0$ and families of functions $(\hat{h}_{\eps,\alpha,\beta})_{0\le\beta\ll1}$ in $C^{\max\(p,2\)}\(M\)$ and $\(\hat{u}_{\eps,\alpha,\beta}\)_{0\le\beta\ll1}$ in $C^2\(M\)$ such that $\hat{h}_{\eps,\alpha,\beta}\to\hat{h}_{\eps,\alpha}$ in $C^p\(M\)$ and $\hat{u}_{\eps,\alpha,\beta}\to\tilde{u}_{\eps,\alpha}$ in $C^2\(M\)$ as $\beta\to0$, and for every $\beta>0$, $\hat{h}_{\eps,\alpha,\beta}=\hat{h}_{\eps,\alpha}$ and $\hat{u}_{\eps,\alpha,\beta}=\tilde{u}_{\eps,\alpha}$ in $\Omega_{\eps,\alpha}$, and $\hat{u}_{\eps,\alpha,\beta}$ is a nondegenerate, positive solution to the equation
$$\Delta_g\hat{u}_{\eps,\alpha,\beta}+\hat{h}_{\eps,\alpha,\beta}\hat{u}_{\eps,\alpha,\beta}=\hat{u}_{\eps,\alpha,\beta}^2\quad\text{in }M.$$
Remark that
$$\hat{h}_{\eps,\alpha,\beta}+2\hat{u}_{\eps,\alpha}-c_n\Scal_g=\hat{h}_{\eps,\alpha}+2\tilde{u}_{\eps,\alpha}-c_n\Scal_g=\widetilde{\varphi}_{\eps,\alpha}\quad\text{in }\Omega_{\eps,\alpha}\cap\Omega_0.$$
Recall that $\widetilde{\varphi}_{\eps,\alpha}\(\xi_0\)=\left|\nabla\widetilde{\varphi}_{\eps,\alpha}\(\xi_0\)\right|=0$ and $D^2\widetilde{\varphi}_{\eps,\alpha}\(\xi_0\)<0$. We can now conclude by applying Theorem~\ref{Th2} together with a straightforward diagonal argument. This ends the proof of Corollary~\ref{Cor1}. 
\endproof

\proof[Proof of Corollary~\ref{Cor2}]
Remark that the equation
$$\Delta_gu_0+h_0u_0=u_0^{2^*-1}\quad\text{in }M$$
can be rewritten as 
\begin{equation}\label{Cor2Eq1}
\Delta_gu_0+\(c_n\Scal_g+\varphi_0\)u_0=\left\{\begin{aligned}&3u_0^2&&\text{if }n=6\\&u_0^{2^*-1}&&\text{if }n\ge7\end{aligned}\right\}\text{ in }M.
\end{equation}
The results by Aubin \cite{A}, Schoen \cite{Sc1} and Trudinger \cite{T} give that  on the one hand, there exists a positive solution $u_0$ to \eqref{Cor2Eq1} provided the operator $\Delta_g+c_n\Scal_g+\varphi_0$ is coercive in $H^1\(M\)$ and 
\begin{equation}\label{Cor2Eq2}
\inf_{u\in C^\infty\(M\)\backslash\left\{0\right\}}\frac{\displaystyle\int_M\(\left|\nabla u\right|_g^2+\(c_n\Scal_g+\varphi_0\)u^2\)dv_g}{\displaystyle\(\int_M\left|u\right|^{2^*}dv_g\)^{2/2^*}}<K_n^{-2},
\end{equation}
where $K_n$ is the best constant for the embedding $D^{1,2}\(\R^n\)\hookrightarrow  L^{2^*}\(\R^n\)$, and on the other hand, \eqref{Cor2Eq2} is satisfied when either $\varphi_0<0$ at some point in $M$ or $\varphi_0=0$ in $M$ and $\(M,g\)$ is not conformally diffeomorphic to the standard sphere. In the latter case, by remarking that
$$\int_M\varphi_0u^2dv_g\le \left\|\max\(\varphi_0,0\)\right\|_{L^1\(M\)}\left\|u\right\|_{C^0\(M\)}^2,$$
we obtain that \eqref{Cor2Eq2} is satisfied when $\left\|\max\(\varphi_0,0\)\right\|_{L^1\(M\)}<\eps_0$ for some constant $\eps_0>0$ depending only on $\(M,g\)$. Therefore, in both cases (i) and (ii), we obtain that there exist functions $h_0\in C^p\(M\)$ and $u_0\in C^2\(M\)$ such that $u_0$ is a positive solution to \eqref{IntroEq1} and $\varphi_0\(\xi_0\)=\left|\nabla\varphi_0\(\xi_0\)\right|=0$. We can then conclude by applying Corollary~\ref{Cor1}.
\endproof

\section{Proofs of Theorems~\ref{Th3} and~\ref{Th4} and Corollaries~\ref{Cor3} and~\ref{Cor4}}\label{Sec5}

This section is devoted to our constructions of blowing-up solutions of type \eqref{IntroEq3} in the case of dimensions $n\in\left\{3,4,5\right\}$. We begin with proving Theorem~\ref{Th3}.

\proof[Proof of Theorem~\ref{Th3}]
We begin with constructing a suitable ansatz. Since $u_0 \not \equiv 0$ is a nondegenerate solution to \eqref{IntroEq1}, we obtain the existence of a family of solutions $\(w_\eps\)_{0<\eps\ll1}$ in $C^2\(M\)$ to the equations
\begin{equation}\label{Th3Eq2}
\Delta_gw_\eps+h_\eps w_\eps=\left|w_\eps\right|^{2^*-2}w_\eps\quad\text{in }M
\end{equation}
such that $w_\eps\to u_0$ in $C^2\(M\)$ as $\eps\to0$ and $w_\eps$ is differentiable with respect to $\eps$. Moreover, since $\(M,g\)$, $h_0$ and $u_0$ are symmetric with respect to $\xi_0$, we can assume that for each $\eps$, the function $w_\eps$ is also symmetric with respect to $\xi_0$. By differentiating \eqref{Th3Eq2}, we obtain 
\begin{equation}\label{Th3Eq3}
\left.\frac{dw_\eps}{d\eps}\right|_{\eps=0}=-\big(\Delta_g+\tilde{h}_0\big)^{-1}\(hu_0\),
\end{equation}
where $\tilde{h}_0 = h_0 - \(2^*-1\)\left|u_0\right|^{2^*-2}$. We let $g_{\xi_0}:=\Lambda_{\xi_0}^{2^*-2}g$ be a conformal metric to $g$ satisfying \eqref{Sec3Eq1} and $\chi$ be a smooth, nonnegative cutoff function in $\[0,\infty\)$ such that $\chi=1$ in $\[0,r_0/2\)$ and $\chi=0$ in $\[r_0,\infty\)$ for some number $r_0>0$ smaller than the injectivity radius of the metric $g_{\xi_0}$. For every $\mu>0$ and $x\in M$, we then define 
$$U_{\mu}\(x\):=\left\{\begin{aligned}&B_{\mu}\(x\)+k_n\mu^{\frac{n-2}{2}}\big(\chi\(d_{g_{\xi_0}}\(x,\xi_0\)\)H_{g,\tilde{h}_0}\(x,\xi_0\)&&\\&\quad+\(1-\chi\(d_{g_{\xi_0}}\(x,\xi_0\)\)\)G_{g,\tilde{h}_0}\(x,\xi_0\)\big)&&\text{in Case  1}\\&B_{\mu}\(x\)&&\text{in Case 2},\end{aligned}\right.$$
where $k_n$ is as in \eqref{kn},
\begin{align}
&\left\{\begin{aligned}&\text{Case 1: }n=3\text{ or } \[n=4,5 \text{ and }  h\(\xi_0\)=c_n\Scal_g\(\xi_0\)\]\\&\text{Case 2: }n\in\left\{4,5\right\}\text{ and }h\(\xi_0\)\ne c_n\Scal_g\(\xi_0\),\end{aligned}\right.\allowdisplaybreaks\nonumber\\
&B_{\mu}\(x\):=\Lambda_{\xi_0}\(x\)\chi\(d_{g_{\xi_0}}\(x,\xi_0\)\)\(\frac{\sqrt{n\(n-2\)}\mu}{\mu^2+d_{g_{\xi_0}}\(x,\xi_0\)^2}\)^{\frac{n-2}{2}}\label{Th3Eq4}
\end{align}
and
$$H_{g,\tilde{h}_0}\(x,\xi_0\):=G_{g,\tilde{h}_0}\(x,\xi_0\)-\(\(n-2\)\omega_{n-1}\)^{-1}\Lambda_{\xi_0}\(x\)d_{g_{\xi_0}}\(x,\xi_0\)^{2-n}.$$
By \eqref{IntroEq6}, the mass $m_{g,\tilde{h}_0}(\xi_0)$ is equal to $H_{g,\tilde{h}_0}(\xi_0,\xi_0)$. For every $\eps>0$ and $u\in H^1\(M\)$, we then define 
$$J_\eps\(u\):=\frac{1}{2}\int_M\(\left|\nabla u\right|_g^2+h_\eps u^2\)dv_g-\frac{1}{2^*}\int_M\left|u\right|^{2^*}dv_g.$$
To apply a Lyapunov--Schmidt reduction method based on the ansatz $w_\eps-U_\mu$ (see for example the article by Robert and V\'etois \cite{RV2}*{Theorem~1.1}), we need to estimate the asymptotic behavior as $\eps,\mu\to0$ of $J_\eps\(w_\eps-U_\mu\)$ together with an error term of the form
\begin{equation}\label{Th3Eq5}
R_{\eps,\mu}:=\left\|w_\eps-U_{\mu}-\(\Delta_g+1\)^{-1}\(\(\left|w_\eps-U_{\mu}\right|^{2^*-2}+1-h_\eps\)\(w_\eps-U_{\mu}\)\)\right\|_{H^1\(M\)}.
\end{equation}
We begin with estimating $J_\eps\(w_\eps-U_\mu\)$. By integrating by parts and using \eqref{Th3Eq2} together with straightforward estimates, we obtain 
\begin{multline}\label{Th3Eq6}
J_\eps\(w_\eps-U_{\mu}\)=\frac{1}{n}\int_M\left|w_\eps\right|^{2^*}dv_g+\frac{1}{2}\int_M\(\left|\nabla U_{\mu}\right|_g^2+\tilde{h}_\eps U_{\mu}^2\)dv_g-\frac{1}{2^*}\int_M\left|U_{\mu}\right|^{2^*}dv_g\\
+\int_M\left|U_{\mu}\right|^{2^*-2}U_{\mu}w_\eps\,dv_g+\bigO\Bigg(\int_{M_{\eps,\mu}}\(\left|U_{\mu}\right|^{2^*-2}w_\eps^2+\left|w_\eps\right|^{2^*-2}U_{\mu}^2\)\,dv_g\\
+\int_{M\backslash M_{\eps,\mu}}\(\left|U_{\mu}\right|^{2^*-1}\left|w_\eps\right|+\left|w_\eps\right|^{2^*-3}\left|U_{\mu}\right|^3\)dv_g\Bigg)
\end{multline}
uniformly with respect to $\eps,\mu>0$, where 
$$\tilde{h}_\eps:=h_\eps-\(2^*-1\)\left|w_\eps\right|^{2^*-2}\quad\text{and}\quad M_{\eps,\mu}:=\left\{x\in M:\,\left|w_\eps\(x\)\right|\le\left|U_{\mu}\(x\)\right|\right\}.$$ 
Remark that 
\begin{equation}\label{Th3Eq7}
U_{\mu}=\bigO\(\mu^{\frac{n-2}{2}}\)\,\text{and}\,\big(\Delta_g+\tilde{h}_\eps\big)U_{\mu}=k_n\mu^{\frac{n-2}{2}}\big(\tilde{h}_\eps-\tilde{h}_0\big)G_{g,\tilde{h}_0}\(\cdot,\xi_0\)=\smallo\(\mu^{\frac{n-2}{2}}\)
\end{equation}
as $\eps,\mu\to0$, uniformly in $M\backslash B_{g_{\xi_0}}\(\xi_0,r_0\)$, where $B_{g_{\xi_0}}\(\xi_0,r_0\)$ is the ball of center $\xi_0$ and radius $r_0$ with respect to the metric $g_{\xi_0}$. By identifying $T_{\xi_0}M$ with $\R^n$ and using \eqref{Th3Eq7} together with the conformal invariance property of the conformal Laplacian and the change of variable $x=\exp_{\xi_0}\(\mu y\)$ for $y\in B\(0,r_0/\mu\)$, where $\exp_{\xi_0}$ is the exponential map at the point $\xi_0$ with respect to the metric $g_{\xi_0}$ and $B\(0,r_0/\mu\)$ is the ball of center 0 and radius $r_0/\mu$ with respect to the Euclidean metric, we obtain
\begin{align}
&\int_M\(\left|\nabla U_{\mu}\right|_g^2+\tilde{h}_\eps U_{\mu}^2\)dv_g=\int_{B\(0,r_0/\mu\)}\widehat{U}_{\mu}\big(\Delta_{\hat{g}_\mu}\widehat{U}_{\mu}+\mu^2\hat{h}_{\eps,\mu}\widehat{U}_{\mu}\big)dv_{\hat{g}_\mu}+\smallo\(\mu^{n-2}\),\label{Th3Eq8}\allowdisplaybreaks\\
&\int_M\left|U_{\mu}\right|^{2^*}dv_g=\int_{B\(0,r_0/\mu\)}\big|\widehat{U}_{\mu}\big|^{2^*}dv_{\hat{g}_\mu}+\smallo\(\mu^{n-2}\),\label{Th3Eq9}\allowdisplaybreaks\\
&\int_M\left|U_{\mu}\right|^{2^*-2}U_{\mu}w_\eps\,dv_g=\int_{B\(0,r_0/\mu\)}\big|\widehat{U}_{\mu}\big|^{2^*-2}\widehat{U}_{\mu}\widehat{w}_{\eps,\mu}\,dv_{\hat{g}_\mu}+\smallo\(\mu^{n-2}\),\allowdisplaybreaks\label{Th3Eq10}\\
&\int_{M_{\eps,\mu}}\(\left|U_{\mu}\right|^{2^*-2}w_\eps^2+\left|w_\eps\right|^{2^*-2}U_{\mu}^2\)\,dv_g\nonumber\\
&\quad=\int_{\Omega_{\eps,\mu}}\(\big|\widehat{U}_{\mu}\big|^{2^*-2}\widehat{w}_{\eps,\mu}^2+\left|\widehat{w}_{\eps,\mu}\right|^{2^*-2}\widehat{U}_{\mu}^2\)dv_{\hat{g}_\mu}+\smallo\(\mu^{n-2}\)\label{Th3Eq11}
\end{align}
and
\begin{multline}\label{Th3Eq12}
\int_{M\backslash M_{\eps,\mu}}\(\left|U_{\mu}\right|^{2^*-1}\left|w_\eps\right|+\left|w_\eps\right|^{2^*-3}\left|U_{\mu}\right|^3\)dv_g\\
=\int_{B\(0,r_0/\mu\)\backslash \Omega_{\eps,\mu}}\(\big|\widehat{U}_{\mu}\big|^{2^*-1}\left|\widehat{w}_{\eps,\mu}\right|+\left|\widehat{w}_{\eps,\mu}\right|^{2^*-3}\big|\widehat{U}_{\mu}\big|^3\)dv_{\hat{g}_\mu}+\smallo\(\mu^{n-2}\)
\end{multline}
as $\mu\to0$, uniformly with respect to $\eps>0$, where
\begin{align*}
&\Omega_{\eps,\mu}:=\left\{y\in B\(0,r_0/\mu\):\,\left|\hat{w}_{\eps,\mu}\(y\)\right|\le\big|\widehat{U}_{\mu}\(y\)\big|\right\},\allowdisplaybreaks\\
&\hat{g}_\mu\(y\):=\exp_{\xi_0}^*g_{\xi_0}\(\mu y\),\quad\widehat{w}_{\eps,\mu}:=\mu^{\frac{n-2}{2}}w_\eps\big(\exp_{\xi_0}\(\mu y\)\big),\allowdisplaybreaks\\
&\hat{h}_{\eps,\mu}:=\(c_n\Scal_{g_{\xi_0}}+\Lambda_{\xi_0}^{2-2^*}\big(\tilde{h}_\eps-c_n\Scal_g\big)\)\big(\exp_{\xi_0}\(\mu y\)\big)
\end{align*}
and
\begin{multline*}
\widehat{U}_{\mu}\(y\):=\\
\left\{\begin{aligned}&\chi\(\mu\left|y\right|\)U_0\(y\)+k_n\mu^{n-2}\Lambda_{\xi_0}\(\exp_{\xi_0}\(\mu y\)\)^{-1}\Big(\chi\(\mu\left|y\right|\)H_{g,\tilde{h}_0}\(\exp_{\xi_0}\(\mu y\),\xi_0\)\\
&\qquad+\(1-\chi\(\mu\left|y\right|\)\)G_{g,\tilde{h}_0}\(\exp_{\xi_0}\(\mu y\),\xi_0\)\Big)\hspace{104pt}\text{in Case 1}\\&\chi\(\mu\left|y\right|\)U_0\(y\)\hspace{235pt}\text{in Case 2}\end{aligned}\right.
\end{multline*}
for all points $y\in B\(0,r_0/\mu\)$, where $\left|y\right|:=\left|y\right|_{\delta_0}$ and {$U_0$ is as in \eqref{Th5Eq3}.} Since $w_\eps\to u_0$ in $C^2\(M\)$ as $\eps\to0$, $u_0\(\xi_0\)=0$ and $w_\eps$ is symmetric with respect to $\xi_0$, a Taylor expansion gives
\begin{equation}\label{Th3Eq13}
\widehat{w}_{\eps,\mu}\(y\)=\frac{d}{d\eps}\[w_\eps\(\xi_0\)\]_{\eps=0}\eps\mu^{\frac{n-2}{2}}+\bigO\(\eps^2\mu^{\frac{n-2}{2}}+\mu^{\frac{n+2}{2}}\left|y\right|^2\)
\end{equation}
uniformly with respect to $\eps,\mu>0$ and $y\in B\(0,r_0/\mu\)$. By using \eqref{IntroEq5} and \eqref{IntroEq6} (together with \eqref{Sec3Eq1} and a result from the book of Habermann \cite{Ha}*{Lemma~5.2.10} in the case where $n=3$), we obtain
\begin{equation}\label{Th3Eq14}
\widehat{U}_{\mu}\(y\)=\left\{\begin{aligned}&U_0\(y\)+k_nm_{g,\tilde{h}_0}\(\xi_0\)\mu^{n-2}+\bigO\(\mu^{n-2}\eta\(\mu y\)\)&&\text{in Case 1}\\&U_0\(y\)+\bigO\(\mu^{n-2}\)&&\text{in Case 2}\end{aligned}\right.
\end{equation}
uniformly with respect to $\mu>0$ and $y\in B\(0,r_0/\mu\)$, for some function $\eta\in C^0\(\R^n\)$ such that $\eta\(0\)=0$. To estimate the right-hand side of \eqref{Th3Eq8}, we remark that in Case 1, for every $\mu>0$, the function $\widehat{U}_{\mu}$ can be rewritten as 
\begin{equation}\label{Th3Eq15}
\widehat{U}_{\mu}=W_\mu\(y\)+k_n\mu^{n-2}\Lambda_{\xi_0}\(\exp_{\xi_0}\(\mu y\)\)^{-1}G_{g,\tilde{h}_0}\(\exp_{\xi_0}\(\mu y\),\xi_0\)
\end{equation}
for all points $y\in B\(0,r_0/\mu\)\backslash\left\{0\right\}$, where 
$$W_\mu\(y\):=\chi\(\mu\left|y\right|\)\(U_0\(y\)-\(n\(n-2\)\)^{\frac{n-2}{4}}\left|y\right|^{2-n}\).$$
Recall moreover that $\Delta_{\delta_0}U_0=U_0^{2^*-1}$ in $\R^n$ and $\Delta_{\delta_0}\left|\cdot\right|^{2-n}=0$ in $\R^n\backslash\left\{0\right\}$. Since the functions $\chi\(\mu\,\cdot\)U_0$ and $W_\mu$ are radially symmetric about the origin, by using \eqref{Sec3Eq1} and \eqref{Th3Eq15} together with straightforward estimates, we then obtain
\begin{align}\label{Th3Eq16}
\Delta_{\hat{g}_\mu}W_{\mu}\(y\)&=\Delta_{\delta_0}W_{\mu}\(y\)+\bigO\big(\mu^{N}\left|y\right|^{N-1}\left|\nabla W_\mu\(y\)\right|\big)\nonumber\\
&=U_0\(y\)^{2^*-1}+\bigO\big(\mu^{N}\left|y\right|^{N-n-2}\big)\quad\text{in Case 1}
\end{align}
and
\begin{align}\label{Th3Eq17}
\Delta_{\hat{g}_\mu}\widehat{U}_{\mu}\(y\)&=\Delta_{\delta_0}\widehat{U}_{\mu}\(y\)+\bigO\big(\mu^{N}\left|y\right|^{N-1}\big|\nabla \widehat{U}_\mu\(y\)\big|\big)\nonumber\\
&=U_0\(y\)^{2^*-1}+\bigO\big(\mu^{N}\left|y\right|^{N-n}\big)\quad\text{in Case 2}
\end{align}
uniformly with respect to $\eps,\mu>0$ and $y\in B\(0,r_0/\mu\)\backslash\left\{0\right\}$. On the other hand, in Case 1, the conformal invariance property of the conformal Laplacian gives
\begin{equation}\label{Th3Eq18}
\big(\Delta_{\hat{g}_\mu}+\mu^2\hat{h}_{0,\mu}\big)\big(\widehat{U}_{\mu}-W_{\mu}\big)\(y\)=\mu^n\[\Lambda_{\xi_0}^{1-2^*}\big(\Delta_g+\tilde{h}_0\big)G_{g,\tilde{h}_0}\(\cdot,\xi_0\)\]\(\exp_{\xi_0}\(\mu y\)\)=0
\end{equation}
for all $y\in B\(0,r_0/\mu\)\backslash\left\{0\right\}$. Moreover, by using \eqref{Sec3Eq1} and \eqref{Sec3Eq2} together with Taylor expansions, we obtain
\begin{align}\label{Th3Eq19}
&\hat{h}_{\eps,\mu}\(y\)W_{\mu}\(y\)+\big(\hat{h}_{\eps,\mu}-\hat{h}_{0,\mu}\big)\(y\)\big(\widehat{U}_{\mu}-W_{\mu}\big)\(y\)\nonumber\\
&=\bigO\(\left\{\begin{aligned}&\left|y\right|^{-1}\big(\(1+\left|y\right|\)^{-2}+\eps\big)&&\text{if }n=3\\&\left|y\right|^{2-n}\big(\mu^2\left|y\right|^2\(1+\left|y\right|\)^{-2}+\eps\big)&&\text{if }n=4,5 \text{ and } h_0\(\xi_0\)=c_n\Scal_g\(\xi_0\)\end{aligned}\right\}\)
\end{align}
and
\begin{equation}\label{Th3Eq20}
\hat{h}_{\eps,\mu}\(y\)\widehat{U}_{\mu}\(y\)=\big(h_0\(\xi_0\)-c_n\Scal_g\(\xi_0\)+\bigO\big(\mu^2\left|y\right|^2+\eps\big)\big)U_0\(y\)+\bigO\(\mu^{n-2}\)\,\,\,\,\text{in Case 2}
\end{equation}
uniformly with respect to $\eps,\mu>0$ and $y\in B\(0,r_0/\mu\)\backslash\left\{0\right\}$. It follows from \eqref{Th3Eq16}--\eqref{Th3Eq20} that
\begin{multline}\label{Th3Eq21}
\big(\Delta_{\hat{g}_\mu}\widehat{U}_{\mu}+\mu^2\hat{h}_{\eps,\mu}\widehat{U}_{\mu}\big)\(y\)=U_0\(y\)^{2^*-1}\\
+\left\{\begin{aligned}&\bigO\big(\mu^N\left|y\right|^{N-5}+\mu^2\left|y\right|^{-1}\big(\(1+\left|y\right|\)^{-2}+\eps\big)\big)\hspace{113pt}\text{if }n=3\\&\bigO\big(\mu^N\left|y\right|^{N-n-2}+\mu^2\left|y\right|^{2-n}\big(\mu^2+\eps\big)\big)\hspace{8pt}\text{if }n=4,5 \text{ and } h_0\(\xi_0\)=c_n\Scal_g\(\xi_0\)\\&\big(h_0\(\xi_0\)-c_n\Scal_g\(\xi_0\)+\bigO\big(\mu^2\left|y\right|^2+\eps\big)\big)\mu^2U_0\(y\)+\bigO\big(\mu^N\left|y\right|^{N-n}\big)\hspace{1pt}\text{in Case 2}\end{aligned}\right.
\end{multline}
uniformly with respect to $\eps,\mu>0$ and $y\in B\(0,r_0/\mu\)$. For every $\eps,t>0$, we define
\begin{equation}\label{Th3Eq22}
\mu_\eps\(t\):=\left\{\begin{aligned}&t\eps^{\frac{2}{n-2}}&&\text{in Case 1}\\&t\nu_\eps&&\text{in Case 2.1}\\&t\eps^2&&\text{in Case 2.2},\end{aligned}\right.
\end{equation}
where Case 2.1 and Case 2.2 are subcases of Case 2 defined as 
$$\left\{\begin{aligned}&\text{Case 2.1: }n=4\text{ and }h\(\xi_0\)\ne c_n\Scal_g\(\xi_0\)\\&\text{Case 2.2: }n=5\text{ and }h\(\xi_0\)\ne c_n\Scal_g\(\xi_0\)\end{aligned}\right.$$
and $\nu_\eps>0$ is such that $\nu_\eps\ln\(1/\nu_\eps\)=\eps$ and $\nu_\eps\to0$ as $\eps\to0$. It follows from \eqref{Sec3Eq1} and \eqref{Th3Eq8}--\eqref{Th3Eq22} together with standard integral estimates (see for instance the article by Esposito, Pistoia and V\'etois~\cite{EPV}*{Section 4}) that
\begin{align}
&\int_M\(\left|\nabla U_{\mu_\eps\(t\)}\right|_g^2+\tilde{h}_\eps U_{\mu_\eps\(t\)}^2\)dv_g=\int_{\R^n} U_0^{2^*}dv_{\delta_0}\nonumber\\
&\qquad+\left\{\begin{aligned}&k_nm_{g,\tilde{h}_0}\(\xi_0\)t^{n-2}\eps^2\int_{\R^n}U_0^{2^*-1}dv_{\delta_0}+\smallo\(\eps^2\)&&\text{in Case 1}\\&8\,\omega_3\(h_0\(\xi_0\)-c_n\Scal_g\(\xi_0\)+\smallo\(1\)\)t^2\eps\nu_\eps&&\text{in Case 2.1}\\&\(h_0\(\xi_0\)-c_n\Scal_g\(\xi_0\)+\smallo\(1\)\)t^2\eps^4\int_{\R^5}U_0^2\,dv_{\delta_0}&&\text{in Case 2.2},\end{aligned}\right.\label{Th3Eq23}\allowdisplaybreaks\\
&\int_M\left|U_{\mu_\eps\(t\)}\right|^{2^*}dv_g=\int_{\R^n}U_0^{2^*}dv_{\delta_0} \nonumber\\
&\qquad+\left\{\begin{aligned}&2^*k_nm_{g,\tilde{h}_0}\(\xi\)t^{n-2}\eps^2\int_{\R^n}U_0^{2^*-1}dv_{\delta_0}+\smallo\(\eps^2\)&&\text{in Case 1}\\&\smallo\(\eps\nu_\eps\)&&\text{in Case 2.1}\\&\smallo\(\eps^4\)&&\text{in Case 2.2}\end{aligned}\right.\label{Th3Eq24}
\end{align}
and
\begin{multline}\label{Th3Eq25}
\int_M\left|U_{\mu_\eps\(t\)}\right|^{2^*-2}U_{\mu_\eps\(t\)}w_\eps\,dv_g\\
=\left\{\begin{aligned}&\frac{d}{d\eps}\[w_\eps\(\xi_0\)\]_{\eps=0}t^{\frac{n-2}{2}}\eps^2\int_{\R^n}U_0^{2^*-1}dv_{\delta_0}+\smallo\(\eps^2\)&&\text{in Case 1}\\&\frac{d}{d\eps}\[w_\eps\(\xi_0\)\]_{\eps=0}t\eps\nu_\eps\int_{\R^n}U_0^3\,dv_{\delta_0}+\smallo\(\eps\nu_\eps\)&&\text{in Case 2.1}\\&\frac{d}{d\eps}\[w_\eps\(\xi_0\)\]_{\eps=0}t^{3/2}\eps^4\int_{\R^n}U_0^{7/3}\,dv_{\delta_0}+\smallo\(\eps^4\)&&\text{in Case 2.2}\end{aligned}\right.
\end{multline}
as $\eps\to0$, uniformly with respect to $t$ in compact subsets of $\(0,\infty\)$. Moreover, by using \eqref{Th3Eq13} and \eqref{Th3Eq22}, we obtain that for every compact set $A\subset\(0,\infty\)$, there exists a constant $r_A>0$ such that for every $\eps>0$ and $t\in A$,
\begin{equation}\label{Th3Eq26}
B\(0,r_A\varrho_\eps\)\subset\Omega_{\eps,\mu_\eps\(t\)},\quad\text{where }\varrho_\eps:=\left\{\begin{aligned}&\eps^{-\frac{n+2}{n\(n-2\)}}&&\text{in Case 1}\\&\nu_\eps^{-3/4}&&\text{in Case 2.1}\\&\eps^{-4/3}&&\text{in Case 2.2.}\end{aligned}\right.
\end{equation}
It follows from \eqref{Th3Eq13}, \eqref{Th3Eq14} and \eqref{Th3Eq26} that
\begin{multline}\label{Th3Eq27}
\int_{\Omega_{\eps,\mu_\eps\(t\)}}\(\big|\widehat{U}_{\mu_\eps\(t\)}\big|^{2^*-2}\widehat{w}_{\eps,\mu_\eps\(t\)}^2+\left|\widehat{w}_{\eps,\mu_\eps\(t\)}\right|^{2^*-2}\widehat{U}_{\mu_\eps\(t\)}^2\)dv_{\hat{g}_\mu}=\bigO\Bigg(\int_{B\(0,r_A\varrho_\eps\)}\\
\bigg(\(\eps\mu_\eps\(t\)^{\frac{n-2}{2}}+\mu_\eps\(t\)^{\frac{n+2}{2}}\left|\cdot\right|^2\)^2U_0^{2^*-2}+\(\eps\mu_\eps\(t\)^{\frac{n-2}{2}}+\mu_\eps\(t\)^{\frac{n+2}{2}}\left|\cdot\right|^2\)^{2^*-2}U_0^2\bigg)dv_{\delta_0}\\
+\int_{B\(0,r_0/\mu_\eps\(t\)\)\backslash B\(0,r_A\varrho_\eps\)}U_0^{2^*}dv_{\delta_0}\Bigg)
\end{multline}
and
\begin{multline}\label{Th3Eq28}
\int_{B\(0,r_0/\mu_\eps\(t\)\)\backslash \Omega_{\eps,\mu_\eps\(t\)}}\Big(\big|\widehat{U}_{\mu_\eps\(t\)}\big|^{2^*-1}\left|\widehat{w}_{\eps,\mu_\eps\(t\)}\right|+\left|\widehat{w}_{\eps,\mu_\eps\(t\)}\right|^{2^*-3}\big|\widehat{U}_{\mu_\eps\(t\)}\big|^3\Big)dv_{\hat{g}_\mu}\\
=\bigO\bigg(\int_{B\(0,r_0/\mu_\eps\(t\)\)\backslash B\(0,r_A\varrho_\eps\)}\bigg(\(\eps\mu_\eps\(t\)^{\frac{n-2}{2}}+\mu_\eps\(t\)^{\frac{n+2}{2}}\left|\cdot\right|^2\)U_0^{2^*-1}\\
+\(\eps\mu_\eps\(t\)^{\frac{n-2}{2}}+\mu_\eps\(t\)^{\frac{n+2}{2}}\left|\cdot\right|^2\)^{2^*-3}U_0^3\bigg)dv_{\delta_0}\bigg)
\end{multline}
uniformly with respect to $\eps>0$ and $t\in A$. By using \eqref{Th3Eq11}, \eqref{Th3Eq12}, \eqref{Th3Eq22}, \eqref{Th3Eq27} and \eqref{Th3Eq28} together with standard integral estimates, we then obtain
\begin{multline}\label{Th3Eq29}
\int_{M_{\eps,\mu_\eps\(t\)}}\(\left|U_{\mu_\eps\(t\)}\right|^{2^*-2}w_\eps^2+\left|w_\eps\right|^{2^*-2}U_{\mu_\eps\(t\)}^2\)dv_g+\int_{M\backslash M_{\eps,\mu_\eps\(t\)}}\Big(\left|U_{\mu_\eps\(t\)}\right|^{2^*-1}\left|w_\eps\right|\\
+\left|w_\eps\right|^{2^*-3}\left|U_{\mu_\eps\(t\)}\right|^3\Big)dv_g=\smallo\(\left\{\begin{aligned}&\eps^2&&\text{in Case 1}\\&\eps\nu_\eps&&\text{in Case 2.1}\\&\eps^4&&\text{in Case 2.2}\end{aligned}\right\}\)
\end{multline}
as $\eps\to0$, uniformly with respect to $t$ in compact subsets of $\(0,\infty\)$. Note that in Case 2.1, we used that $\nu_\eps = \smallo \(\eps\)$ by \eqref{Th3Eq22}. It follows from \eqref{Th3Eq6}, \eqref{Th3Eq23}--\eqref{Th3Eq25} and \eqref{Th3Eq29} that 
\begin{multline}\label{Th3Eq30}
J_\eps\(w_\eps-U_{\mu_\eps\(t\)}\)=\frac{1}{n}\int_{M}\left|w_\eps\right|^{2^*}dv_g+\frac{1}{n}\int_{\R^n}U_0^{2^*}dv_{\delta_0}\\
+\left\{\begin{aligned}&F\(t\)\eps^2+\smallo\(\eps^2\)&&\text{in Case 1}\\&F\(t\)\eps\nu_\eps+\smallo\(\eps\nu_\eps\)&&\text{in Case 2.1}\\&F\(t\)\eps^4+\smallo\(\eps^4\)&&\text{in Case 2.2}\end{aligned}\right.
\end{multline}
as $\eps\to0$, uniformly with respect to $t$ in compact subsets of $\(0,\infty\)$,where
\begin{multline*}
F\(t\):=\frac{d}{d\eps}\[w_\eps\(\xi_0\)\]_{\eps=0}t^{\frac{n-2}{2}}\int_{\R^n}U_0^{2^*-1}dv_{\delta_0}\\
+\left\{\begin{aligned}&-\frac{k_n}{2}\,m_{g,\tilde{h}_0}\(\xi_0\)t^{n-2}\int_{\R^n}U_0^{2^*-1}dv_{\delta_0}&&\text{in Case 1}\\&4\,\omega_3\(h_0\(\xi_0\)-c_n\Scal_g\(\xi_0\)\)t^2&&\text{in Case 2.1}\\&\frac{1}{2}\(h_0\(\xi_0\)-c_n\Scal_g\(\xi_0\)\)t^2\int_{\R^5}U_0^2\,dv_{\delta_0}&&\text{in Case 2.2.}\end{aligned}\right.
\end{multline*}
We now consider the error term $R_{\eps,\mu}$ defined in \eqref{Th3Eq5}. By continuity of the embedding $H^1\(M\)\hookrightarrow  L^{2^*}\(M\)$, we obtain
\begin{equation}\label{Th3Eq31}
R_{\eps,\mu}=\bigO\(\left\|\(\Delta_g+h_\eps-\left|w_\eps-U_{\mu}\right|^{2^*-2}\)\(w_\eps-U_{\mu}\)\right\|_{L^{\frac{2n}{n+2}}\(M\)}\)
\end{equation}
uniformly with respect to $\eps,\mu>0$. By using \eqref{Th3Eq2} and \eqref{Th3Eq31} together with similar estimate as in \eqref{Th3Eq6}--\eqref{Th3Eq12}, we obtain
\begin{align}\label{Th3Eq32}
R_{\eps,\mu}&=\bigO\bigg(\left\|\Delta_gU_{\mu}+\(h_\eps-\(2^*-1\)\left|w_\eps\right|^{2^*-2}\)U_{\mu}-\left|U_{\mu}\right|^{2^*-2}U_{\mu}\right\|_{L^{\frac{2n}{n+2}}\(M\)}\allowdisplaybreaks\nonumber\\
&\qquad+\left\|\left|U_{\mu}\right|^{2^*-2}\left|w_\eps\right|+\left|w_\eps\right|^{2^*-2}\left|U_{\mu}\right|\right\|_{L^{\frac{2n}{n+2}}\(M_{\eps,\mu}\)}\allowdisplaybreaks\nonumber\\
&\qquad+\left\|\left|U_{\mu}\right|^{2^*-1}+\left|w_\eps\right|^{2^*-3}U_{\mu}^2\right\|_{L^{\frac{2n}{n+2}}\(M\backslash M_{\eps,\mu}\)}\bigg)\allowdisplaybreaks\nonumber\\
&=\bigO\bigg(\Big\|\Delta_{\hat{g}_\mu}\widehat{U}_{\mu}+\mu^2\hat{h}_{\eps,\mu}\widehat{U}_{\mu}-\big|\widehat{U}_{\mu}\big|^{2^*-2}U_{\mu}\Big\|_{L^{\frac{2n}{n+2}}\(B\(0,r_0/\mu\)\)}\allowdisplaybreaks\nonumber\\
&\qquad+\Big\|\big|\widehat{U}_{\mu}\big|^{2^*-2}\left|\widehat{w}_{\eps,\mu}\right|+\left|\widehat{w}_{\eps,\mu}\right|^{2^*-2}\big|\widehat{U}_{\mu}\big|\Big\|_{L^{\frac{2n}{n+2}}\(\Omega_{\eps,\mu}\)}\nonumber\\
&\qquad+\Big\|\big|\widehat{U}_{\mu}\big|^{2^*-1}+\left|\widehat{w}_{\eps,\mu}\right|^{2^*-3}\widehat{U}_{\mu}^2\Big\|_{L^{\frac{2n}{n+2}}\(B\(0,r_0/\mu\)\backslash \Omega_{\eps,\mu}\)}\bigg)+\smallo\(\mu^{\frac{n-2}{2}}\)
\end{align}
as $\eps,\mu\to0$. By using \eqref{Th3Eq14} and \eqref{Th3Eq21}, we obtain
\begin{equation}\label{Th3Eq33}
\Big\|\Delta_{\hat{g}_\mu}\widehat{U}_{\mu}+\mu^2\hat{h}_{\eps,\mu}\widehat{U}_{\mu}-\big|\widehat{U}_{\mu}\big|^{2^*-2}U_{\mu}\Big\|_{L^{\frac{2n}{n+2}}\(B\(0,r_0/\mu\)\)}=\left\{\begin{aligned}&\smallo\(\mu^{\frac{n-2}{2}}\)&&\text{in Case 1}\\&\bigO\(\mu^{\frac{n-2}{2}}\)&&\text{in Case 2}\end{aligned}\right.
\end{equation}
as $\eps,\mu\to0$. By proceeding in a similar way as in \eqref{Th3Eq27}--\eqref{Th3Eq29}, by using \eqref{Th3Eq13}, \eqref{Th3Eq14}, \eqref{Th3Eq22} and \eqref{Th3Eq26}, we obtain
 \begin{multline}\label{Th3Eq34}
\Big\|\big|\widehat{U}_{\mu_\eps\(t\)}\big|^{2^*-2}\left|\widehat{w}_{\eps,\mu_\eps\(t\)}\right|+\left|\widehat{w}_{\eps,\mu_\eps\(t\)}\right|^{2^*-2}\big|\widehat{U}_{\mu_\eps\(t\)}\big|\Big\|_{L^{\frac{2n}{n+2}}\(\Omega_{\eps,\mu_\eps\(t\)}\)}+\Big\|\big|\widehat{U}_{\mu_\eps\(t\)}\big|^{2^*-1}\\
+\left|\widehat{w}_{\eps,\mu_\eps\(t\)}\right|^{2^*-3}\widehat{U}_{\mu_\eps\(t\)}^2\Big\|_{L^{\frac{2n}{n+2}}\(B\(0,r_0/\mu_\eps\(t\)\)\backslash \Omega_{\eps,\mu_\eps\(t\)}\)}=\smallo\(\mu_\eps\(t\)^{\frac{n-2}{2}}\) 
\end{multline}
as $\eps\to0$, uniformly with respect to $t$ in compact subsets of $\(0,\infty\)$. It follows from \eqref{Th3Eq22} and \eqref{Th3Eq32}--\eqref{Th3Eq34} that
\begin{equation}\label{Th3Eq35}
R_{\eps,\mu_\eps\(t\)}=\left\{\begin{aligned}&\smallo\(\eps\)&&\text{in Case 1}\\&\bigO\(\nu_\eps\)&&\text{in Case 2.1}\\&\bigO\(\eps^3\)&&\text{in Case 2.2}\end{aligned}\right.
\end{equation}
as $\eps\to0$, uniformly with respect to $t$ in compact subsets of $\(0,\infty\)$. By using \eqref{Th3Eq30} and \eqref{Th3Eq35} and performing a Lyapunov--Schmidt reduction in the subspace of $H^1\(M\)$ which consists of symmetric functions with respect to $\xi_0$ (see the article by Morabito, Pistoia and Vaira \cite{MPV}; see also the article by Robert and V\'etois \cite{RV2}*{Theorem~1.1}, where the proof of the main result can easily be adapted to this case), we obtain that for small $\eps,\mu>0$, there exists a function $\psi_{\eps,\mu}\in H^1\(M\)$ such that $\psi_{\eps,\mu}\to0$ in $H^1\(M\)$ as $\eps,\mu\to0$,
\begin{multline}\label{Th3Eq36}
J_\eps\(w_\eps-U_{\mu_\eps\(t\)}+\psi_{\eps,\mu_\eps\(t\)}\)=\frac{1}{n}\int_{M}\left|w_\eps\right|^{2^*}dv_g+\frac{1}{n}\int_{\R^n}U_0^{2^*}dv_{\delta_0}\\
+\left\{\begin{aligned}&F\(t\)\eps^2+\smallo\(\eps^2\)&&\text{in Case 1}\\&F\(t\)\eps\nu_\eps+\smallo\(\eps\nu_\eps\)&&\text{in Case 2.1}\\&F\(t\)\eps^4+\smallo\(\eps^4\)&&\text{in Case 2.2}\end{aligned}\right.
\end{multline}
as $\eps\to0$, uniformly with respect to $t$ in compact subsets of $\(0,\infty\)$, and
\begin{equation}\label{Th3Eq37}
D_uJ_\eps\(w_\eps-U_{\mu}+\psi_{\eps,\mu}\)=0\quad\Longleftrightarrow\quad \frac{d}{d\mu}\[J_\eps\(u_0-U_{\mu}+\psi_{\eps,\mu}\)\]=0.
\end{equation}
By using \eqref{Th3Eq1} and \eqref{Th3Eq3}, we obtain that the function $F$ has a strict maximum or minimum point in $\(0,\infty\)$, that we call $t_0$. By using \eqref{Th3Eq36}, we then obtain that for large $C>0$, there exists a family of maximum or minimum points $\(t_\eps\)_{0<\eps\ll1}$ of the function $J_\eps\(w_\eps-U_{\mu_\eps\(\cdot\)}+\psi_{\eps,\mu_\eps\(\cdot\)}\)$ in $\(1/C,C\)$ that converges to $t_0$ as $\eps \to 0$. It then follows from \eqref{Th3Eq37} that for small $\eps>0$, the function $u_\eps:=w_\eps-U_{\mu_\eps\(t_\eps\)}+\psi_{\eps,\mu_\eps\(t_\eps\)}$ is a solution to \eqref{IntroEq2}. It is easy to verify that $u_\eps$ satisfies \eqref{IntroEq3}. This ends the proof of Theorem~\ref{Th3}.
\endproof

To prove Theorem~\ref{Th4}, we also need the following:

\begin{proposition}\label{Pr1}
Let $\(M,g\)$ be a smooth, closed Riemannian manifold of dimension $n\ge3$. Assume that $\(M,g\)$ is $2$-symmetric with respect to a point $\xi_0\in M$. Let $\sigma_0$, $\sigma_1$, $\sigma_2$, $\Gamma_{\sigma_0}$, $\Gamma_{\sigma_1}$ and $\Gamma_{\sigma_2}$ be as in the definition of $2$-symmetric manifolds. Then the following assertions hold true:
\begin{enumerate}
\item[(i)]For every point $\xi\in\Gamma_{\sigma_1}$ or $\xi\in\Gamma_{\sigma_2}$, the connected component of $\Gamma_{\sigma_1}$ or $\Gamma_{\sigma_2}$ which contains $\xi$ is a totally geodesic submanifold of dimension $n-1$ of $M$.
\item[(ii)] For every point $\xi\in\Gamma_{\sigma_1}\cap\Gamma_{\sigma_2}$, in geodesic normal coordinates, the maps $d\sigma_1\(\xi\)$ and $d\sigma_2\(\xi\)$ are two reflections about orthogonal hyperplanes $H_1$ and $H_2$ of $T_\xi M$. If moreover $\xi=\xi_0$, then in geodesic normal coordinates, the map $d\sigma_0\(\xi_0\)$ is the reflection about the plane which passes through the origin and is spanned by the normal vectors to the hyperplanes $H_1$ and $H_2$.
\end{enumerate}
In particular, it follows from the point (ii) that $\Omega\cap\Gamma_{\sigma_0}\ne\emptyset$.
\end{proposition}

\proof[Proof of Proposition~\ref{Pr1}]
We begin with proving the point (i). We let $\xi$ be a point in $\Gamma_{\sigma_1}\backslash \Gamma_{\sigma_2}$ (the case where $\xi\in\Gamma_{\sigma_1}\cap\Gamma_{\sigma_2}$ is treated in the point (ii) and the case where $\xi\in\Gamma_{\sigma_2}\backslash \Gamma_{\sigma_1}$ is similar to the case where $\xi\in\Gamma_{\sigma_1}\backslash \Gamma_{\sigma_2}$). We let $C_\xi$ be the connected component of $\Gamma_{\sigma_1}$ which contains the point $\xi$. Since $\sigma_1$ is an isometry of $M$, we obtain that $C_\xi$ is a totally geodesic submanifold of $M$. Let $n_1$ be the dimension of $C_\xi$. Since $M$ is connected and $M\backslash\Gamma_{\sigma_1}\ne\emptyset$, we obtain that $n_1<n$. We now assume by contradiction that $n_1<n-1$. We then obtain that the set $B_g\(\xi,r\)\backslash\Gamma_{\sigma_1}$ is connected for small $r>0$, where $B_g\(\xi,r\)$ is the ball of center $\xi$ and radius $r$ with respect to the metric $g$. Since $\xi\in\Gamma_{\sigma_1}\backslash\Gamma_{\sigma_2}$, $M=\Omega\sqcup\sigma_1\(\Omega\)\sqcup\sigma_2\(\Omega\)\sqcup\sigma_1\circ\sigma_2\(\Omega\)\sqcup\(\Gamma_{\sigma_1}\cup\Gamma_{\sigma_2}\)$, the sets $B_g\(\xi,r\)\backslash\Gamma_{\sigma_1}$, $\Omega$, $\sigma_1\(\Omega\)$, $\sigma_2\(\Omega\)$ and $\sigma_1\circ\sigma_2\(\Omega\)$ are open and the set $\Gamma_{\sigma_1}\cup\Gamma_{\sigma_2}$ is closed, it follows that for small $r>0$, the set $B_g\(\xi,r\)\backslash\Gamma_{\sigma_1}$ is included in only one set $S\in\left\{\Omega,\sigma_1\(\Omega\),\sigma_2\(\Omega\),\sigma_1\circ\sigma_2\(\Omega\)\right\}$. On the other hand , for every sequence of points $\(x_k\)_{k\in\N}$ in $S$ such that $x_k\to\xi$ as $k\to\infty$, by continuity of $\sigma_1$, we obtain that $\sigma_1\(x_k\)\to\sigma_1\(\xi\)=\xi$ as $k\to\infty$, which implies that $\sigma_1\(S\)\cap\(S\cup\Gamma_{\sigma_1}\)\ne\emptyset$. This clearly contradicts the definitions of $\sigma_1$, $\sigma_2$ and $\Omega$, thus proving that $n_1=n-1$, which ends the proof of the point (i).

\smallskip
We now prove the point (ii). We let $\xi$ be a point in $\Gamma_{\sigma_1}\cap\Gamma_{\sigma_2}$. By similar considerations as in the proof of the point (i), we obtain that the connected components of $\Gamma_{\sigma_1}$ and $\Gamma_{\sigma_2}$ which contain the point $\xi$ are totally geodesic submanifold of $M$ of dimensions $n_1$ and $n_2$, respectively, such that $n_1<n$ and $n_2<n$. We assume by contradiction that $n_1<n-1$ or $n_2<n-1$. We then obtain that for small $r>0$, the set $B_g\(\xi,r\)\backslash\(\Gamma_{\sigma_1}\cup\Gamma_{\sigma_2}\)$ has at most two distinct connected components. Similarly as in the proof of the point (i), it follows that for small $r>0$, $B_g\(\xi,r\)\backslash\(\Gamma_{\sigma_1}\cup\Gamma_{\sigma_2}\)\subset S_1\cup S_2$ for some sets $S_1,S_2\in\left\{\Omega,\sigma_1\(\Omega\),\sigma_2\(\Omega\),\sigma_1\circ\sigma_2\(\Omega\)\right\}$. On the other hand , for every sequence of points $\(x_k\)_{k\in\N}$ in $S_1\cup S_2$ such that $x_k\to\xi$ as $k\to\infty$, by continuity of $\sigma_1$ and $\sigma_2$, we obtain that $\sigma_1\(x_k\)\to\sigma_1\(\xi\)=\xi$, $\sigma_2\(x_k\)\to\sigma_2\(\xi\)=\xi$ and $\sigma_1\circ\sigma_2\(x_k\)\to\sigma_1\circ\sigma_2\(\xi\)=\xi$ as $k\to\infty$. This together with the properties of $\sigma_1$ and $\sigma_2$ implies that $\xi\in\partial\Omega\cap\partial\(\sigma_1\(\Omega\)\)\cap\partial\(\sigma_2\(\Omega\)\)\cap\partial\(\sigma_1\circ\sigma_2\(\Omega\)\)$, which is in contradiction with $B_g\(\xi,r\)\backslash\(\Gamma_{\sigma_1}\cup\Gamma_{\sigma_2}\)\subset S_1\cup S_2$. This proves that $n_1=n_2=n-1$. Since the reflections about hyperplanes are the only Euclidean isometries which fixed point sets are hyperplane, it follows that in geodesic normal coordinates, the maps $d\sigma_1\(\xi\)$ and $d\sigma_2\(\xi\)$ are two reflections about hyperplanes $H_1$ and $H_2$ of $T_\xi M$. Moreover, since $\sigma_1\circ\sigma_2=\sigma_2\circ\sigma_1$ in $M$, we obtain that $d\sigma_1\(\xi\)\circ d\sigma_2\(\xi\)=d\sigma_2\(\xi\)\circ d\sigma_1\(\xi\)$, which implies that either $H_1=H_2$ or $H_1$ is orthogonal to $H_2$. If $H_1=H_2$, then we obtain that $\Gamma_{\sigma_1\circ\sigma_2}$ contains a neighborhood of the point $\xi$, which is in contradiction with $M=\Omega\sqcup\sigma_1\(\Omega\)\sqcup\sigma_2\(\Omega\)\sqcup\sigma_1\circ\sigma_2\(\Omega\)\sqcup\(\Gamma_{\sigma_1}\cup\Gamma_{\sigma_2}\)$. Therefore, we obtain that $H_1$ is orthogonal to $H_2$. Finally, in the case where $\xi=\xi_0$, for every vector $v\in T_{\xi_0} M$, since $d\(\sigma_0\circ\sigma_1\circ\sigma_2\)_{\xi_0}\(v\)=-v$, $\sigma_0=\(\sigma_0\circ\sigma_1\circ\sigma_2\)\circ\sigma_2\circ\sigma_1$ and $\xi_0\in\Gamma_{\sigma_0}\cap\Gamma_{\sigma_1}\cap\Gamma_{\sigma_2}$, we obtain that $d\(\sigma_0\)_{\xi_0}\(v\)=-d\(\sigma_2\)_{\xi_0}\circ d\(\sigma_1\)_{\xi_0}\(v\)$. It follows that in geodesic normal coordinates, the map $d\sigma_0\(\xi_0\)$ is the reflection about the plane which passes through the origin and is spanned by the normal vectors to the hyperplanes $H_1$ and $H_2$. This ends the proof of the point (ii).
\endproof

We can now use Theorem~\ref{Th3} and Proposition~\ref{Pr1} to prove Theorem~\ref{Th4}.

\proof[Proof of Theorem~\ref{Th4}]

\begin{step}\label{Th4Step1}
The function $u_0$ can be extended as a solution to \eqref{IntroEq1} such that $u_0$ is symmetric with respect to $\xi_0$, $u_0>0$ in $\Omega$ and $u_0\circ\sigma_1=u_0\circ\sigma_2=-u_0$ in $M$.
\end{step}

\proof[Proof of Step~\ref{Th4Step1}]
Since $M=\Omega\sqcup\sigma_1\(\Omega\)\sqcup\sigma_2\(\Omega\)\sqcup\sigma_1\circ\sigma_2\(\Omega\)\sqcup\(\Gamma_{\sigma_1}\cup\Gamma_{\sigma_2}\)$, we can extend $u_0$ to $M$ by letting 
$$\left\{\begin{aligned}&u_0:=-u_0\circ \sigma_1\text{ in }\sigma_1\(\Omega\),\quad u_0:=-u_0\circ \sigma_2\text{ in }\sigma_2\(\Omega\),\\
&u_0:=u_0\circ \sigma_1\circ \sigma_2\text{ in }\sigma_1\circ \sigma_2\(\Omega\)\quad\text{and}\quad u_0:=0\text{ in }\Gamma_{\sigma_1}\cup\Gamma_{\sigma_2}.\end{aligned}\right.$$
By using the composition properties of $\sigma_0$, $\sigma_1$ and $\sigma_2$ and since $u_0\circ\sigma_0=u_0$ in $\Omega$, we obtain that $u_0$ is symmetric with respect to $\xi_0$ for $\sigma:= \sigma_0 \circ \sigma_1 \circ \sigma_2$, $u_0\in H^1\(M\)$ and $u_0\circ\sigma_1=u_0\circ\sigma_2=-u_0$ in $M$. To prove that $u_0$ is a solution to \eqref{IntroEq1}, for every function $v\in C^\infty\(M\)$, we define 
$$v_{\sym}:=v-v\circ \sigma_1-v\circ \sigma_2+v\circ \sigma_1\circ \sigma_2.$$
Since $\partial\Omega\subset\Gamma_{\sigma_1}\cup\Gamma_{\sigma_2}$, we obtain that $v_{\sym}=0$ on $\partial\Omega$. By using \eqref{IntroEq7} together with a straightforward change of variable and the composition properties of $\sigma_0$, $\sigma_1$ and $\sigma_2$, we then obtain
\begin{align}\label{Th4Step1Eq}
\int_{M}\(\<\nabla u_0,\nabla v\>_{g}+h_0 u_0 v\)dv_g&=\int_{\Omega}\(\<\nabla u_0,\nabla v_{\sym}\>_{g}+h_0 u_0 v_{\sym}\)dv_{g}\nonumber\\
&=\int_{\Omega}u_0^{2^*-1}v_{\sym}\,dv_{g}=\int_M\left|u_0\right|^{2^*-2}u_0 v\,dv_{g}.
\end{align}
Since \eqref{Th4Step1Eq} holds true for all functions $v\in C^\infty\(M\)$, we obtain that $u_0$ is a solution to \eqref{IntroEq1}. Since $h_0\in C^1\(M\)$, by standard elliptic regularity, we then obtain that $u_0\in C^2\(M\)$. This ends the proof of Step~\ref{Th4Step1}.
\endproof

\begin{step}\label{Th4Step2}
There exist families of functions $\(h_{\gamma}\)_{0<\gamma\ll1}$ in $C^p\(M\)$ and $\(u_{\gamma}\)_{0<\gamma\ll1}$ in $C^2\(M\)$ such that $h_{\gamma}\to h_0$ in $C^p\(M\)$ and $u_{\gamma}\to u_0$ in $C^2\(M\)$ as $\gamma\to0$ and for each $\gamma$, $h_{\gamma}$ is $2$-symmetric with respect to $\xi_0$, $u_{\gamma}$ is symmetric with respect to $\xi_0$, $u_{\gamma}\circ\sigma_1=u_{\gamma}\circ\sigma_2=-u_{\gamma}$ in $M$ and $u_{\gamma}$ is a nondegenerate solution to the equation
\begin{equation}\label{Th4Step2Eq0}
\Delta_gu_{\gamma}+h_{\gamma}u_{\gamma}=\left|u_{\gamma}\right|^{2^*-2}u_{\gamma}\quad\text{in }M.
\end{equation}
\end{step}

\proof[Proof of Step~\ref{Th4Step2}]
We use a similar idea as in the article by Robert and V\'etois \cite{RV4}*{Lemma~13.1}, but the fact that $u_0$ vanishes requires to make some nontrivial modifications. Since $u_0\in C^2\(M\)$, a standard density argument gives that there exists a sequence $\(v_\beta\)_{\beta\in\(0,1\)}$ in $C^\infty\(M\)$ such that $v_\beta\to u_0$ in $C^2\(M\)$ as $\beta\to0$. Moreover, since $u_0$ is symmetric with respect to $\xi_0$ and $u_0\circ\sigma_1=u_0\circ\sigma_2=-u_0$ in $M$, we can assume that for every $\beta\in\(0,1\)$, $v_\beta$ is also symmetric with respect to $\xi_0$ and $v_\beta\circ\sigma_1=v_\beta\circ\sigma_2=-v_\beta$ in $M$. For every $\alpha\in\(0,1\)$, we then let $\beta_\alpha>0$ be small enough so that
\begin{equation}\label{Th4Step2Eq1}
\left\{x\in M:\,v_\beta\(x\)>\alpha\right\}\subset\left\{x\in M:\,u_0\(x\)>\alpha/2\right\}
\end{equation}
and
\begin{equation}\label{Th4Step2Eq2}
\left\{x\in M:\,v_\beta\(x\)<-\alpha\right\}\subset\left\{x\in M:\,u_0\(x\)<-\alpha/2\right\}
\end{equation}
for all $\beta\in\(0,\beta_\alpha\)$. We let $\eta$ be a smooth, increasing cutoff function in $\R$ such that 
\begin{equation}\label{Th4Step2Eq3}
\eta\(t\)=0\quad\forall t\in\[-1,1\]\quad\text{and}\quad\eta\(t\)=t\quad\forall t\in\(-\infty,-2\]\cup\[2,\infty\).
\end{equation}
For every $\alpha,\beta>0$ and $t\in\R$, we define  
\begin{equation}\label{Th4Step2Eq4}
\eta_{\alpha,\beta}\(t\):=\left\{\begin{aligned}&\beta\,\eta\(\(t+\alpha\)/\beta\)&&\text{if }t<-\alpha\\&0&&\text{if }-\alpha\le t\le\alpha\\&\beta\,\eta\(\(t-\alpha\)/\beta\)&&\text{if }t>\alpha.\end{aligned}\right.
\end{equation}
For every $\alpha,\gamma\in\(0,1\)$ and $\beta\in\(0,\beta_\alpha\)$, we then define 
\begin{equation}\label{Th4Step2Eq5}
u_{\alpha,\beta,\gamma}:=u_0+\gamma w_{\alpha,\beta},\quad w_{\alpha,\beta}:=\eta_{\alpha,\beta}\(v_\beta\),
\end{equation}
and 
\begin{equation}\label{Th4Step2Eq6}
h_{\alpha,\beta,\gamma}:=\left|u_{\alpha,\beta,\gamma}\right|^{2^*-2}-\frac{\Delta_gu_{\alpha,\beta,\gamma}}{u_{\alpha,\beta,\gamma}}=\left|u_{\alpha,\beta,\gamma}\right|^{2^*-2}-\frac{\left|u_0\right|^{2^*-2}u_0-h_0u_0+\gamma\Delta_gw_{\alpha,\beta}}{u_{\alpha,\beta,\gamma}}\,,
\end{equation}
so that 
\begin{equation}\label{Th4Step2Eq7}
\Delta_gu_{\alpha,\beta,\gamma}+h_{\alpha,\beta,\gamma}u_{\alpha,\beta,\gamma}=\left|u_{\alpha,\beta,\gamma}\right|^{2^*-2}u_{\alpha,\beta,\gamma}.
\end{equation}
Remark that \eqref{Th4Step2Eq1} and \eqref{Th4Step2Eq2} give 
\begin{equation}\label{Th4Step2Eq8}
\left\{\begin{aligned}&u_{\alpha,\beta,\gamma}\le u_0<-\alpha/2&&\text{if }v_\beta<-\alpha\\&u_{\alpha,\beta,\gamma}=u_0\text{ and }h_{\alpha,\beta,\gamma}=h_0&&\text{if }-\alpha\le v_\beta\le\alpha\\&u_{\alpha,\beta,\gamma}\ge u_0>\alpha/2&&\text{if }v_\beta>\alpha,\end{aligned}\right.
\end{equation}
so that in particular $h_{\alpha,\beta,\gamma}$ is well defined. Since $h_0\in C^p\(M\)$ and $p\ge1$, standard elliptic regularity gives that $u_0\in C^2\(M\)\cap C^{p+1}\(M_0\)$, where $M_0:=\left\{x\in M:\,u_0\(x\)\ne0\right\}$. Since $\eta\in C^\infty\(\R\)$ and $\eta=0$ in $\[-1,1\]$, it then follows from \eqref{Th4Step2Eq6} and \eqref{Th4Step2Eq8} that $h_{\alpha,\beta,\gamma}\in C^p\(M\)$, $u_{\alpha,\beta,\gamma}\in C^2\(M\)$, $h_{\alpha,\beta,\gamma}\to h_0$ in $C^p\(M\)$  and $u_{\alpha,\beta,\gamma}\to u_0$ in $C^2\(M\)$ as $\gamma\to0$. Moreover, it follows from Step~\ref{Th4Step1}, \eqref{Th4Step2Eq5} and \eqref{Th4Step2Eq6} together with the symmetry properties of $h_0$, $u_0$ and $v_\beta$ that $h_{\alpha,\beta,\gamma}$ is 2-symmetric with respect to $\xi_0$, $u_{\alpha,\beta,\gamma}$ is symmetric with respect to $\xi_0$ and $u_{\alpha,\beta,\gamma}\circ\sigma_1=u_{\alpha,\beta,\gamma}\circ\sigma_2=-u_{\alpha,\beta,\gamma}$ in $M$. We claim that there exist $\alpha\in\(0,1\)$ and $\beta\in\(0,\beta_\alpha\)$ such that for small $\gamma>0$, $u_{\alpha,\beta,\gamma}$ is a non-degenerate solution to \eqref{Th4Step2Eq7}. Assume by contradiction that this is not the case. Then for every $\alpha\in\(0,1\)$ and $\beta\in\(0,\beta_\alpha\)$, there exist sequences $\(\gamma_k\)_{k\in\N}$ in $\(0,1\)$ and $\(\psi_k\)_{k\in\N}$ in $C^2\(M\)\backslash\left\{0\right\}$ such that $\gamma_k\to0$ as $k\to\infty$ and for each $k\in\N$, $\psi_k$ is a solution to the equation
\begin{equation}\label{Th4Step2Eq9}
\Delta_g\psi_k+h_{\alpha,\beta,\gamma_k}\psi_k=\(2^*-1\)\left|u_{\alpha,\beta,\gamma_k}\right|^{2^*-2}\psi_k\quad\text{in }M.
\end{equation}
By renormalizing, we can assume moreover that 
$$\psi_k\in\S_{H^1\(M\)}:=\big\{v\in H^1\(M\):\,\left\|v\right\|_{H^1\(M\)}=1\big\}.$$ 
Since $u_{\alpha,\beta,\gamma_k}\to u_0$ in $C^2\(M\)$ and $h_{\alpha,\beta,\gamma_k}\to h_0$ in $C^p\(M\)$ as $k\to\infty$, we then obtain that $\psi_k\to\psi_{\alpha,\beta}$ in $C^2\(M\)$ as $k\to\infty$ for some solution $\psi_{\alpha,\beta}\in C^2\(M\)\cap \S_{H^1\(M\)}$ to the equation
\begin{equation}\label{Th4Step2Eq10}
\Delta_g\psi_{\alpha,\beta}+h_0\psi_{\alpha,\beta}=\(2^*-1\)\left|u_0\right|^{2^*-2}\psi_{\alpha,\beta}\quad\text{in }M.
\end{equation}
By using \eqref{Th4Step2Eq5}, \eqref{Th4Step2Eq6}, \eqref{Th4Step2Eq9}, \eqref{Th4Step2Eq10} and the equation satisfied by $u_0$, we obtain
\begin{align}\label{Th4Step2Eq11}
&\(\Delta_g+h_0-\(2^*-1\)\left|u_0\right|^{2^*-2}\)\frac{\psi_k-\psi_{\alpha,\beta}}{\gamma_k}\nonumber\\
&\quad=\(\frac{2^*-2}{\gamma_k}\(\left|u_{\alpha,\beta,\gamma_k}\right|^{2^*-2}-\left|u_0\right|^{2^*-2}\)+\frac{u_0\Delta_gw_{\alpha,\beta}-w_{\alpha,\beta}\Delta_gu_0}{u_0 u_{\alpha,\beta,\gamma_k}}\)\psi_k\nonumber\allowdisplaybreaks\\
&\quad=\(\Delta_gw_{\alpha,\beta}+h_0w_{\alpha,\beta}-\big(1-\(2^*-2\)^2\big)\left|u_0\right|^{2^*-2}w_{\alpha,\beta}\)u_0^{-1}\psi_{\alpha,\beta}+\smallo\(1\)
\end{align}
as $k\to\infty$, uniformly in $M$. It then follows from \eqref{Th4Step2Eq10} and \eqref{Th4Step2Eq11} that
\begin{equation}\label{Th4Step2Eq12}
\int_M\(\Delta_gw_{\alpha,\beta}+h_0w_{\alpha,\beta}-\big(1-\(2^*-2\)^2\big)\left|u_0\right|^{2^*-2}w_{\alpha,\beta}\)u_0^{-1}\psi_{\alpha,\beta}^2\,dv_g=0.
\end{equation}
On the other hand, by integrating by parts and using \eqref{Th4Step2Eq8}, we obtain 
\begin{align}\label{Th4Step2Eq13}
&\int_M u_0^{-1}\psi_{\alpha,\beta}^2\,\Delta_gw_{\alpha,\beta}\,dv_g=\int_M w_{\alpha,\beta}\,\Delta_g\(u_0^{-1}\psi_{\alpha,\beta}^2\)dv_g\nonumber\allowdisplaybreaks\\
&\qquad=\int_M w_{\alpha,\beta}\(\psi_{\alpha,\beta}^2\divergence_g\(u_0^{-2}\nabla u_0\)-2\,u_0^{-1}\divergence_g\(\psi_{\alpha,\beta}\nabla \psi_{\alpha,\beta}\)\)dv_g\nonumber\\
&\qquad\qquad+2\int_M u_0^{-2}w_{\alpha,\beta}\<\nabla u_0,\nabla\(\psi_{\alpha,\beta}^2\)\>_gdv_g\nonumber\allowdisplaybreaks\\
&\qquad=\int_M u_0^{-1}w_{\alpha,\beta}\Big(u_0^{-1}\psi_{\alpha,\beta}^2\,\Delta_gu_0+2\psi_{\alpha,\beta}\Delta_g\psi_{\alpha,\beta}+2u_0^{-2}\left|\nabla u_0\right|^2_g\psi_{\alpha,\beta}^2\nonumber\\
&\qquad\qquad-2\left|\nabla\psi_{\alpha,\beta}\right|^2_g\Big)dv_g-2\int_M\<\nabla w_{\alpha,\beta},\nabla u_0\>_gu_0^{-2}\psi_{\alpha,\beta}^2\,dv_g.
\end{align}
It follows from \eqref{Th4Step2Eq12}, \eqref{Th4Step2Eq13} and the equation satisfied by $u_0$ that
\begin{multline*}
\int_M u_0^{-1}w_{\alpha,\beta}\Big(\(2^*-2\)^2\left|u_0\right|^{2^*-2}\psi_{\alpha,\beta}^2+2\psi_{\alpha,\beta}\Delta_g\psi_{\alpha,\beta}+2u_0^{-2}\left|\nabla u_0\right|^2_g\psi_{\alpha,\beta}^2\\
-2\left|\nabla\psi_{\alpha,\beta}\right|^2_g\Big)dv_g-2\int_M\<\nabla w_{\alpha,\beta},\nabla u_0\>_gu_0^{-2}\psi_{\alpha,\beta}^2\,dv_g=0,
\end{multline*}
which gives, by definition of $w_{\alpha, \beta}$,
\begin{multline}\label{Th4Step2Eq14}
\int_M u_0^{-1}\eta_{\alpha,\beta}\(v_\beta\)\Big(\(2^*-2\)^2\left|u_0\right|^{2^*-2}\psi_{\alpha,\beta}^2+2\psi_{\alpha,\beta}\Delta_g\psi_{\alpha,\beta}+2u_0^{-2}\left|\nabla u_0\right|^2_g\psi_{\alpha,\beta}^2\\
-2\left|\nabla\psi_{\alpha,\beta}\right|^2_g\Big)dv_g-2\int_M\eta'_{\alpha,\beta}\(v_\beta\)\<\nabla v_\beta,\nabla u_0\>_gu_0^{-2}\psi_{\alpha,\beta}^2\,dv_g=0.
\end{multline}
Since $\psi_{\alpha,\beta}\in\S_{H^1\(M\)}$ and $\psi_{\alpha,\beta}$ is a solution to \eqref{Th4Step2Eq10} for all $\alpha\in\(0,1\)$ and $\beta\in\(0,\beta_\alpha\)$, we obtain that up to a subsequence, $\psi_{\alpha,\beta}\to\psi_\alpha$ in $C^2\(M\)$ as $\beta\to0$ and $\psi_\alpha\to\psi_0$ in $C^2\(M\)$ as $\alpha\to0$ for some solutions $\(\psi_\alpha\)_{\alpha\in\[0,1\)}$ in $C^2\(M\)\cap\S_{H^1\(M\)}$ to the equation 
\begin{equation}\label{Th4Step2Eq15}
\Delta_g\psi_\alpha+h_0\psi_\alpha=\(2^*-1\)\left|u_0\right|^{2^*-2}\psi_\alpha\quad\text{in }M.
\end{equation}
On the other hand, by using \eqref{Th4Step2Eq3} and \eqref{Th4Step2Eq4} and since $v_\beta\to u_0$ in $C^2\(M\)$ as $\beta\to0$, we obtain
\begin{equation}\label{Th4Step2Eq16}
\left|\eta_{\alpha,\beta}\(v_\beta\(x\)\)\right|+\left|\eta'_{\alpha,\beta}\(v_\beta\(x\)\)\nabla v_\beta\right|\le C
\end{equation}
for some constant $C>0$ independent of $\alpha$ and $\beta$ and
\begin{equation}\label{Th4Step2Eq17}
\lim_{\beta\to0}\eta_{\alpha,\beta}\(v_\beta\(x\)\)=\eta_\alpha\(u_0\(x\)\):=\left\{\begin{aligned}&u_0\(x\)+\alpha&&\text{if }u_0\(x\)<-\alpha\\&0&&\text{if }-\alpha\le u_0\(x\)\le\alpha\\&u_0\(x\)-\alpha&&\text{if }u_0\(x\)>\alpha\end{aligned}\right.
\end{equation}
and
\begin{equation}\label{Th4Step2Eq18}
\lim_{\beta\to0}\eta'_{\alpha,\beta}\(v_\beta\(x\)\)\nabla v_\beta\(x\)=\left\{\begin{aligned}&0&&\text{if }\left|u_0\(x\)\right|\le\alpha\\&\nabla u_0\(x\)&&\text{if }\left|u_0\(x\)\right|>\alpha\end{aligned}\right.
\end{equation}
for all points $x\in M$. By dominated convergence, it follows from \eqref{Th4Step2Eq8} and \eqref{Th4Step2Eq14}--\eqref{Th4Step2Eq18} that
\begin{multline}\label{Th4Step2Eq19}
2\alpha\int_{M_\alpha}\left|\nabla u_0\right|^2_g\left|u_0\right|^{-3}\psi_\alpha^2\,dv_g\\
=\int_M u_0^{-1}\eta_\alpha\(u_0\)\(\(2^*-2\)^2\left|u_0\right|^{2^*-2}\psi_\alpha^2+2\psi_\alpha\Delta_g\psi_\alpha-2\left|\nabla\psi_\alpha\right|^2_g\)dv_g,
\end{multline}
where 
$$M_\alpha:=\left\{x\in M:\,\left|u_0\(x\)\right|>\alpha\right\}.$$ 
By passing to the limit as $\alpha\to0$ and using again dominated convergence, it then follows from \eqref{Th4Step2Eq19} that 
\begin{multline}\label{Th4Step2Eq20}
\lim_{\alpha\to0}\(\alpha\int_{M_\alpha}\left|\nabla u_0\right|^2_g\left|u_0\right|^{-3}\psi_\alpha^2\,dv_g\)=\frac{1}{2}\int_M\Big(\(2^*-2\)^2\left|u_0\right|^{2^*-2}\psi_0^2+2\psi_0\Delta_g\psi_0\\
-2\left|\nabla\psi_0\right|^2_g\Big)dv_g=\frac{1}{2}\(2^*-2\)^2\int_M\left|u_0\right|^{2^*-2}\psi_0^2dv_g.
\end{multline}
We claim that for small $\alpha>0$, there exists a constant $C>0$ independent of $\alpha$ such that 
\begin{equation}\label{Th4Step2Eq21}
\left|\psi_\alpha\right|\le C\left|u_0\right|\quad\text{in }M.
\end{equation}
Remark that by dominated convergence, it follows from \eqref{Th4Step2Eq21} that
\begin{equation}\label{Th4Step2Eq22}
\lim_{\alpha\to0}\(\alpha\int_{M_\alpha}\left|\nabla u_0\right|^2_g\left|u_0\right|^{-3}\psi_\alpha^2\,dv_g\)=0.
\end{equation}
To prove \eqref{Th4Step2Eq21}, we begin with showing that $\psi_0=0$ on $\Gamma_{\sigma_1}\cup\Gamma_{\sigma_2}$. Assume by contradiction that there exists a point $x_0\in \Gamma_{\sigma_1}\cup\Gamma_{\sigma_2}$ such that $\psi_0\(x_0\)\ne0$. By using Proposition~\ref{Pr1} together with the continuity of $\psi_0$, we can then assume that $\Gamma_{\sigma_1}\cup\Gamma_{\sigma_2}$ is smooth at $x_0$. Since $\psi_\alpha\to\psi_0$ in $C^0\(M\)$ as $\alpha\to0$, we then obtain that there exist constants $C_1,r_1>0$ such that $\Gamma_{\sigma_1}\cup\Gamma_{\sigma_2}$ is smooth in $B_g\(x_0,2r_1\)$ and
\begin{equation}\label{Th4Step2Eq23}
\psi_\alpha^2\ge C_1\quad\text{in }B_g\(x_0,r_1\),
\end{equation}
where $B_g\(x_0,r_1\)$ is the ball of center $x_0$ and radius $r_1$ with respect to the metric $g$. On the other hand, since $\Gamma_{\sigma_1}\cup\Gamma_{\sigma_2}$ is smooth in $B_g\(x_0,2r_1\)$ and $\Gamma_{\sigma_1}\cup\Gamma_{\sigma_2}=\left\{x\in M:\,u_0\(x\)=0\right\}$, by using Hopf's lemma, we obtain that there exists a constant $C_2>0$ such that
\begin{equation}\label{Th4Step2Eq24}
\left|\nabla u_0\right|^2\ge\frac{1}{C_2}\quad\text{and}\quad\frac{1}{C_2}\le\frac{\left|u_0\(x\)\right|}{d_g\(x,\Gamma_{\sigma_1}\cup\Gamma_{\sigma_2}\)}\le C_2\quad\text{in }B_g\(x_0,r_1\).
\end{equation}
It follows from \eqref{Th4Step2Eq23} and \eqref{Th4Step2Eq24} that
\begin{align}\label{Th4Step2Eq25}
&\int_{M_\alpha}\left|\nabla u_0\right|^2_g\left|u_0\right|^{-3}\psi_\alpha^2\,dv_g\ge\int_{B_g\(x_0,r_1\)\cap M_\alpha\backslash M_{2\alpha}}\left|\nabla u_0\right|^2_g\left|u_0\right|^{-3}\psi_\alpha^2\,dv_g\nonumber\allowdisplaybreaks\\
&\ge\frac{C_1}{8\,C_2\,\alpha^3}\Vol_g\big(B_g\(x_0,r_1\)\cap M_\alpha\backslash M_{2\alpha}\big),\nonumber\allowdisplaybreaks\\
&\ge\frac{C_1}{8\,C_2\,\alpha^3}\Vol_g\(\left\{x\in B_g\(x_0,r_1\):\,\frac{\alpha}{C_2}\le d_g\(x,\Gamma_{\sigma_1}\cup\Gamma_{\sigma_2}\)\le2\,C_2\,\alpha\right\}\),
\end{align}
where $\Vol_g$ is the volume with respect to the metric $g$. Since $\Gamma_{\sigma_1}\cup\Gamma_{\sigma_2}$ is smooth in $B_g\(x_0,2r_1\)$, we obtain that there exists a constant $C_3>0$ such that 
\begin{equation}\label{Th4Step2Eq26}
\Vol_g\(\left\{x\in B_g\(x_0,r_1\):\,\frac{\alpha}{C_2}\le d_g\(x,\Gamma_{\sigma_1}\cup\Gamma_{\sigma_2}\)\le2C_2\alpha\right\}\)\ge C_3\alpha.
\end{equation}
We then obtain a contradiction by putting together \eqref{Th4Step2Eq20}, \eqref{Th4Step2Eq25} and \eqref{Th4Step2Eq26}. This proves that $\psi_0=0$ on $\Gamma_{\sigma_1}\cup\Gamma_{\sigma_2}$. Assume now that \eqref{Th4Step2Eq21} does not hold true. Then there exist sequences $\(x_k\)_{k\in\N}$ in $M$ and $\(\alpha_k\)_{k\in\N}$ in $\(0,1\)$ such that
\begin{equation}\label{Th4Step2Eq27}
\lim_{k\to\infty}\alpha_k=0\quad\text{and}\quad\left|\psi_{\alpha_k}\(x_k\)\right|>k\left|u_0\(x_k\)\right|\quad\forall k\in\N.
\end{equation}
Since $\psi_{\alpha_k}=0$ on $\Gamma_{\sigma_1}\cup\Gamma_{\sigma_2}$, by symmetry and passing to a subsequence, we can assume that $x_k\in\Omega$ and $x_k\to x_0$ as $k\to\infty$ for some point $x_0\in \overline\Omega$. Since $\psi_\alpha\to\psi_0$ in $C^0\(M\)$ as $\alpha\to0$, it follows from \eqref{Th4Step2Eq27} that $u_0\(x_0\)=0$, i.e. $\xi_0\in\partial\Omega$. We first consider the case where $\partial\Omega$ is smooth at the point $x_0$. In this case, since $u_0=\psi_{\alpha_k}=0$ on $\partial\Omega$, Taylor expansions give
\begin{equation}\label{Th4Step2Eq28}
u_0\(x_k\)=-\partial_\nu u_0\(x_k^*\)d_g\(x_k,\partial\Omega\)+\bigO\big(d_g\(x_k,\partial\Omega\)^2\big)
\end{equation}
and 
\begin{equation}\label{Th4Step2Eq29}
\psi_{\alpha_k}\(x_k\)=-\partial_\nu\psi_{\alpha_k}\(x_k^*\)d_g\(x_k,\partial\Omega\)+\bigO\big(d_g\(x_k,\partial\Omega\)^2\big)
\end{equation}
as $k\to\infty$, where $x_k^*$ is the closest point to $x_k$ on $\partial\Omega$ and $\nu$ is the outward unit normal vector induced by the metric $g$ on $\partial\Omega$. Moreover, still assuming that $\partial\Omega$ is smooth at the point  $x_0$, since $u_0>0$ in $\Omega$, $x_k^*\in\partial\Omega$ and $x_k^*\to x_0$ as $k\to\infty$, Hopf's lemma gives that there exists a constant $C>0$ such that 
\begin{equation}\label{Th4Step2Eq30}
-\partial_\nu u_0\(x_k^*\)\ge C\quad\forall k\in\N.
\end{equation}
It follows from \eqref{Th4Step2Eq27}--\eqref{Th4Step2Eq30} that $\partial_\nu\psi_{\alpha_k}\(x_k^*\)\to-\infty$ as $k\to\infty$, which is a contradiction since $\psi_\alpha\to\psi_0$ in $C^1\(M\)$ as $\alpha\to0$. This proves that $\partial\Omega$ is not smooth at the point $x_0$. Therefore, according to Proposition~\ref{Pr1}, we are in the case where $\Gamma_{\sigma_1}$ and $\Gamma_{\sigma_2}$ intersect at the point $x_0$. In this case, by identifying $T_{x_0}M$ with $\R^n$, we obtain that 
$$\left\{\begin{aligned}&\exp_{x_0}^{-1}\(\Gamma_{\sigma_1}\)\cap B\(0,r\)=\left\{\(y_1,\dotsc,y_n\)\in B\(0,r\):\,y_1=0\right\}\\&\exp_{x_0}^{-1}\(\Gamma_{\sigma_2}\)\cap B\(0,r\)=\left\{\(y_1,\dotsc,y_n\)\in B\(0,r\):\,y_2=0\right\}\\&\exp_{x_0}^{-1}\(\Omega\)\cap B\(0,r\)=\left\{\(y_1,\dotsc,y_n\)\in B\(0,r\):\,y_1>0\quad\text{and}\quad y_2>0\right\}\end{aligned}\right.$$
for small $r>0$. Since $u_0>0$ in $\Omega$ and $u_0$ satisfies \eqref{IntroEq1} and \eqref{IntroEq7}, by applying a result by Bers \cite{Be}*{Theorem~1}, we obtain that 
\begin{equation}\label{Th4Step2Eq31}
u_0\(\exp_{x_0}\(y\)\)=\lambda y_1y_2+\smallo\big(\left|y\right|^2\big)
\end{equation}
as $y\to0$, for some number $\lambda>0$. Moreover, since $u_0=\psi_{\alpha_k}=0$ on $\Gamma_{\sigma_1}\cup\Gamma_{\sigma_2}$ and $u_0$ satisfies \eqref{Th4Step2Eq31}, Taylor expansions give 
\begin{equation}\label{Th4Step2Eq32}
u_0\(\exp_{x_0}\(y\)\)=\lambda y_1y_2+\smallo\big(\left|y_1y_2\right|\big)
\end{equation}
as $y\to0$ and
\begin{equation}\label{Th4Step2Eq33}
\left|\psi_{\alpha_k}\(\exp_{x_0}\(y\)\)\right|\le\left\|\partial_{y_1}\partial_{y_2}\(\psi_{\alpha_k}\circ\exp_{x_0}\)\right\|_{C^0\(B\(0,r\)\)}\left|y_1y_2\right|
\end{equation}
for all $k\in\N$ and $y\in B\(0,r\)$. Since $x_k\in\Omega$ and $\lambda>0$, it follows from \eqref{Th4Step2Eq27} and \eqref{Th4Step2Eq31}--\eqref{Th4Step2Eq33} that $\left\|\partial_{y_1}\partial_{y_2}\(\psi_{\alpha_k}\circ\exp_{x_0}\)\right\|_{C^0\(B\(0,r\)\)}\to\infty$ as $k\to\infty$, which is a contradiction since $\psi_\alpha\to\psi_0$ in $C^2\(M\)$ as $\alpha\to0$.. This proves our claim, namely that \eqref{Th4Step2Eq21} and \eqref{Th4Step2Eq22} hold true. It then follows from \eqref{Th4Step2Eq20} and \eqref{Th4Step2Eq22} that 
$$\int_M\left|u_0\right|^{2^*-2}\psi_0^2\,dv_g=0,$$
which gives $u_0\,\psi_0=0$ in $M$. Since $u_0$ is nontrivial and $M$ is connected, by unique continuation, we obtain that $u_0\ne0$ in a dense subset of $M$. By continuity of $\psi_0$, we then obtain that $\psi_0=0$ in $M$. This is a contradiction since $\psi_0\in \S_{H^1\(M\)}$. This proves our claim, namely that there exist $\alpha\in\(0,1\)$ and $\beta\in\(0,\beta_\alpha\)$ such that for small $\gamma>0$, $u_{\alpha,\beta,\gamma}$ is a non-degenerate solution to \eqref{Th4Step2Eq7}. By letting $h_\gamma:=h_{\alpha,\beta,\gamma}$ and $u_\gamma:=u_{\alpha,\beta,\gamma}$, this ends the proof of Step~2.
\endproof

\begin{step}\label{Th4Step3}
For each $\gamma$, there exist families of functions $\(h_{\gamma,\delta}\)_{0<\delta\ll1}$ in $C^2\(M\)\cap C^p\(M\)$ and $\(u_{\gamma,\delta}\)_{0<\delta\ll1}$ in $C^2\(M\)$ such that $h_{\gamma,\delta}\to h_\gamma$ in $C^p\(M\)$ and $u_{\gamma,\delta}\to u_\gamma$ in $C^2\(M\)$ as $\delta\to0$ and for each $\delta$, $u_{\gamma,\delta}$ and $h_{\gamma,\delta}$ are symmetric with respect to $\xi_0$, $u_{\gamma,\delta}$ is a nondegenerate solution to the equation
\begin{equation}\label{Th4Step3Eq1}
\Delta_gu_{\gamma,\delta}+h_{\gamma,\delta}u_{\gamma,\delta}=\left|u_{\gamma,\delta}\right|^{2^*-2}u_{\gamma,\delta}\quad\text{in }M
\end{equation}
and
\begin{equation}\label{Th4Step3Eq2}
\left\{\begin{aligned}&m_{\tilde{h}_{\gamma,\delta},g}\(\xi_0\)\ne0&&\text{if }n=3\\&h_{\gamma,\delta}\(\xi_0\)\ne c_n\Scal_g\(\xi_0\)&&\text{if }n\in\left\{4,5\right\},\end{aligned}\right.
\end{equation}
where $\tilde{h}_{\gamma,\delta}:=h_{\gamma,\delta}-5u_{\gamma,\delta}^4$.
\end{step}

\proof[Proof of Step~\ref{Th4Step3}]
By using a straightforward regularization argument, we may assume that $p\ge2$. We first consider the case where $n\in\left\{4,5\right\}$. In this case, we define $h_{\gamma,\delta}:=h_\gamma+\delta$, so that $h_{\gamma,\delta}\(\xi_0\)\ne c_n\Scal_g\(\xi_0\)$ for small $\delta>0$. Since $h_\gamma\in C^p\(M\)$ and $h_\gamma$ is $2$-symmetric with respect to $\xi_0$, we obtain that $h_{\gamma,\delta}\in C^p\(M\)$, $h_{\gamma,\delta}$ is $2$-symmetric with respect to $\xi_0$ and $h_{\gamma,\delta}\to h_\gamma$ in $C^p\(M\)$ as $\delta\to0$. Since $u_\gamma$ is a nondegenerate solution to \eqref{Th4Step2Eq0}, we obtain that for small $\delta>0$, there exists a solution $u_{\gamma,\delta}\in C^2\(M\)$ to \eqref{Th4Step3Eq1} such that $u_{\gamma,\delta}\to u_\gamma$ in $C^2\(M\)$ as $\delta\to0$. Moreover, since $h_{\gamma,\delta}$ is $2$-symmetric with respect to $\xi_0$, $u_\gamma$ is symmetric with respect to $\xi_0$ and $u_{\gamma}\circ\sigma_1=u_{\gamma}\circ\sigma_2=-u_{\gamma}$ in $M$, we can assume that $u_{\gamma,\delta}$ is symmetric with respect to $\xi_0$ and $u_{\gamma,\delta}\circ\sigma_1=u_{\gamma,\delta}\circ\sigma_2=-u_{\gamma,\delta}$ in $M$. In particular, since $\xi_0\in\Gamma_{\sigma_1}\cap\Gamma_{\sigma_2}$, we obtain that $u_{\gamma,\delta}\(\xi_0\)=0$. We now consider the case where $n=3$. In this case, we define 
$$u_{\gamma,\delta}:=\(1+\delta\)u_\gamma\quad\text{and}\quad h_{\gamma,\delta}:=h_\gamma+\big(\(1+\delta\)^4-1\big)u_\gamma^4,$$
so that \eqref{Th4Step3Eq1} holds true. Since $h_\gamma\in C^p\(M\)$, standard elliptic regularity theory gives that $u_\gamma\in C^{p+1}\(M\)$, and so $h_{\gamma,\delta}\in C^p\(M\)$, $u_{\gamma,\delta}\in C^2\(M\)$, $h_{\gamma,\delta}\to h_\gamma$ in $C^p\(M\)$ and $u_{\gamma,\delta}\to u_\gamma$ in $C^2\(M\)$ as $\delta\to0$. Moreover, since $h_\gamma$ and $u_\gamma$ are symmetric with respect to $\xi_0$ and $u_{\gamma}\(\xi_0\)=0$ in $M$, we obtain that $h_{\gamma,\delta}$ and $u_{\gamma,\delta}$ are symmetric with respect to $\xi_0$ and $u_{\gamma,\delta}\(\xi_0\)=0$ in $M$. Now, for small $\delta>0$, by remarking that
$$\(\Delta_g+h_\gamma-5u_\gamma^4\)\big(G_{g,\tilde{h}_{\gamma,\delta}}\(\xi_0,\cdot\)-G_{g,\tilde{h}_{\gamma}}\(\xi_0,\cdot\)\big)=4\big(\(1+\delta\)^4-1\big)u_\gamma^4G_{g,h_{\gamma,\delta}}\(\xi_0,\cdot\)$$
in $M$, where $\tilde{h}_\gamma:=h_\gamma-5u_\gamma^4$ and $\tilde{h}_{\gamma,\delta}:=h_{\gamma,\delta}-5u_{\gamma,\delta}^4$, we obtain
\begin{align}\label{Th4Step3Eq3}
m_{g,\tilde{h}_{\gamma,\delta}}\(\xi_0\)-m_{g,\tilde{h}_\gamma}\(\xi_0\)&=4\big(\(1+\delta\)^4-1\big)\int_Mu_\gamma^4G_{g,\tilde{h}_\gamma}\(\xi_0,\cdot\)G_{g,\tilde{h}_{\gamma,\delta}}\(\xi_0,\cdot\)dv_g\nonumber\\
&=16\delta\int_Mu_\gamma^4G_{g,\tilde{h}_\gamma}\(\xi_0,\cdot\)^2dv_g+\smallo\(\delta\)
\end{align}
as $\delta\to0$. Since $u_\gamma$ is nontrivial and $M$ is connected, by unique continuation, we obtain that $u_\gamma\ne0$ and $G_{g,\tilde{h}_\gamma}\(\xi_0,\cdot\)\ne0$ in a dense subset of $M$. It then follows from \eqref{Th4Step3Eq3} that $m_{g,\tilde{h}_{\gamma,\delta}}\(\xi_0\)\ne0$ for small $\delta>0$. This ends the proof of Step~\ref{Th4Step3}.
\endproof

\proof[End of proof of Theorem~\ref{Th4}]
To prove Theorem~\ref{Th4}, it remains to apply Theorem~\ref{Th3}. Observe that
\begin{equation}\label{Th4Eq}
\big[\big(\Delta_g+\tilde{h}_{\gamma,\delta}\big)^{-1}\(hu_{\gamma,\delta}\)\big]\(\xi_0\)=\int_Mhu_{\gamma,\delta}G_{g,\tilde{h}_{\gamma,\delta}}\(\xi_0,\cdot\)dv_g
\end{equation}
for all functions $h\in C^{0,\vartheta}\(M\)$. Since $u_{\gamma,\delta}$ is nontrivial and $M$ is connected, by unique continuation, we obtain that $u_{\gamma,\delta}\ne0$ and $G_{g,\tilde{h}_{\gamma,\delta}}\(\xi_0,\cdot\)\ne0$ in a dense subset of $M$. In particular, we can choose a function $\hat{h}_{\gamma,\delta}\in C^p\(M\)$ such that $\hat{h}_{\gamma,\delta}=0$ in some neighborhood of the point $\xi_0$, $\hat{h}_{\gamma,\delta}u_{\gamma,\delta}G_{g,\tilde{h}_{\gamma,\delta}}\(\xi_0,\cdot\)\ge0$ at all points in $M$ and $\hat{h}_{\gamma,\delta}u_{\gamma,\delta}G_{g,\tilde{h}_{\gamma,\delta}}\(\xi_0,\cdot\)>0$ at some point in $M$. It then follows from \eqref{Th4Step3Eq2} and \eqref{Th4Eq} that \eqref{Th3Eq1} holds true with $h=\pm\hat{h}_{\gamma,\delta}$. By applying Theorem~\ref{Th3}, we then obtain that there exists a family of solutions $\(u_{\gamma,\delta,\eps}\)_{0<\eps\ll1}$ to \eqref{IntroEq2} with $h_\eps=h_0+\eps\hat{h}_{\gamma,\delta}$ or $h_\eps=h_0-\eps\hat{h}_{\gamma,\delta}$ which satisfies \eqref{IntroEq3} with $u_0=u_{\gamma,\delta}$ and $\xi_\eps\to\xi_0$ as $\eps\to0$. We can now conclude by applying a straightforward diagonal argument. This ends the proof of Theorem~\ref{Th4}.
\endproof

We can now prove Corollaries~\ref{Cor3} and~\ref{Cor4}.

\proof[Proof of Corollary~\ref{Cor3}]
Since the operator $\Delta_g+h_0$ is coercive in $H^1_{0,\sigma_0}\(\Omega\)$, by proceeding in a similar way as in the article by Aubin \cite{A}*{Theorem~4}, we obtain that there exists a positive solution $u_0$ to \eqref{IntroEq7} provided we can show that
\begin{equation}\label{Th4Step1Eq2}
\min_{u\in H^1_{0,\sigma_0}\(\Omega\)\backslash\left\{0\right\}}I_{g,h_0}\(u\)<K_n^{-2},
\end{equation}
where $K_n$ is the best constant for the embedding $D^{1,2}\(\R^n\)\hookrightarrow  L^{2^*}\(\R^n\)$ and
\begin{equation}\label{Th4Step1Eq3}
I_{g,h_0}\(u\):=\frac{\displaystyle\int_\Omega\(\left|\nabla u\right|_g^2+h_0u^2\)dv_g}{\displaystyle\(\int_\Omega\left|u\right|^{2^*}dv_g\)^{2/2^*}}\,.
\end{equation}
To prove \eqref{Th4Step1Eq3}, we let $\(B_\mu\)_{\mu>0}$ be as in \eqref{Th3Eq4} with $\xi_1$ instead of $\xi_0$. Since $\xi_1\in\Omega\cap\Gamma_{\sigma_0}$, by letting the outer radius $r_0$ of the support of the function $\chi$ be small enough, we obtain that $B_\mu\in H^1_{0,\sigma_0}\(\Omega\)$. Since $h_0\(\xi_1\)<c_n\Scal_g\(\xi_1\)$, standard computations as in \cite{A} then give $I_{g,h_0}\(B_\mu\)<K_n^{-2}$ for small $\mu>0$, thus proving \eqref{Th4Step1Eq2}. We can then conclude by applying Theorem~\ref{Th4}.
\endproof 

\proof[Proof of Corollary~\ref{Cor4}]
In this case, by using \eqref{IntroEq8} and proceeding in a similar way as in the article by Hebey \cite{He1}*{Theorem~1} (see also the article by Hebey and Vaugon \cite{HV}), we obtain the existence of a $G$-invariant, positive solution $u_0$ to \eqref{IntroEq7}. We can then conclude by applying Theorem~\ref{Th4}.
\endproof

\end{document}